\documentclass[11pt,reqno]{amsart}
\usepackage{ytableau}
\usepackage{amsthm}
\usepackage{amssymb}
\usepackage{latexsym}
\usepackage{multicol}
\usepackage{etoolbox} 
\usepackage{verbatim,enumerate}
\usepackage{accents}
\usepackage[utf8]{inputenc}
\usepackage{hyperref}
\usepackage{amsmath, amscd}
\usepackage{soul}
\usepackage{tikz} 
\usepackage{mathtools}
\usetikzlibrary{matrix,arrows,decorations.pathmorphing,arrows.meta,calc}
\newlength\cellsize 

\savebox2{%
\begin{picture}(10,10)
\put(0,0){\line(1,0){10}}
\put(0,0){\line(0,1){10}}
\put(10,0){\line(0,1){10}}
\put(0,10){\line(1,0){10}}
\end{picture}}

\newcommand\cellify[1]{\def\thearg{#1}\def\nothing{}%
\ifx\thearg\nothing\vrule width0pt height\cellsize depth0pt%
  \else\hbox to 0pt{\usebox2\hss}\fi%
  \vbox to 10\unitlength{\vss\hbox to 10\unitlength{\hss$#1$\hss}\vss}}

\newlength{\bibitemsep}
\newlength{\bibparskip}
\let\oldthebibliography\thebibliography
\renewcommand\thebibliography[1]{%
  \oldthebibliography{#1}%
  \setlength{\parskip}{\bibitemsep}%
  \setlength{\itemsep}{\bibparskip}%
}

\newcommand\tableau[1]{\vtop{\let\\=\cr
\setlength\baselineskip{-10000pt}
\setlength\lineskiplimit{10000pt}
\setlength\lineskip{0pt}
\halign{&\cellify{##}\cr#1\crcr}}}

\savebox7{
\begin{picture}(10,10)
  \put(5,5){\circle{9.4}}
\end{picture}}

\newcommand{\cirfy}[1]{\def\thearg{#1}\def\nothing{}%
\ifx\thearg\nothing\vrule width0pt height\cellsize depth0pt%
  \else\hbox to 0pt{\usebox7\hss}\fi%
  \vbox to 10\unitlength{\vss\hbox to 10\unitlength{\hss$#1$\hss}\vss}}

\newcommand\cirtab[1]{\vtop{\let\\=\cr
\setlength\baselineskip{-10000pt}
\setlength\lineskiplimit{10000pt}
\setlength\lineskip{0pt}
\halign{&\cirfy{##}\cr#1\crcr}}}

\theoremstyle{plain}
\newtheorem*{prop}{Proposition}
\newtheorem{thm}{Theorem}
\newtheorem*{lem}{Lemma}

\theoremstyle{definition}
\newtheorem*{example}{Example}
\newtheorem*{defn}{Definition}

\newtheorem*{rem}{Remark}

\theoremstyle{remark}

\newcommand{\lie}[1]{\mathfrak{#1}}

\newcommand\bc{\mathbb C}

\newcommand\bz{\mathbb Z}

\def\gr{\operatorname{gr}}

\advance\textwidth by 1.2in  \advance\oddsidemargin by -.6in \advance\evensidemargin by -.6in

\newcounter{cnt}
 \makeatletter
\def\mydggeometry{\makeatletter\dg@YGRID=1\dg@XGRID=20\unitlength=0.003pt\makeatother}
\makeatother \theoremstyle{remark}

\numberwithin{equation}{section}
\makeatletter
\def\section{\def\@secnumfont{\mdseries}\@startsection{section}{1}%
  \z@{.7\linespacing\@plus\linespacing}{.5\linespacing}%
  {\normalfont\scshape\centering}}
\def\subsection{\def\@secnumfont{\bfseries}\@startsection{subsection}{2}%
  {\parindent}{.5\linespacing\@plus.7\linespacing}{-.5em}%
  {\normalfont\bfseries}}
\def\subsubsection{\def\@subsecnumfont{\bfseries}\@startsection{subsubsection}{3}%
  {\parindent}{.5\linespacing\@plus.7\linespacing}{-.5em}%
  {\normalfont\bfseries}}
  \makeatother

\begin{document}

\title[Pieri formulas, higher level Demazure crystals and numerical multiplicities]{Pieri formulas, higher level Demazure crystals and numerical multiplicities of excellent filtrations}
\author{Deniz Kus}
\address{University of Bochum, Faculty of Mathematics, Universit{\"a}tsstr. 150, Bochum 44801, 
Germany}
\email{deniz.kus@rub.de}
\thanks{}
\author{Valentin Rappel}
\address{University of Cologne, Department of Mathematics and Computer Science, Division of Mathematics, Weyertal 86-90, K{\"o}ln 50931, Germany}
\email{math@vrappel.de}

\thanks{D.K. and V.R. were partially funded by the Deutsche Forschungsgemeinschaft (DFG, German Research Foundation)-- grant 446246717.}
\subjclass[2020]{}

\begin{abstract}
The classical Pieri formula gives a multiplicity free expansion of an irreducible module with a fundamental one for the complex general linear group. In this article we replace the tensor product by the fusion product and prove an analogue Pieri formula for higher level Demazure modules for the affine Lie algebra $\widehat{\mathfrak{sl}}_{n+1}$. To be more precise, we show that the fusion product of an arbitrary stable Demazure module with a fundamental module admits a multiplicity free excellent filtration and the successive quotients are described explicitly. As a consequence, we derive recurrence relations for the generating series encoding the numerical multiplicities in Demazure flags of level one Demazure modules.
\end{abstract}
\maketitle

\section{Introduction}
The Pieri rule gives a multiplicity free decomposition formula of the tensor product of two irreducible finite-dimensional representations for the complex general linear group when one of the factors has a special form. Originally stated in the context of Schubert calculus \cite{schubert79}, it has many important consequences and applications. For example, the Pieri formula gives the Giambelli formula and this in turn implies that the cohomology ring of the Grassmannian is generated by special classes of Schubert cycles.

Since then, many generalizations of this formula have been considered in various contexts. The Pieri rule for flag manifolds and Schubert polynomials has been studied in \cite{sottile96} and recently a Pieri formula for other complex algebraic groups has been investigated in \cite{vandiejen22}. Another possible direction of generalizing the Pieri formula is to replace the tensor product by the fusion product (also called graded tensor product) \cite{FL99} and study the graded multiplicity free decomposition of the fusion product of two irreducible finite-dimensional representations when one of the factors has again a special form \cite{BK22a,fourier15}. The fusion product gives a graded representation for the current algebra $\mathfrak{g}[t]$ which is the Lie algebra of polynomial maps $\mathbb{C}\rightarrow \lie g$ of a simple Lie algebra $\lie g$ (see Section~\ref{section26} for details).

Many interesting classes of representations for current algebras have been studied in the past decades due to their rich connection to combinatorics, representation theory of quantum affine algebras, number theory and mathematical physics. In this article we will focus mainly on the class of higher level affine Demazure modules and their fusion products. Demazure \cite{demazure74}
generalized the Kac-Weyl character formula to certain submodules for the  Borel subalgebra of integrable, highest weight representations of the affine Lie algebra. If the center of the affine algebra acts by a positive integer $\ell$, these modules are referred to as the Demazure modules of level $\ell$. In this article we are interested in Demazure modules which admit the structure of a $\lie g[t]$-module; these are known in the literature as \textit{stable} Demazure module and are parametrized by tuples $(\lambda,\ell)$ where $\lambda$ is a dominant integral weight for $\lie g$ and $\ell\in\mathbb{N}$ is the level. To avoid confusion with other articles we emphasize that the highest weight $\lambda$ does not need to be a multiple of $\ell$ and all prime Demazure modules are permitted.

The celebrated Demazure character formula \cite{kumar87} can be used to identify level one Demazure characters with certain specializations of non-symmetric Macdonald polynomials \cite{I03, sanderson00}. Level two Demazure modules appear as classical limits of a family of irreducible representations of the quantum affine algebra \cite{BCMo15} and their graded decompositions are given in \cite{BK20b} by using the combinatorics of convex polytopes. Moreover, there is a deep connection to Kirillov-Reshetikhin modules and consequently they satisfy certain functional relations which are known as $Q$-systems \cite{CV13,KV14}. However, the most basic questions like closed dimension formulas are, in general, unknown for higher level Demazure modules. One of the approaches in the past to better understand them was to observe tensor product factorizations \cite{CSVW16,FoL06}; however the prime higher level Demazure modules seem to be still mysterious.

Another natural approach is to consider the expansion of a product of a level $\ell$ stable Demazure module with a fundamental module into irreducibles like in the Pieri formula. However the coefficients appearing are not, in general, zero or one and can grow arbitrarily. The more natural alternative approach to this problem is to expand the tensor product into level $\ell$ stable Demazure modules again instead of irreducibles; similar questions were considered in the finite case for example in \cite{assaf19}. This leads to the notion of an \textit{excellent filtration}. 

Objects in the category of finite-dimensional $\mathbb{Z}$-graded modules for the current algebra $\lie g[t]$ admit an \textit{excellent filtration} if there exists a flag whose successive quotients are isomorphic to stable Demazure modules. Naoi proved for an affine Lie algebra associated to a simply laced simple Lie algebra \cite{Na11} that for a given $\ell\geq m\geq 1$ any stable Demazure module of level $m$ admits a filtration such that the successive quotients are isomorphic to level $\ell$ Demazure modules. The corresponding multiplicities are refered to as \textit{graded} and \textit{numerical} multiplicities respectively.

In rank one meany deep connections were observed by various authors. The numerical multiplicities are closely related to Chebyshev polynomials \cite[Corollary 1.3]{BCSV15}, several specializations of the generting series associated to the graded multiplicities specialize to Ramanujan’s fifth order mock theta functions \cite[Theorem 1.6]{BCSV15}, certain weighted versions of the generating series give Carlitz q-Fibonacci polynomials \cite[Proposition 2.5.3]{BCK18} and are limits of hypergeometric series \cite[Section 2.5.4]{BCK18} and there is a combinatorial formula using Dyck path \cite{BK19a}.

All these approaches have one common characteristics: They use the Pieri type decomposition, i.e. the decomposition of the graded character of the tensor product of a level $\ell$ Demazure module and a fundamental module into the graded characters of level $\ell$ Demazure modules; see \cite[Proposition 6.1]{BCSV15} and \cite[Section 6.3]{BCK18}. However the expansion involves polynomials in $\mathbb{Z}[q]$ as coefficients with possibly negative coefficients and computer based calculations show that they become quite complicated when the rank of the Lie algebra increases.

The more natural approach is to replace the tensor product by the fusion product and ask the following question. Does the Pieri type fusion product admit an excellent filtration by level $\ell$ Demazure modules? What are the precise successive quotients appearing in a flag? In this article we answer this question in the affirmative and give the explicit constituents of the flag (see Theorem~\ref{mainthmres1}(1)), determine generators and relations for the fusion product (see Theorem~\ref{mainthmres1}(2)) and derive recurrence relations for the generating series encoding the numerical multiplicities of a level $\ell$ flag in a level one Demazure module. 

An important ingredient in the proof of the main theorem is crystal theory; in particular the fact that a Demazure crystal can be embedded into the tensor product of Kirillov-Reshetikhin crystals \cite{LS19a,ST12} with an affine highest weight crystal (see Theorem~\ref{schuting} and Theorem~\ref{crysesti1}.

\textit{Our paper is organized as follows:} In Section~\ref{section2} we define the main notions of the article, e.g. affine Lie algebras, fusion products and prove an elementary result on the affine Weyl group orbit. In Section~\ref{section3} we recall the notion of excellent filtrations, Demazure modules and their crystal bases, KR crystals and prove a lower bound for the combinatorial Pieri formula (see Theorem~\ref{crysesti1}) using crystal theory. In Section~\ref{mainresultsection} we state our main results (see Theorem~\ref{mainthmres1}) and discuss an alternative approach in rank two. In Section~\ref{section5} we introduce a poset which we use in Section~\ref{section6} to construct a filtration whose successive quotients are homomorphic images of Demazure modules. In Section~\ref{appendix} we prove the transitivity of the poset introduced in Section~\ref{section5}.

\textit{Acknowledgement:} The first author thanks Christian Lenart for many helpful discussions.

\section{Preliminaries}\label{section2}
\subsection{}
Throughout this paper we denote by $\mathbb{C}$ the field of complex numbers and by $\mathbb{Z}$ (resp. $\mathbb{Z}_{+}$, $\mathbb{N}$) the subset of integers (resp. non-negative, positive integers). For a Lie algebra $\mathfrak{a}$, let $\mathbf{U}(\mathfrak{a})$ be the universal enveloping algebra of $\mathfrak{a}$ and let $\lie a[t^{\pm}]=\lie a\otimes \mathbb{C}[t^{\pm}]$ the loop algebra with Lie bracket
$$[x\otimes t^r, y\otimes t^s]=[x,y]\otimes t^{r+s},\ \ x,y\in\lie a,\ \ r,s\in \mathbb{Z}.$$

\subsection{} We denote by $\mathfrak{g}$ the complex special linear Lie algebra $\mathfrak{sl}_{n+1}$ of order $(n+1)$. Let $\mathfrak{h}$ be the Cartan subalgebra of trace zero diagonal matrices, $\mathfrak{b}$ the Borel subalgebra of trace zero upper triangular matrices and
$$\Pi=\{\alpha_1=\epsilon_1-\epsilon_2,\dots,\alpha_n=\epsilon_n-\epsilon_{n+1}\},$$
$$R=\{\epsilon_i-\epsilon_{j}: 1\leq i\neq j \leq n+1\},\ R^+=\{\alpha_{i,j}:=\epsilon_i-\epsilon_{j+1}: 1\leq i\leq j \leq n\},$$  be the corresponding set of simple roots, roots and positive roots respectively. Let $(\cdot,\cdot)$ be the non-degenerate bilinear form on $\lie h^{*}$ with $(\epsilon_i,\epsilon_j)=0$, $i\neq j$, induced by the restriction of the (suitably normalized) Killing form of $\lie g$ to $\lie h$. Furthermore, we denote by
$\{\varpi_1,\dots,\varpi_n\}$ the set of fundamental weights, $\varpi_j=\epsilon_1+\cdots+\epsilon_j, \ 1\leq j\leq n$, the $\mathbb{Z}$ (resp. $\mathbb{Z}_+$)-span of the simple roots by $Q$ (resp. $Q^+$) and the $\mathbb{Z}$ (resp. $\mathbb{Z}_+$)-span of the fundamental weights by $P$ (resp. $P^+$). The dominance ordering for $\alpha,\beta\in R$ is given by
$$\alpha\geq \beta :\iff \alpha-\beta\in Q^+.$$
For a positive root \(\alpha_{i,j}\) we denote by \(x_{i,j}^{-}\) (resp. \(x_{i,j}^{+}\)) the elementary matrix \(E_{j+1,i}\) (resp. $E_{i,j+1}$) which is a root vector of weight $\pm \alpha_{i,j}$ and $h_{i,j}:=E_{i,i}-E_{j+1,j+1}$ the corresponding coroot. In particular we have
\[
    \left[x^-_{j+1,k}, x^-_{i,j}\right] = x^-_{i,k}.
\]
In the rest of the paper we abbreviate
$$x^{\pm}_{\alpha_{i,j}}=x^{\pm}_{i,j},\ x^{\pm}_{i}=x^{\pm}_{i,i},\ h_{\alpha_{i,j}}:=h_{i,j},\ h_{i}=h_{i,i}.$$
We have a triangular decomposition 
$$\mathfrak{g}=\mathfrak{n}^{-}\oplus\mathfrak{h}\oplus \mathfrak{n}^{+},\ \ \mathfrak{n}^{\pm}=\bigoplus_{\alpha\in R^+} \mathbb{C} \cdot x_{\alpha}^{\pm}.$$
The unique irreducible representation of $\mathfrak{g}$ of highest weight $\lambda\in P^+$ is denoted by $V(\lambda)$.
\subsection{}\label{section13} Given $\lambda\in P^+$, $\ell\in\mathbb{N}$ and $\alpha\in R^+$ we write
\begin{equation}\label{teilmitell}\lambda(h_{\alpha})=(s_{\alpha}-1)\ell+m_{\alpha},\ 1\leq m_{\alpha}\leq \ell \end{equation}
The remainder depends on $\ell$ and a more appropriate notation would be $m_{\alpha}(\ell)$ and $s_{\alpha}(\ell)$ respectively. However, we will use the shorter notation since the divisor $\ell$ remains most of the time the same throughout the paper and we abbreviate $m_{i,j}:=m_{\alpha_{i,j}}$ and $s_{i,j}:=s_{\alpha_{i,j}}$ respectively. We introduce the sets
$$R_i^+=\left\{\alpha\in R^+: \varpi_i(h_{\alpha})=1\right\},\ \ R_{\lambda,i}^+=\left\{\alpha\in R_i^+: \lambda(h_{\alpha})\notin \mathbb{Z}\ell\right\}$$
 The height and support respectively of a positive root $\alpha_{i,j}$ is defined as $$\mathrm{ht}(\alpha_{i,j}):=j-i+1,\ \ \mathrm{supp}(\alpha_{i,j}):=\{i,i+1,\dots,j\}.$$ Furthermore, we denote by $W\cong \Sigma_{n+1}$ the Weyl group of $\lie g$ which is the group of permutations on $n+1$ letters and let $w_0\in W$ be the longest word in $W$. Note that the Weyl group orbit of a fundamental weight is given by
$$W(\varpi_i)=\left\{\epsilon_{j_1}+\cdots+\epsilon_{j_i}: 1\leq j_1< \cdots < j_i\leq n+1\right\}.$$
Note that $\varpi_i-\chi\in Q^+$ for all $\chi\in W(\varpi_i)$ and hence $\varpi_i-\chi=\sum_{j=1}^n k_j^{\chi}\alpha_j$ for suitable $k_1^{\chi},\dots,k_n^{\chi}\in\mathbb{Z}_+$.
The tensor product of two irreducible finite-dimensional $\lie g$-modules is semi-simple and the following decomposition rule is known as the classical Pieri formula 
$$V(\lambda)\otimes V(\varpi_i)=\bigoplus_{\substack{\mu\in W(\varpi_i)\\ \lambda+\mu\in P^+}}V(\lambda+\mu)$$
The aim of this article is to prove a Pieri formular for higher level Demazure modules where the tensor product is replaced by the fusion product (see Theorem~\ref{mainthmres1} for details).
\subsection{}Let $\widehat{\mathfrak{g}}$ the untwisted affine Lie algebra associated to $\mathfrak{g}$ which is realized as 
$$\widehat{\mathfrak{g}}=\mathfrak{g}\otimes\mathbb{C}[t^{\pm}]\oplus \mathbb{C}c\oplus \mathbb{C} d$$ where $c$ is required to be central and the Lie bracket is defined as
$$[x\otimes t^{r},y\otimes t^s]=[x,y]\otimes t^{r+s}+\mathrm{tr}(xy)c,\ \ [d,x\otimes t^r]=r(x\otimes t^r),\ \ \ x,y\in\mathfrak{g},\ \ r,s\in\mathbb{Z}.$$
The Cartan subalgebra and standard Borel subalgebra respectively is realized as 
$$\widehat{\mathfrak{h}}=(\mathfrak{h}\otimes 1)\oplus \mathbb{C}c\oplus \mathbb{C} d,\ \ \widehat{\mathfrak{b}}=\mathfrak{g}\otimes t\mathbb{C}[t]\oplus (\mathfrak{b}\otimes 1)\oplus \mathbb{C}c\oplus \mathbb{C} d$$
We denote by $\delta$ the non-divisible positive imaginary root and recall that the affine fundamental weights are given by 
$$\Lambda_i=\varpi_i+\Lambda_0,\ \ 1\leq i\leq n$$
where $\Lambda_0$ is the zeroth affine fundamental weight determined by
$\Lambda_0(c)=1,\  \Lambda_0(\lie h\oplus \mathbb{C}d)=0.$
The $\mathbb{Z}$-span (resp. $\mathbb{Z}_+$-span) of the affine fundamental weights is denoted by $\widehat{P}$ (resp. $\widehat{P}^+$). A simple system compatible with $\Pi$ is given by the set $\{\alpha_0:=\delta-\alpha_{1,n},\alpha_1,\dots,\alpha_n\}$. Note that the affine Weyl group $W^{\mathrm{aff}}$ which is the group generated by the simple reflections $\mathbf{s}_j, 0\leq j\leq n$, can be realized as $W^{\mathrm{aff}}=W\ltimes t_{Q}$ where a translation $t_{\mu}, \mu\in Q$ acts by 
 $$ t_{\mu}(\beta) = \beta \mod \delta,\ \ \beta\in \mathfrak{h}^*\oplus \mathbb{C}\delta,\ \ t_{\mu}(\Lambda_0)=\Lambda_0+\mu \mod\delta.$$ 
Note that an element $b_1\epsilon_1+\cdots+b_n\epsilon_n$ with integer coefficients is contained in $Q$ if and only if the sum over the coefficients $b_1+\cdots+b_n$ is divisible by $(n+1)$.
It will also be necessary to introduce the extended affine Weyl group $W^{\mathrm{ext}}= W^{\mathrm{aff}}\ltimes t_P$ which admits an alternative description as follows
$$W^{\mathrm{ext}}=\mathcal{T}\ltimes W^{\mathrm{aff}}$$
where $\mathcal{T}$ is the subgroup of $W^{\mathrm{aff}}$ stabilizing the dominant Weyl chamber; these elements correspond to automorphisms of the affine Dynkin diagram. We denote by $|\cdot |: W^{\mathrm{ext}}\rightarrow \mathbb{Z}_+$ the length function and recall that elements in $\mathcal{T}$ have length zero. The unique irreducible integrable representation of $\widehat{\mathfrak{g}}$ of highest weight $\Lambda\in \widehat{P}^+$ is denoted by $V(\Lambda)$.
\begin{lem}\label{stabwcv2}
Let $\Lambda\in \widehat{P}^+$ and $1\leq k_1<\cdots <k_i\leq n+1$. Then there exists $z\in W^{\mathrm{aff}}_{\Lambda}$ in the stabilizer of $\Lambda$ such that 
$$\Lambda+z(\epsilon_{k_1}+\cdots+\epsilon_{k_i})\in \widehat{P}^+.$$
\begin{proof}
Let $\Lambda=a_0\Lambda_0+\cdots+a_n\Lambda_n$ and note that $$\epsilon_{k_1}+\cdots+\epsilon_{k_i}=\Lambda_{k_1}+\cdots+\Lambda_{k_i}-\Lambda_{k_1-1}-\cdots-\Lambda_{k_i-1}$$
where the indices of the affine fundamental weights are understood modulo $(n+1)$. Let $j$ be the minimal index such that $a_{j}>0$ and $p\in\{1,\dots,i\}$ such that $k_{p-1}\leq j<k_p$ where we understand $k_0=0$. Note that we can exclude the case $j\geq k_i$ for the following reason. If $j\geq k_i$ we can act with $\mathbf{s}_0\cdots \mathbf{s}_{k_1-1}\in W^{\mathrm{aff}}_{\Lambda}$ and obtain an element where we have replaced $\Lambda_{k_1}-\Lambda_{k_1-1}$ by $\Lambda_{n+1}-\Lambda_{n}$. So we could continue the proof with the newly obtained element. \par

Hence we suppose $j<k_i$ and denote the maximal one among all such indices by $j_0$, i.e. $a_{j_0}>0$, $k_{p-1}\leq j_0<k_p$ and $a_{j_0+1}=\cdots=a_{k_p-1}=0$. By acting with 
$$\mathbf{s}_{j_0+1}\cdots \mathbf{s}_{k_p-1}\in W^{\mathrm{aff}}_{\Lambda}$$
we obtain an element where we have replaced $\Lambda_{k_p}-\Lambda_{k_p-1}$ by $\Lambda_{j_0+1}-\Lambda_{j_0}$. Hence we can assume without loss of generality from the beginning that $a_{k_p-1}\neq 0$. To keep the notation as simple as possible let $p=1$. Now acting with 
 $$\mathbf{s}_{t+1}\mathbf{s}_{t+2}\cdots \mathbf{s}_{k_2-1}\in W^{\mathrm{aff}}_{\Lambda}$$ 
we obtain an element where we have further replaced 
 $\Lambda_{k_2}-\Lambda_{k_2-1}$ by $\Lambda_{t+1}-\Lambda_{t}$ 
where $t=k_1$ if $a_{k_1+1}=\cdots=a_{k_2-1}=0$ or $t\in\{k_1,\dots, k_2-1\}$ is the maximal index with $a_{t}>0$. By repeating this process with the subsequent indices $3,\dots,i$ we obtain the desired dominant weight.
\end{proof}
\end{lem}
\begin{example}
Let \(\Lambda=\Lambda_0+\Lambda_1+\Lambda_6,\ n+1=8\) and \(k_1=4,k_2=5,k_3=8\). Thus 
\[\Lambda+\varepsilon_4+\varepsilon_5+\varepsilon_8=(\Lambda_0+\Lambda_1+\Lambda_6)+(\Lambda_4-\Lambda_3+\Lambda_5-\Lambda_4+\Lambda_8-\Lambda_7).\]
In accordance to the proof of Lemma~\ref{stabwcv2} we have \(j=0\) and \(p=1\), but \(a_1\neq 0\) and \(k_0=0\leq 1\leq k_1=4 \) and thus \(j_0\neq 0\). Instead \(j_0=1\) and we act by the element \(\mathbf{s}_2\mathbf{s}_3\in W^{\mathrm{aff}}_{\Lambda}\) and obtain
\[\mathbf{s}_2\mathbf{s}_3(\Lambda+\varepsilon_4+\varepsilon_5+\varepsilon_8)=(\Lambda_0+\Lambda_1+\Lambda_5)+(\Lambda_2-\Lambda_1+\Lambda_5-\Lambda_4+\Lambda_7-\Lambda_6).\]
We observe that the negative coefficient of \(\Lambda_1\) in the second term gets canceled by the first term. 
We continue with the next index \(k_2=5\) and act by the element \(\mathbf{s}_3\mathbf{s}_4\in W^{\mathrm{aff}}_{\Lambda}\) and obtain 
\[\mathbf{s}_3\mathbf{s}_4\mathbf{s}_2\mathbf{s}_3(\Lambda+\varepsilon_4+\varepsilon_5+\varepsilon_8)=(\Lambda_0+\Lambda_1+\Lambda_5)+(\Lambda_2-\Lambda_1+\Lambda_3-\Lambda_2+\Lambda_7-\Lambda_6).\]
This illustrates the case \(a_{k_1+1}=\cdots = a_{k_2-1}=0\). 
The next index to consider is \(k_3=8\) and we act by the element \(\mathbf{s}_7\in W^{\mathrm{aff}}_{\Lambda}\) to finally obtain
\begin{align*}
\mathbf{s}_7\mathbf{s}_3\mathbf{s}_4\mathbf{s}_2\mathbf{s}_3(\Lambda+\varepsilon_4+\varepsilon_5+\varepsilon_8)&=(\Lambda_0+\Lambda_1+\Lambda_5)+(\Lambda_2-\Lambda_1+\Lambda_3-\Lambda_2+\Lambda_6-\Lambda_5)\\
&=\Lambda_0+\Lambda_3+\Lambda_6,
\end{align*}
which is dominant and \(z=\mathbf{s}_7\mathbf{s}_3\mathbf{s}_4\mathbf{s}_2\mathbf{s}_3\in W^{\mathrm{aff}}_{\Lambda}\).
\end{example}

\subsection{}The commutator subalgebra $[\widehat{\mathfrak{g}},\widehat{\mathfrak{g}}]$ modulo the center is the loop algebra $\mathfrak{g}[t^{\pm}]$ and note that the element $d$ defines a grading on the loop algebra. The $\mathbb{Z}_+$-graded subalgebra  $\mathfrak{g}[t]:=\mathfrak{g}\otimes \mathbb{C}[t]$ of the loop algebra is the current algebra associated to $\mathfrak{g}$. 
Then $\mathbf{U}(\mathfrak{g}[t])$ inherits a grading where an element $(a_1\otimes t^{r_1})\cdots (a_s\otimes t^{r_s})$, $a_j\in\lie g$, $r_j\in\mathbb Z_+$ for $1\le j\le s$ will have grade $r_1+\cdots+r_s$. We denote by $\mathbf U(\lie g[t])_k$ be the homogeneous component of degree $k$ and recall that it is a $\lie g$--module for all $k\in\mathbb{Z}_+$. In the rest of the paper we abbreviate  $$\mathbf{U}=\mathbf{U}(\mathfrak{g}[t]),\ \ \mathbf{U}^{\pm}=\mathbf{U}(\mathfrak{n}^{\pm}[t]),\ \ \mathbf{U}^0=\mathbf{U}(\mathfrak{h}[t]).$$ 
So as vector spaces,
$$\mathbf{U}\cong \mathbf{U}^-\otimes \mathbf{U}^{0} \otimes \mathbf{U}^+.$$
A finite-dimensional $\mathbb{Z}$-graded $\lie g[t]$-module is a $\mathbb{Z}$-graded vector space admitting a compatible graded action of $\lie g[t]$:
$$V=\bigoplus_{k\in \mathbb{Z}} V[k],\ \ (a\otimes t^r) V[k]\subseteq V[k+r]$$
In particular each graded piece $V[k]$ is a $\lie g$-module. If $\dim V[k]<\infty$ for all $k\in \mathbb{Z}$ we define the graded character as a formal sum 
$$\mathrm{ch}_{\mathrm{gr}} V=\sum_{k\in \mathbb{Z}} \mathrm{ch}_{\lie h} V[k]\  q^k$$
where 
$$\mathrm{ch}_{\lie h}(V')=\sum_{\mu\in P} (\dim V'_\mu) e_{\mu}\in\mathbb{Z}[P],\ \ V_{\mu}'=\{v\in V': h.v=\mu(h)v,\ \forall h\in \lie h\}$$ refers to the usual $\lie h$ -character of a $\lie g$-module $V'$ and $\mathbb{Z}[P]$ is the group algebra of $P$ with basis $e_{\mu},\ \mu\in P$.
Given a $\mathbb{Z}$-graded space $V$ let $\tau_p V$ be the graded space whose $r$-th graded piece is given by $V[r-p]$. 

\subsection{}\label{section26} Suppose now that we are given a cyclic $\lie g[t]$--module $V$ generated by a vector $v$. Define an increasing filtration 
$$0\subseteq V^0\subseteq V^1\subseteq\cdots $$ of $\lie g$--submodules of $V$ by $$V^k=\sum_{s=0}^k \mathbf U(\lie g[t])_s v.$$ The associated graded vector space $\gr V$ admits an action of $\lie g[t]$ given by: 
$$
x(v+V^{k})= x.v+ V^{k+s},\ \ x\in\lie g[t]_s,\ \ v\in V^{k+1}.
$$
Furthermore, $\gr V$ is a cyclic $\lie g[t]$--module with cyclic generator $\bar v$, the image of $v$ in $\gr  V$. The fusion product defined in \cite{FL99} is a $\lie g[t]$--module of the form $\gr V$ for a special choice of a cyclic $\lie g[t]$--module $V$ which we define next.\medskip 

\begin{defn}
Given a finite-dimensional, graded and cyclic $\lie g[t]$-module $W^1$ and $z\in\mathbb \bc$, let $W^1_z$ be the $\lie g[t]$--module with action 
$$(x\otimes t^r)u=(x\otimes (t+z)^r).u,\ x\in \lie g,\ u\in W^1,\ r\in \mathbb Z_+.$$
\end{defn}
So if we have two such modules $W^1$ and $W^2$ (as in the above definition) with cyclic generators $u$ and $u'$ the module 
$W^1_{z_1}\otimes W^2_{z_2}$ is cyclic for all $z_1\neq z_2$ with cyclic generator $u\otimes u'$ (see for example \cite[Proposition 1.4]{FL99}).
The \textit{fusion product} is defined as follows
$$W^1_{z_1}*W^1_{z_2}:=\gr (W^1_{z_1}\otimes W^2_{z_2}).$$ 

\begin{rem}Clearly the definition depends on the parameters $z_1$ and $z_2$. However, it is conjectured in \cite{FL99} (and proved in several special cases, see \cite{CV13,FoL07,KL14,Na15aa,Ra14} for example and references therein) that the structure is independent of the choice. The independence for the twofold fusion product of finite-dimensional simple $\lie g$-modules can be proven directly and a presentation is obtained for all rank two Lie algebras (see \cite{BK22a}). 
\end{rem}
\section{Kirillov-Reshetikhin crystals, Demazure crystals and modules}\label{section3}
In this section we collect some needed results on Demazure modules and crystals. 
\subsection{} An important class of representations for the standard Borel subalgebra $\widehat{\mathfrak{b}}$ is given by the so-called \textit{Demazure modules}. We restrict ourselves here to the $\lie g$-stable Demazure modules and use the simplified presentation of \cite{CV13} as a general definition (see also \cite{KV21,KV14} for the simplified presentation in general). 
\begin{defn} Let $\lambda\in P^+$ and $\ell\in\mathbb{N}$. The $\lie g$-stable Demazure module $\mathbf{D}_{\lambda}^{\ell}$ is the cyclic $\lie g[t]$-module with cyclic generator $v$ subject to the defining relations ($\alpha\in R^+$, $h\in \lie h$)
$$\lie n^{+}[t]v=(h\otimes t^{r+1})v=0,\ \left(x_{\alpha}^{-}\otimes 1\right)^{\lambda(h_{\alpha})+1}v,\ \ r\in \mathbb{Z}_+$$
$$(h\otimes 1)v=\lambda(h)v,\ \ \left(x_{\alpha}^{-}\otimes t^{s_{\alpha}}\right)v=0,\ \ \left(x_{\alpha}^{-}\otimes t^{s_{\alpha-1}}\right)^{m_{\alpha}+1}v=0,\ \text{ if $m_{\alpha}<\ell$}$$
\end{defn}
The Demazure module $\mathbf{D}_{\lambda}^{\ell}$ can be embedded into a level $\ell$ integrable irreducible highest weight module for the corresponding affine Lie algebra. Their structure has been intensively studied and various connections to number theory and algebraic combinatorics have been made. We list a few examples.
\begin{itemize}
    \item The character of a level one Demazure module is given by the specialization of a non-symmetric Macdonald polynomial \cite{I03}.
    \item The characters of certain generalized Demazure modules are given by graded Euler characteristics of certain vector bundles on the flag variety \cite{BMP20}.
    \item The level two Demazure modules appear as graded limits of irreducible finite-dimensional representations of quantum affine algebras \cite{BCMo15}. 
    \item They can be used to give upper bounds for the dimension of various representations of quantum affine algebras \cite{KV21}.

    \item The generating function for the graded multiplicities in excellent filtrations (see below for a definition) specialize to Mock-Theta functions, hypergeometric series \cite{BCK18,BCSV15} and admit a combinatorial formula in terms of two dimensional lattice path \cite{BK19a}.
   \item Demazure modules admit a fusion product structure \cite{CSVW16,FoL07,KV14}. For example, if we write \(\lambda=\lambda_0+\ell\lambda_1\) for some \(\lambda_0,\lambda_1\in P^+\) we get (which is known as the Steinberg type decomposition formula; see \cite[Theorem 1]{CSVW16}) 
\begin{align}\label{stde}
    \mathbf{D}^\ell_\lambda= \mathbf{D}^\ell_{\ell\lambda_1} * \mathbf{D}^\ell_{\lambda_0}.
\end{align}
\end{itemize}
Nevertheless, many basic questions remain unanswered (especially for higher level Demazure modules) e.g. their classical decompositions are unexplored except in a few special cases (see for example \cite{BK20b,FoL06}) or there is no closed dimension formula. The following can be derived from \cite{Na11} and is known in the literature as an \textit{excellent filtration}.
\begin{thm} \label{flagc} For any $k\geq \ell$, the Demazure module $\mathbf{D}_{\lambda}^{\ell}$ admits a filtration 
$$0=V_0\subseteq V_1\subseteq \cdots \subseteq V_m=\mathbf{D}_{\lambda}^{\ell}$$
such that 
$$V_i/V_{i-1}\cong \tau_{p_i}\mathbf{D}_{\mu_i}^{k},\ \ 1\leq i\leq m,\ (p_i,\mu_i)\in \mathbb{Z}_+\times P^+$$
\qed
\end{thm}
We denote by $[\mathbf{D}_{\lambda}^{\ell }:\tau_p\mathbf{D}_{\mu}^{k}]$ the multiplicity of $\tau_p\mathbf{D}_{\mu}^{k}$ in a level $k$ flag of $\mathbf{D}_{\lambda}^{\ell}$. Define further the polynomials and numerical multiplicities 
$$[\mathbf{D}_{\lambda}^{\ell}:\mathbf{D}_{\mu}^{k}]_q=\sum_{p\geq 0} \ [\mathbf{D}_{\lambda}^{\ell}:\tau_p\mathbf{D}_{\mu}^{k}] \ q^p$$
and the generating series by
$$\mathbf{A}^{\ell\rightarrow k}_{\mu}(x_1,\dots,x_n,q)=\sum_{(k_1,\dots,k_n)\in\mathbb{Z}_+^n}[\mathbf{D}_{\mu+\sum_{i=1}^nk_i\alpha_i}^{\ell}:\mathbf{D}_{\mu}^{k}]_q \ x_1^{k_1}\cdots x_n^{k_n}$$
$$\mathbf{A}^{\ell\rightarrow k}_{\mu}(x_1,\dots,x_n):=\mathbf{A}^{\ell\rightarrow k}_{\mu}(x_1,\dots,x_n,1).$$

\subsection{}\label{section32}  Given $\lambda\in P^+$ and $\ell\in\mathbb{N}$, let $w\in W^{\mathrm{aff}}$ be such that 
\begin{equation}\label{verw1}\Lambda=w^{-1}(w_0(\lambda)+\ell\Lambda_0)\in\widehat P^+.\end{equation}
Without loss of generality we assume $w$ to be the element of smallest length in the coset $w W^{\mathrm{aff}}_{\Lambda}$, i.e. we have
$$|w\mathbf{s}_j|=|w|+1,\ \ \forall \mathbf{s}_j\in W^{\mathrm{aff}}_{\Lambda}.$$
The \emph{socle} of $\mathbf{D}_{\lambda}^{\ell}$ which is by definition the largest semi-simple submodule is in fact irreducible (see for example \cite[Section 5.3.3]{BCKV22}) and isomorphic to $V(\mathrm{soc}(\lambda,\ell))$ where 
$$\mathrm{soc}(\lambda,\ell):=\Lambda_{|_{\lie h}},\ \ \Lambda=\mathrm{soc}(\lambda,\ell)+\ell\Lambda_0 \mod \delta.$$ 
We show in the rest of this subsection how the socle can be computed explicitly from $\lambda$ and $\ell$. Recall the equation from \eqref{teilmitell} and write 
$$s_{1,1}+\cdots+s_{1,n}\equiv h \hspace{-0,2cm} \mod (n+1),\ \  0\leq h\leq n.$$
Since $b_1\epsilon_1+\cdots+b_n\epsilon_n\in Q$ if and only if the sum over the coefficients is divisible by $(n+1)$, we get that 
$$t_{s_{1,n}\epsilon_1+\cdots+s_{1,1}\epsilon_n}\in \tau^{-1} \cdot  W^{\mathrm{aff}}$$
where $\tau$ is the Dynkin diagram automorphism given by $\tau(r)=r-h \mod (n+1)$ for $0\leq r\leq n+1$. Furthermore, we can order the remainders
\begin{equation}\label{sort5}0\leq m_{1,\sigma(1)}\leq m_{1,\sigma(2)}\leq \cdots \leq m_{1,\sigma(n)}\leq \ell\end{equation}
for a suitable $\sigma\in \Sigma_{n+1}$ with $\sigma(n+1)=n+1$.
Since $w_0(\epsilon_j)=\epsilon_{n+2-j}$ we get 
$$ w_0(\lambda)=\lambda(h_{1,n})\epsilon_{n+1}+\lambda(h_{2,n})\epsilon_n+\cdots+\lambda(h_{n,n})\epsilon_2$$
and hence (the indices are again considered modulo $n+1$)
\begin{align*}w_0(\lambda)+\ell\Lambda_0 &\xmapsto[]{t_{s_{1,n}\epsilon_1+\cdots+s_{1,1}\epsilon_n}}(\ell-m_{1,n})\epsilon_{1}+(\ell-m_{1,n-1})\epsilon_2+\cdots+(\ell-m_{1,1})\epsilon_n+\ell\Lambda_0&\\&
\xmapsto[]{ \ \  \ \tilde{\sigma} \ \ \ \  } (\ell-m_{1,\sigma(1)})\epsilon_{1}+(\ell-m_{1,\sigma(2)})\epsilon_2+\cdots+(\ell-m_{1,\sigma(n)})\epsilon_n+\ell\Lambda_0&\\&
\xmapsto[]{\ \ \  \tau \ \ \ \  } (\ell-m_{1,\sigma(1)})\epsilon_{1-h}+(\ell-m_{1,\sigma(2)})\epsilon_{2-h}+\cdots+(\ell-m_{1,\sigma(n)})\epsilon_{n-h}+\ell\Lambda_{-h}
\end{align*}
where $\tilde{\sigma}\in \Sigma_{n+1}$ is the permutation determined by $\tilde{\sigma}(\epsilon_{n+1-j})=\epsilon_{\sigma^{-1}(j)}$. Since the image lies in $\widehat{P}^+$ it coincides with $\Lambda$ and we have derived the following. 

\begin{prop}\label{stabg439}
Let $\lambda\in P^+$, $\ell\in\mathbb{N}$ and $w\in W^{\mathrm{aff}}$ as in \eqref{verw1}. Then $\mathrm{soc}(\lambda,\ell)=\Lambda_{|_{\lie h}}$ where 
    \begin{equation}\label{gl3}\tau \tilde{\sigma} t_{s_{1,n}\epsilon_1+\cdots+s_{1,1}\epsilon_n}\in W^{\mathrm{aff}}_{\Lambda} w^{-1},\ \ \ \Lambda=\sum_{j=0}^n a_j \Lambda_j,\ \ a_{j-h}=m_{1,\sigma(j+1)}-m_{1,\sigma(j)},\end{equation}
and we understand $m_{1,\sigma(0)}=0$ and $m_{1,\sigma(n+1)}=\ell$.
\qed
\end{prop}
\begin{example}
    We consider \(n+1=6,\ \ell=5\) and \(\lambda=4\varpi_1+3\varpi_2+5\varpi_3+\varpi_4+3\varpi_5\). 
    We first compute
    \[s_{1,1}+s_{1,2}+s_{1,3}+s_{1,4}+s_{1,5}=1+2+3+3+4=13 =6\cdot 2+1 \implies h=1.\]
    Next we compute the remainders
    \[m_{1,1}=4,m_{1,2}=2,m_{1,3}=2,m_{1,4}=3,m_{1,5}=1\]
    and order these
    \[0\leq m_{1,5}\leq m_{1,2}\leq m_{1,3}\leq m_{1,4} \leq m_{1,1}\leq \ell.\]
    So the permutation \(\sigma\in\Sigma_{n+1}\) can be chosen to be the transposition \( (1\ 5)\).
    According to Proposition~\ref{stabg439} the differences between subsequent pairs in this chain give the coefficients of \(\Lambda\) shifted by \(-h=-1\) and thus 
    \begin{align*}  
    \Lambda&=(m_{1,5}-0) \Lambda_5+(m_{1,2}-m_{1,5}) \Lambda_0+ \cdots+(\ell-m_{1,1}) \Lambda_4 \\
    &=\Lambda_0+\Lambda_2+\Lambda_3+\Lambda_4+\Lambda_5
    \end{align*}
    which gives \(\operatorname{soc}(\lambda,5)=\varpi
    +\varpi_3+\varpi_4+\varpi_5\).
    
\end{example}

\subsection{} We keep the notation established in Section~\ref{section32}. In particular, recall the definition of $\sigma$ and $\tilde{\sigma}$ as elements of $\Sigma_{n+1}$ and the Dynkin diagram automorphism $\tau$ which maps $r$ to $r-h$. Given $u\in W$ with
$u(\varpi_i)=\epsilon_{k_1}+\cdots+\epsilon_{k_i}\in W(\varpi_i)$ we define the \textit{shifted element} of $u(\varpi_i)$ by
$$u^{\mathrm{sh}}(\varpi_i):=\epsilon_{\sigma^{-1}(k_1-1)-h}+\cdots+\epsilon_{\sigma^{-1}(k_i-1)-h}$$
Recall the equation in \eqref{verw1}. 
\begin{lem}\label{orbitcond49}
Suppose that $u\in W, \tilde{w}\in W^{\mathrm{aff}}$ are such that
$$\widetilde{\Lambda}=\tilde{w}^{-1}(w_0(\lambda+u(\varpi_i))+\ell \Lambda_0)\in \widehat{P}^+.$$   
Then $\widetilde{\Lambda}$ is the unique dominant integral weight in the orbit $W^{\mathrm{aff}}_{\Lambda}\left(\Lambda+u^{\mathrm{sh}}(\varpi_i)\right)$ and $\tilde{w}\in w \cdot W^{\mathrm{aff}}_{\Lambda}W^{\mathrm{aff}}_{\widetilde{\Lambda}}$. In particular, given $u,\tilde{u}\in W$ we have
$$\mathrm{soc}(\lambda+u(\varpi_i),\ell)=\mathrm{soc}(\lambda+\tilde{u}(\varpi_i),\ell)\iff \exists z\in W^{\mathrm{aff}}_{\Lambda}: z(u^{\mathrm{sh}}(\varpi_i))=\tilde{u}^{\mathrm{sh}}(\varpi_i).$$
\begin{proof} Let $u(\varpi_i)=\epsilon_{k_1}+\cdots+\epsilon_{k_i}$. We obtain
\begin{align*}\hspace{1cm}  w^{-1}\tilde{w}\widetilde{\Lambda}&=\Lambda+w^{-1}w_0u(\varpi_i)&\\&
=\Lambda +y\tau \tilde{\sigma} t_{s_{1,n}\epsilon_1+\cdots+s_{1,1}\epsilon_n} w_0u(\varpi_i)\ \ \ \ \ \ \  \text{ (for some $y\in W^{\mathrm{aff}}_{\Lambda}$ by Proposition~\ref{stabg439})} &\\&
=\Lambda+y\tau \tilde{\sigma} (\epsilon_{n+2-k_1}+\cdots+\epsilon_{n+2-k_i})&\\&
=\Lambda+y\epsilon_{\sigma^{-1}(k_1-1)-h}+\cdots+\epsilon_{\sigma^{-1}(k_i-1)-h}&\\&
=\Lambda+yu^{\mathrm{sh}}(\varpi_i)
\end{align*}
Now with Lemma~\ref{stabwcv2} there exists $y'\in W^{\mathrm{aff}}_{\Lambda}$ with $$\Lambda+y'y(u^{\mathrm{sh}}(\varpi_i))\in \widehat{P}^{+}\implies \widetilde{\Lambda}=\Lambda+y'y(u^{\mathrm{sh}}(\varpi_i))$$ and the first part follows. The second part of the lemma is immediate.
\end{proof}
\end{lem}
\subsection{} In this subsection we introduce the following subset of $W(\varpi_i)$ which will play an important role in this paper:
\begin{equation}\label{mainset}\mathbf{R}^{\ell}_{\lambda,i}=\left\{\epsilon_{k_1}+\cdots+\epsilon_{k_i}: m_{\epsilon_r-\epsilon_{k_s}}<\ell \ \ \forall  s\in\{1,\dots,i\},\  \forall r\notin\{k_1,\dots,k_i\},\ r<k_s\right\}.\end{equation}
Our aim in this subsection is to show that the socle of a Demazure module can distinguish the elements of $\mathbf{R}^{\ell}_{\lambda,i}$ (see Proposition~\ref{diffsoc}). We first record the following lemma. 
\begin{lem} \label{lemsoc5} Let $u(\varpi_i)=\epsilon_{k_1}+\cdots+\epsilon_{k_i}\in W(\varpi_i)$ and $z\in W^{\mathrm{aff}}_{\Lambda}$ such that
$$z(u^{\mathrm{sh}}(\varpi_i))=\epsilon_{\sigma^{-1}(\tilde{k}_1-1)-h}+\cdots+\epsilon_{\sigma^{-1}(\tilde{k}_i-1)-h}$$
Then there is a permutation $\pi\in \Sigma_i$ such that for all $d\in\{1,\dots,i\}$ we have
$$m_{\epsilon_{k_d}-\epsilon_{\tilde{k}_{\pi(d)}}}=\ell\ \ \text{or} \ \ 
m_{\epsilon_{\tilde{k}_{\pi(d)}}-\epsilon_{k_d}}=\ell\ \ \text{or} \ \  k_d=\tilde{k}_{\pi(d)}.$$
\begin{proof}
All calculations are modulo $(n+1)$ and the proof is by induction on the length of $z$. If $z=\mathrm{Id}$ there is nothing to show. So write $z=\mathbf{s}_{j}z'$ and suppose by induction that 
$$z'(u^{\mathrm{sh}}(\varpi_i))=\epsilon_{\sigma^{-1}(g_1-1)-h}+\cdots+\epsilon_{\sigma^{-1}(g_i-1)-h}$$
with the desired properties:
\begin{equation}\label{hhzt24} m_{\epsilon_{k_d}-\epsilon_{g_{\kappa(d)}}}=\ell\ \  \text{or} \ \ m_{\epsilon_{g_{\kappa(d)}}-\epsilon_{k_{d}}}=\ell\ \ \text{or}\ \ k_d=g_{\kappa(d)},\ \ \kappa\in \Sigma_i.\end{equation}
If $\mathbf{s}_{j}$ acts trivially there is nothing to show, so let $j=\sigma^{-1}(g_{\kappa(r)}-1)-h-v$ for some $r\in\{1,\dots,i\}$ and $v\in\{0,1\}$; we set $v'=(1-2\delta_{v,0})$. Note that the application of $\mathbf{s}_j$ replaces $\epsilon_{\sigma^{-1}(g_{\kappa(r)}-1)-h}$ by $\epsilon_{\sigma^{-1}(g_{\kappa(r)}-1)-h+v'}$ and thus the newly obtained element (after applying the reflection $\mathbf{s}_j$) is the shifted element of 
$$\epsilon_{g_{\kappa(1)}}+\cdots+\epsilon_{g_{\kappa(r-1)}}+\epsilon_{g^1_{\kappa(r)}}+\epsilon_{g_{\kappa(r+1)}}+\cdots+\epsilon_{g_{\kappa(i)}}\ \ g^1_{\kappa(r)}=\sigma(\sigma^{-1}(g_{\kappa(r)}-1)+v')+1$$
Since $\mathbf{s}_j\in W^{\mathrm{aff}}_{\Lambda}$ we get with \eqref{gl3}
$$m_{\epsilon_{g^1_{\kappa(r)}}-\epsilon_{g_{\kappa(r)}}}=\ell \ \ \text{or} \ \ m_{\epsilon_{g_{\kappa(r)}}-\epsilon_{g^1_{\kappa(r)}}}=\ell$$
and thus the desired property
$m_{\epsilon_{g^1_{\kappa(r)}}-\epsilon_{k_{r}}}=\ell$ or $m_{\epsilon_{k_{r}}-\epsilon_{g^1_{\kappa(r)}}}=\ell$ follows together with our induction hypothesis \eqref{hhzt24}.
\end{proof}
\end{lem}
Now we are able to prove the main result of this subsection.
\begin{prop}\label{diffsoc}
Let $u(\varpi_i),\tilde{u}(\varpi_i)\in \mathbf{R}^{\ell}_{\lambda,i}$ two different elements. Then we have 
$$\mathrm{soc}(\lambda+u(\varpi_i),\ell)\neq \mathrm{soc}(\lambda+\tilde{u}(\varpi_i),\ell)$$
\begin{proof} Assume that the socle of both elements coincides and set
$$u(\varpi_i)=\epsilon_{k_1}+\cdots+\epsilon_{k_i},\ \ \tilde{u}(\varpi_i)=\epsilon_{\tilde{k}_1}+\cdots+\epsilon_{\tilde{k}_i}$$ By Lemma~\ref{orbitcond49} and Lemma~\ref{lemsoc5} there exists $\pi\in \Sigma_i$ such that for all $d\in\{1,\dots,i\}$ we have
$$k_d=\tilde{k}_{\pi(d)} \ \ \text{or} \ \ m_{\epsilon_{k_d}-\epsilon_{\tilde{k}_{\pi(d)}}}=\ell \ \ \text{or} \ \ m_{\epsilon_{\tilde{k}_{\pi(d)}}-\epsilon_{k_d}}=\ell$$
By removing the indices with $k_d=\tilde{k}_{\pi(d)}$ we can assume without loss of generality that we are left with two sets
$$\{k_1,\dots,k_i\}\ \ \text{and}\ \ \{\tilde{k}_{\pi(1)},\dots,\tilde{k}_{\pi(i)}\}$$
such that for all $d\in \{1,\dots,i\}$
$$m_{\epsilon_{k_d}-\epsilon_{\tilde{k}_{\pi(d)}}}=\ell\ \ \text{or} \ \ 
m_{\epsilon_{\tilde{k}_{\pi(d)}}-\epsilon_{k_d}}=\ell.$$
Since both elements lie in $\mathbf{R}^{\ell}_{\lambda,i}$ we also have the property 
$$m_{\epsilon_{\tilde{k}_p}-\epsilon_{k_s}}=\ell \implies 
\tilde{k}_p\in\{k_1,\dots,k_i\}\ \ \ \  \ m_{\epsilon_{k_s}-\epsilon_{\tilde{k}_p}}=\ell \implies k_s\in\{\tilde{k}_1,\dots,\tilde{k}_i\}.$$
So the following claim finishes the proof of the proposition.

\textit{Claim:} Let $\{k_1,\dots,k_i\},\{g_1,\dots,g_i\}$ be two sets with
$$m_{\epsilon_{k_d}-\epsilon_{g_d}}=\ell \ \ \text{or} \ \ m_{\epsilon_{g_d}-\epsilon_{k_d}}=\ell\ \ \ \forall d\in\{1,\dots,i\}.$$
Moreover, assume that $m_{\epsilon_{g_p}-\epsilon_{k_s}}=\ell$ (resp. $m_{\epsilon_{k_s}-\epsilon_{g_p}}=\ell$) implies $g_p\in\{k_1,\dots,k_i\}$ (resp. $k_s\in\{g_1,\dots,g_i\}$). Then both sets are equal. 

\vspace{0,2cm}

\textit{Proof of the Claim:} The proof is by induction on the cardinality. If $g_d\in\{k_1,\dots,k_i\}$ for all $d$, we are done. Otherwise there exists $r\in\{1,\dots,i\}$ such that $g_r\notin \{k_1,\dots,k_i\}$ and therefore
\begin{equation}\label{kkk2}k_r=g_j \text{ for some } j\in\{1,\dots,i\}\ \text{ and } \ m_{\epsilon_{k_r}-\epsilon_{g_r}}=\ell\end{equation}
We consider the sets
\begin{equation}\label{ggrew}\{k_1,\dots,k_i\}\backslash\{k_r\}\ \ \text{and}\ \ \{g_1,\dots,g_i\}\backslash\{g_j=k_r\}\end{equation}
and show that they satisfy the initial conditions where we match $k_q$ with $g_q$ for $q\neq j$ and $k_j$ with $g_r$. Since $m_{\epsilon_{k_r}-\epsilon_{k_j}}=m_{\epsilon_{g_j}-\epsilon_{k_j}}=\ell$ 
or $m_{\epsilon_{k_j}-\epsilon_{k_r}}=m_{\epsilon_{k_j}-\epsilon_{g_j}}=\ell$ we conclude with \eqref{kkk2} 
$$m_{\epsilon_{g_r}-\epsilon_{k_j}}=\ell \ \ \text{or} \ \ m_{\epsilon_{k_j}-\epsilon_{g_r}}=\ell.$$
Now suppose that $m_{\epsilon_{g_p}-\epsilon_{k_s}}=\ell$ for suitable $p\in\{1,\dots,i\}\backslash\{j\}$ and $s\in\{1,\dots,i\}\backslash\{r\}$. This gives $g_p\in\{k_1,\dots,k_i\}\backslash\{k_r\}$, since $g_p\in\{k_1,\dots,k_i\}$ by our assumptions and $g_p=k_r$ is impossible because of $j\neq p$. Similarly if $m_{\epsilon_{k_s}-\epsilon_{g_p}}=\ell$ for suitable $p\in\{1,\dots,i\}\backslash\{j\}$ and $s\in\{1,\dots,i\}\backslash\{r\}$ gives $k_s\in\{g_1,\dots,g_i\}\backslash\{g_j\}$. Therefore we can apply induction and the claim is proven.
\end{proof}
\end{prop} 
\subsection{} 
A seminormal abstract crystal is a set $B$ endowed with crystal operators $\tilde{e}_i, \tilde{f}_i \colon B \to B \sqcup \{0\}$, for $i \in I$, and weight function $\mathrm{wt}: B \to P$ that satisfy the following conditions:
\begin{itemize}
\item[(1)] $\widetilde{\varphi}_i(b) = \widetilde{\epsilon}_i(b) + \mathrm{wt}(b)(h_i)$, for all $b \in B$ and $i \in I$,

\item[(2)] $\tilde{f}_i b = b'$ if and only if $b = \tilde{e}_i b'$, for $b, b' \in B$ and $i \in I$,

\item[(3)] $\mathrm{wt}(\tilde{f}_i b) = \mathrm{wt}(b) - \alpha_i$ if $\tilde{f}_i b \neq 0$
\end{itemize}
where $\widetilde{\epsilon}_i, \widetilde{\varphi}_i \colon  B \to \mathbb{Z}_{\geq 0}$ are defined by
$$
\widetilde{\epsilon}_i(b) := \max \{ k : \tilde{e}_i^k b \neq 0 \}\,,
\qquad  \widetilde{\varphi}_i(b) := \max \{ k : \tilde{f}_i^k b \neq 0 \}\,.
$$
We define the tensor product of seminormal abstract crystals $B_1$ and $B_2$, denoted by $B_1 \otimes B_2$, as the Cartesian product $B_1 \times B_2$ with the following crystal structure:
\begin{align*}
\tilde{e}_i(b_1 \otimes b_2) & := \begin{cases}
\tilde{e}_i b_1 \otimes b_2 & \text{if } \widetilde{\epsilon}_i(b_1) > \widetilde{\varphi}_i(b_2)\,, \\
b_1 \otimes \tilde{e}_i b_2 & \text{if } \widetilde{\epsilon}_i(b_1) \leq \widetilde{\varphi}_i(b_2)\,,
\end{cases}
\\ \tilde{f}_i(b_1 \otimes b_2) & := \begin{cases}
\tilde{f}_i b_1 \otimes b_2 & \text{if } \widetilde{\epsilon}_i(b_1) \geq \widetilde{\varphi}_i(b_2)\,, \\
b_1 \otimes \tilde{f}_i b_2 & \text{if } \widetilde{\epsilon}_i(b_1) < \widetilde{\varphi}_i(b_2)\,,
\end{cases}
\end{align*}
$$ \widetilde{\epsilon}_i(b_1 \otimes b_2)  := \max\{\widetilde{\epsilon}_i(b_2), \widetilde{\epsilon}_i(b_1) - \mathrm{wt}(b_2)(h_i)\}$$
$$ \widetilde{\varphi}_i(b_1 \otimes b_2)  := \max\{\widetilde{\varphi}_i(b_1), \widetilde{\varphi}_i(b_2) + \mathrm{wt}(b_1)(h_i)\}$$
$$\mathrm{wt}(b_1 \otimes b_2) := \mathrm{wt}(b_1) + \mathrm{wt}(b_1)$$
Kashiwara has shown in \cite{K91} that the irreducible highest weight module $V(\lambda)$ (resp. $V(\Lambda)$) for $\lambda \in P^+$ (resp. $\Lambda \in \widehat{P}^+$) admits a crystal basis, denoted by $B(\lambda)$ (resp. $B(\Lambda)$). For a precise definition we refer to \cite{K91} and point out that $B(\lambda)$ and $B(\Lambda)$ respectively satisfy the above mentioned axioms with an appropriate index set $I$. Moreover, the corresponding (affine) Weyl group acts on the crystal basis (see for example \cite{K94}) by
\begin{equation}\label{z543e9}\mathbf{S}_ib:=\begin{cases}\tilde{f}_i^{\mathrm{wt}(b)(h_i)}b,& \text{if $\mathrm{wt}(b)(h_i)\geq 0$}\\
\tilde{e}_i^{-\mathrm{wt}(b)(h_i)}b,& \text{if $\mathrm{wt}(b)(h_i)\leq 0$}
    \end{cases}\end{equation}
If $z=\mathbf{s}_{j_1}\cdots \mathbf{s}_{j_r}$ is a reduced expression, we define $\mathbf{S}_{z}:=\mathbf{S}_{j_1}\cdots \mathbf{S}_{j_r}$ which is independent of the choice of the reduced expression and we have $\mathrm{wt}(\mathbf{S}_{z}b)=z\mathrm{wt}(b)$.
\subsection{}\label{section36} We denote by $B^{r,s}$ the \textit{Kirillov--Reshetikhin crystal}, where $r \in \{1,\dots,n\}$ and $s \in \mathbb{N}$. These crystals are the crystal bases of a particular class of representations (called Kirillov-Reshetikhin modules) for the corresponding quantum affine algebra. We omit the details and work instead with a very explicit combinatorial model for KR crystals. We identify $B^{r,s}$ with the set of all tuples $(a_{p,q})_{1\leq p\leq r,\ r\leq q\leq n}$ of non-negative integers 
such that $\sum^k_{j=1}a_{\beta(j)}\leq s$ for all sequences \[
(\beta(1),\dots, \beta(k)), \ k\ge 1
\] satisfying the following: $\beta(1)=(1,r), \beta(k)=(r,n)$ and if $\beta(j)=(p,q)$ then the next element in the sequence is either $\beta(j+1)=(p,q+1)$  or $\beta(j+1)=(p+1,q).$
The crystal structure on these tuples is described in \cite{K13a,K16a} including also the description of the $R$-matrix, ground state path or highest weight elements of twofold tensor products. 
\begin{example}
For $\mathfrak{sl}_3$ the tuples in $B^{1,s}$ can be visualized as elements
$$b_{u,v}:=\ytableausetup{boxsize=1.2em}
\begin{ytableau}
u \\
v 
\end{ytableau}\, \ \ \ \  \ \ u+v\leq s$$
of weight $s\varpi_1-u\alpha_1-v\alpha_{1,2}$ and (whenever the Kashiwara operator acts)
$$\tilde{f}_0 b_{u,v}=b_{u,v-1},\ \ \ \  \tilde{f}_1b_{u,v}=b_{u+1,v},\ \ \ \ \tilde{f}_2b_{u,v}=b_{u-1,v+1}.$$
Similarly the tuples in $B^{2,s}$ can be visualized as elements
$$b_{u,v}:=\ytableausetup{boxsize=1.2em}
\begin{ytableau}
u &
v 
\end{ytableau}\, \ \ \ \  \ \ u+v\leq s$$
of weight $s\varpi_2-u\alpha_{1,2}-v\alpha_{2}$ and (whenever the Kashiwara operator acts)
$$\tilde{f}_0b_{u,v}=b_{u-1,v}, \ \ \ \ \tilde{f}_1b_{u+1,v-1}=b_{u,v}, \ \ \  \ \tilde{f}_2b_{u,v}=b_{u,v+1}.$$
\end{example}
\subsection{}
Let $\lambda \in P^+$, $\ell\in\mathbb{N}$ and $w\in W^{\mathrm{aff}}$ as in \eqref{verw1}. Kashiwara showed in \cite{K93} that the Demazure module $\mathbf{D}^{\ell}_{\lambda}$ admits a crystal basis in a suitable sense that is compatible with the crystal basis $B(\Lambda)$. For a fixed reduced expression 
$w = \mathbf{s}_{j_1} \mathbf{s}_{j_2} \cdots \mathbf{s}_{j_{t}}$, we denote by $\mathcal{F}_w$ the set of crystal operators of the form $\tilde{f}_{j_1}^{n_1} \tilde{f}_{j_2}^{n_2} \cdots \tilde{f}_{j_{t}}^{n_t}$ for $n_1,\dots,n_t\in\mathbb{Z}_+$. The \textit{Demazure crystal} of $\mathbf{D}^{\ell}_{\lambda}$ is the full subcrystal of $B(\Lambda)$ given by
\begin{equation}
\mathbf{B}^{\ell}_{\lambda}:= \mathcal{F}_w b_{\Lambda}\backslash\{0\}
\end{equation}
where $b_{\Lambda}$ is the unique element in $B(\Lambda)$ with $\mathrm{wt}(b_{\Lambda})=\Lambda.$ 
The following can be found in \cite[Theorem 3.2]{LS19a} in its most general form (without the perfectness assumption on KR crystals); see also \cite[Proposition 8.1]{ST12} for the first part. 
\begin{thm}\label{schuting}
Let $B$ be a tensor product of KR crystals of level bounded by $\ell$. Then there exists $\Lambda^{(1)},\dots,\Lambda^{(d)}\in \widehat{P}^+$ of level $\ell$, i.e. $\Lambda^{(j)}(c)=\ell$ and $\xi_1,\dots,\xi_d\in P^+$ such that
$$B\otimes B(\ell \Lambda_0)\cong B(\Lambda^{(1)})\sqcup \cdots \sqcup B(\Lambda^{(d)})$$
and the restriction gives rise to an isomorphism
\begin{equation}\label{decdem0}
B \otimes b_{\ell \Lambda_0} \cong \mathbf{B}^{\ell}_{\xi_1}\sqcup \cdots \sqcup\mathbf{B}^{\ell}_{\xi_d} 
\end{equation}
\qed
\end{thm}
The main result of this article will allow to determine recursively the decomposition in \eqref{decdem0} for tensor products of fundamental KR crystals.
\subsection{} Given $\lambda\in P^+$ with $\lambda(h_1),\dots,\lambda(h_n)\leq \ell$ (the level needs to be bounded by $\ell$) we consider the tensor product of KR crystals 
$$\mathbf{B}_{\otimes \lambda}:=B^{1,\lambda(h_1)}\otimes \cdots \otimes B^{n,\lambda(h_n)}$$
and obtain with Theorem~\ref{schuting}
\begin{equation}\label{decdem01}\mathbf{B}_{\otimes \lambda}\otimes b_{\ell \Lambda_0}\cong \mathbf{B}^{\ell}_{\xi_1}\sqcup \cdots \sqcup\mathbf{B}^{\ell}_{\xi_d}.\end{equation}
 From abstract crystal theory it is clear that $\mathbf{B}^{\ell}_{\lambda}$ appears in \eqref{decdem01}; say $\xi_1=\lambda$.
 \begin{thm}\label{crysesti1} Let $\lambda\in P^+$ whose coefficients are bounded by $\ell\in \mathbb{N}$ and $w\in W^{\mathrm{aff}}$ as in \eqref{verw1}. For any $i\in\{1,\dots,n\}$, the tensor product $(B^{i,1}\otimes \mathbf{B}^{\ell}_{\lambda})$ decomposes into a disjoint union of Demazure crystals 
 $$\mathbf{B}^{\ell}_{\nu_1}\sqcup \cdots \sqcup\mathbf{B}^{\ell}_{\nu_p}.$$
Moreover, for each $\chi\in \mathbf{R}_{\lambda,i}^{\ell}$ there exists an index $j_\chi\in\{1,\dots,p\}$ with $\mathbf{B}^{\ell}_{\lambda+\chi}\subseteq \mathbf{B}^{\ell}_{\nu_{j_\chi}}$ and the assignment
$\chi\mapsto j_{\chi}$ is injective. In particular, 
 \begin{equation}\label{unglcrys12}\left|B^{i,1}\otimes \mathbf{B}^{\ell}_{\lambda}\right|\geq \sum_{\chi\in\mathbf{R}_{\lambda,i}^{\ell}}\left|\mathbf{B}^{\ell}_{\lambda+\chi}\right|.\end{equation}
   \begin{proof}
From Theorem~\ref{schuting} and the discussion preceeding the theorem we can decompose 
\begin{equation*}B^{i,1}\otimes \mathbf{B}_{\otimes \lambda}\otimes b_{\ell \Lambda_0}\cong (B^{i,1}\otimes \mathbf{B}^{\ell}_{\xi_1})\sqcup \cdots \sqcup (B^{i,1}\otimes\mathbf{B}^{\ell}_{\xi_d})\cong  \mathbf{B}^{\ell}_{\nu'_1}\sqcup \cdots \sqcup \mathbf{B}^{\ell}_{\nu'_{d'}}\end{equation*}
and hence 
$$(B^{i,1}\otimes \mathbf{B}^{\ell}_{\lambda})\cong \mathbf{B}^{\ell}_{\nu_{1}}\sqcup \cdots \sqcup \mathbf{B}^{\ell}_{\nu_{p}}$$
for a suitable subset $\{\nu_1,\dots,\nu_p\}\subseteq \{\nu_1',\dots,\nu'_{d'}\}$. Recall that $B^{i,1}$ is a crystal of a minuscule representation and hence there is a unique element $b_{z}$ in $B^{i,1}$ of weight $z\in W(\varpi_i)$. Given $\chi=u(\varpi_i)\in \mathbf{R}_{\lambda,i}^{\ell}$ there exists by Lemma~\ref{orbitcond49}  elements $y,y'\in W^{\mathrm{aff}}_{\Lambda}$ such that
$$\widetilde{\Lambda}:=yw^{-1}\left(w_0(\lambda+u(\varpi_i))+\ell \Lambda_0\right)=\Lambda+y'u^{\mathrm{sh}}(\varpi_i)\in \widehat{P}^+$$
and hence the element 
$$b_{y'u^{\mathrm{sh}}(\varpi_i)}\otimes b_{\mathrm{soc}(\lambda,\ell)}\otimes b_{\ell\Lambda_0}\in B^{i,1}\otimes \mathbf{B}^{\ell}_{\lambda}$$ 
is a highest weight element inside $B(\widetilde{\Lambda})$, where $b_{\mathrm{soc}(\lambda,\ell)}\otimes b_{\ell\Lambda_0}\in \mathbf{B}^{\ell}_{\lambda}$ is the unique element of weight $\Lambda$. In particular
\begin{equation}\label{swwoh}\mathcal{F}_w(b_{\mathrm{soc}(\lambda,\ell)}\otimes b_{\ell\Lambda_0})\in \mathbf{B}^{\ell}_{\lambda}\sqcup \{0\}\end{equation}
Now using the Weyl group action \eqref{z543e9} (recall that we can view all crystals inside a direct sum of highest weight crystals) we obtain with $|w\cdot y^{-1}|=|w|+|y^{-1}|$ 
and \eqref{swwoh}  
\begin{align*}\mathbf{S}_{wy^{-1}}&\left(b_{y'u^{\mathrm{sh}}(\varpi_i)}\otimes b_{\mathrm{soc}(\lambda,\ell)}\otimes b_{\ell\Lambda_0}\right)=\mathbf{S}_{w}\left(\mathbf{S}_{y^{-1}}b_{y'u^{\mathrm{sh}}(\varpi_i)}\otimes b_{\mathrm{soc}(\lambda,\ell)}\otimes b_{\ell\Lambda_0}\right)&\\&\subseteq \left(\mathcal{F}_{w}\mathbf{S}_{y^{-1}}b_{y'u^{\mathrm{sh}}(\varpi_i)}\otimes \mathcal{F}_{w}(b_{\mathrm{soc}(\lambda,\ell)}\otimes b_{\ell\Lambda_0})\right)\backslash\{0\}\subseteq B^{i,1}\otimes \mathbf{B}^{\ell}_{\lambda} \end{align*}
Hence we found an extremal weight element in $B^{i,1}\otimes \mathbf{B}^{\ell}_{\lambda}$ of weight 
$w_0(\lambda+u(\varpi_i))+\ell \Lambda_0$. So the aforementioned element lies inside a Demazure crystal $\mathbf{B}^{\ell}_{\nu_{j_{\chi}}}$ for an appropriate index $j_{\chi}$. In particular, $\mathbf{B}^{\ell}_{\nu_{j_{\chi}}}$ contains $\mathbf{B}^{\ell}_{\lambda+\chi}$ (see for example \cite[Proposition 4.4]{K93}). The fact that $\mathbf{B}^{\ell}_{\nu_{j_{\chi}}}$ contains no other $\mathbf{B}^{\ell}_{\lambda+\chi'}$ with $\chi'\in \mathbf{R}^{\ell}_{\lambda,i}$ follows from Proposition~\ref{diffsoc} since otherwise the socle would coincide.
     \end{proof} 
 \end{thm}
\section{The main results}\label{mainresultsection}
Here we state the main results. 
\subsection{}
For a tuple $\boldsymbol{\alpha}=(\alpha_{r_1,p_1},\dots,\alpha_{r_k,p_k})\in (R_i^+)^{\times k}$ 
we define 
$$s_{\boldsymbol{\alpha}}=\sum_{j=1}^k s_{r_{j},p_j},\ \ s_{\emptyset}=0,\ \ \kappa(\boldsymbol{\alpha})=(\alpha_{r_{\kappa(1)},p_1},\dots,\alpha_{r_{\kappa(k)},p_k}),\ \ s_{\boldsymbol{\alpha}}^{\mathrm{min}}=\min\left\{s_{\kappa(\boldsymbol{\alpha)}}: \kappa\in \Sigma_k\right\}$$
$$X_{\boldsymbol{\alpha}}=(x^-_{r_{1},p_1}\otimes t^{s_{r_{1},p_{1}}})\cdots (x^-_{r_{k},p_k}\otimes t^{s_{r_{k},p_{k}}}).$$
\begin{defn}
Let $\mathbf{M}^{\ell}_{\lambda,i}$ be the quotient of the level one module $\mathbf{D}^1_{\lambda+\varpi_i}$ with cyclic generator $m$ by the submodule generated by the following elements:

\begin{equation}\label{req1}X_{\boldsymbol{\alpha}}m,\ \text{ if $s_{\boldsymbol{\alpha}}\neq s_{\boldsymbol{\alpha}}^{\mathrm{min}}$},\  X_{\boldsymbol{\alpha}}m-\mathrm{sgn}(\kappa)X_{\kappa(\boldsymbol{\alpha})}m,\text{ if $s_{\boldsymbol{\alpha}}=s_{\kappa(\boldsymbol{\alpha})}$}\ \ \forall \ \boldsymbol{\alpha}\in \left(R_{\lambda,i}^+\right)^{\times k} \end{equation}

\begin{equation}\label{req2}
\left(x^-_{\alpha}\otimes t^{s_{\alpha}+\varpi_i(h_{\alpha})}\right)m,\ \ \left(x^-_{\alpha}\otimes t^{s_{\alpha}-1}\right)^{m_{\alpha}+1+\varpi_i(h_{\alpha})}m,\ \ \left(x^-_{\alpha}\otimes t^{s_{\alpha}}\right)^{2}m=0,\ \ \alpha\in R^+\end{equation}

\begin{equation}\label{req3}\left(x_{\alpha}^{-}\otimes t^{s_{\alpha}}\right)\left(x_{\beta}^{-}\otimes t^{s_{\beta}}\right)m,\text{ $\alpha,\beta\in R_{\lambda,i}^+$, } (\alpha,\beta)=1\end{equation}

\begin{equation}\label{req4}\left(x^-_{\gamma}\otimes t^{s_{\gamma}-1}\right)^{m_{\gamma}+\varpi_i(h_{\gamma})}\left(x_{\alpha}^-\otimes t^{s_{\alpha}}\right)m,\ \ \alpha\in R_{\lambda,i}^+,\ \gamma\in R^+,\ (\gamma,\alpha)=1\end{equation}

\begin{equation}\label{req5}\left(x^-_{\gamma}\otimes t^{s_{\gamma}-1}\right)^{m_{\gamma}}\left(x^-_{\alpha}\otimes t^{s_{\alpha}-1+\varpi_i(h_{\gamma})}\right)\left(x^-_{\beta}\otimes t^{s_{\beta}}\right)m,\ \ \alpha,\beta\in R_{\lambda,i}^+,\ \gamma\in R^+\end{equation}

where we assume in \eqref{req5} that
\begin{align*}&\gamma-\alpha,\ \gamma-\beta\in R^+, \text{ if $\varpi_i(h_{\gamma})=1$}\\
&\alpha-\gamma,\ \beta-\gamma\in R^+, \text{ if $\varpi_i(h_{\gamma})=0$}
\end{align*}

and in both equations \eqref{req4} and \eqref{req5} we additionally assume that 
$$m_{\alpha-\gamma}+m_{\gamma}>\ell \text{ if } \varpi_i(h_{\gamma})=0.$$
\end{defn}
\subsection{} 
Recall the definition of $\mathbf{R}_{\lambda,i}^{\ell}$ from \eqref{mainset} and note that 
$$\widetilde{\mathbf{R}}_{\lambda,i}^{\ell}:=\varpi_i-\mathbf{R}_{\lambda,i}^{\ell}=\left\{\mu\in \varpi_i-W(\varpi_i):  \varpi_i-\mu+\beta\notin W(\varpi_i), \ \forall\beta\in R^+ \text{ with } m_{\beta}=\ell\right\}$$
In Section~\ref{section63} we will associate to each $\mu\in \widetilde{\mathbf{R}}_{\lambda,i}^{\ell}\backslash\{0\}$ a tuple $\mathbf{o}(\mu)$ which we will call the \textit{orbit decomposition} of $\mu$ (see Definition~\ref{orbitdec9}). Defining
$$\mathbf{x}^\beta = \prod_{j=1}^n x_j^{b_j}\in\mathbb{C}[x_1,\dots,x_n],\ \ \beta= \sum_{j=1}^n b_j\alpha_j\in Q^+,$$ 
our main result reads as follows. 
\begin{thm}\label{mainthmres1} Let $\lambda\in P^+, \ell\in\mathbb{N}$ and $i\in\{1,\dots,i\}$.
\begin{enumerate}
\item The fusion product $\mathbf{D}_{\lambda}^{\ell}*V(\varpi_i)$ is independent of the fusion parameters and admits a filtration whose successive quotients are given by
$$\tau_{s_{\mathbf{o}(\mu)}}\mathbf{D}_{\lambda+\varpi_i-\mu}^{\ell},\ \ \mu \in \widetilde{\mathbf{R}}_{\lambda,i}^{\ell}$$
each appearing with multiplicity one. In particular, 
$$\mathrm{ch}_{\mathrm{gr}}\left(\mathbf{D}_{\lambda}^{\ell}*V(\varpi_i)\right)=\sum_{\mu\in \widetilde{\mathbf{R}}_{\lambda,i}^{\ell}} \mathrm{ch}_{\mathrm{gr}}(\mathbf{D}_{\lambda+\varpi_i-\mu}^{\ell})\ q^{s_{\mathbf{o}(\mu)}}.$$
\item We have an isomorphism of graded modules $$\mathbf{M}^{\ell}_{\lambda,i}\cong \mathbf{D}_{\lambda}^{\ell}*V(\varpi_i).$$
\item  The following recursion holds $$\mathbf{A}^{1\to \ell}_{\lambda+\varpi_i}(x_1,\dots,x_n) =\sum\limits_{\substack{\mu\in P^+\\ \lambda+\varpi_i-\mu\in\mathbf{R}_{\mu,i}^\ell} } \ \mathbf{x}^{\mu-\lambda} \mathbf{A}^{1\rightarrow \ell}_{\mu}(x_1,\dots,x_n).$$

 \end{enumerate}
\end{thm}
The proof of part (3) will be given in Section~\ref{section43} and in Section~\ref{section44} we provide an alternative proof for the tensor product decomposition for $\mathfrak{sl}_3$ using KR crystals. The proof of part (1) and (2) will be postponed to Section~\ref{section66}; we explain the idea behind the proof now. 

In a first step (see Lemma~\ref{smap}) we construct a surjective map 
\begin{equation}\label{730}\mathbf{M}^{\ell}_{\lambda,i}\rightarrow\mathbf{D}_{\lambda}^{\ell}*V(\varpi_i)\rightarrow 0\end{equation}
and a short exact sequence 
$$0\rightarrow \mathrm{Ker}(\psi)\rightarrow \mathbf{M}^\ell_{\lambda,i}\xrightarrow{\psi}\mathbf{D}^{\ell}_{\lambda+\varpi_i}\rightarrow 0.$$
Subsequently we show in Proposition~\ref{mainprophel} that $\mathrm{Ker}(\psi)$ admits a filtration whose successive quotients are homomorphic images of the level $\ell$ Demazure modules corresponding to elements in $\lambda+\mathbf{R}_{\lambda,i}^{\ell}\backslash\{\varpi_i\}$. Now using Theorem~\ref{crysesti1} we will derive in Section~\ref{section66} that the homomorphic maps have trivial kernels and \eqref{730} is an isomorphism.

\begin{rem}
\begin{enumerate} 
\item The Pieri type formula in part (1) can be used to express the $\lie h$-character of $\mathbf{D}_{\lambda}^{\ell}$ as a sum of products of ``smaller'' Demazure characters with elementary symmetric functions. It would be interesting to see whether this leads to a Jacobi-Trudi type identity for Demazure characters.
\item Part (3) of the above theorem reduces in the $\mathfrak{sl}_2$ case to \cite[Theorem 1.3]{CSVW16}.
\item Note that the description of $\widetilde{\mathbf{R}}_{\lambda,i}^{\ell}$ is type independent and we conjecture that Theorem~\ref{mainthmres1}(1) holds for all non-exceptional Lie algebras with appropriate \eqref{teilmitell} given in \cite[Section 3.2]{CV13}.
\end{enumerate}
\end{rem}

\begin{example}We discuss two examples.
\begin{enumerate}
    \item We consider \(n=3\) and the dominant weight \(\lambda=2\varpi_1+3\varpi_2+4\varpi_3\).
    Then Theorem~\ref{mainthmres1}(1) implies the following character equality
    \begin{align*}
       \hspace{1,5cm} \operatorname{ch}_{\mathrm{gr}}&(\mathbf{D}_{\lambda}^{5} * V(\varpi_2))&\\&= \operatorname{ch}_{\mathrm{gr}}(\mathbf{D}_{\lambda+\varpi_2}^{5})+q\operatorname{ch}_{\mathrm{gr}}(\mathbf{D}_{\lambda+\varpi_2-\alpha_2}^{5} )+q^2
        \operatorname{ch}_{\mathrm{gr}}(\mathbf{D}_{\lambda+\varpi_2-\alpha_{2,3}}^{5} ) +q^2\operatorname{ch}_{\mathrm{gr}}(\mathbf{D}_{\lambda+\varpi_2-\alpha_{1,3}}^{5} )
    \end{align*}
    In contrast to the classical Pieri formula, the terms of weight \(\lambda+\varpi_2-\alpha_{1,2}\) and \(\lambda+\varpi_2-\alpha_{1,3}-\alpha_2\) do not appear due to the fact that \(m_{1,2}=5\).
    \item  For \(n=5\) and the dominant weight \(\lambda=2\varpi_1+3\varpi_2+4\varpi_3+2\varpi_4+2\varpi_5\) we obtain
    \begin{align*}
     \hspace{1,5cm} \operatorname{ch}_{\mathrm{gr}}&(\mathbf{D}_{\lambda}^{4} * V(\varpi_1) )&\\&= \operatorname{ch}_{\mathrm{gr}}(\mathbf{D}_{\lambda+\varpi_1}^{4})+q^1\operatorname{ch}_{\mathrm{gr}}(\mathbf{D}_{\lambda+\varpi_1-\alpha_1}^{4} )+
        q^2\operatorname{ch}_{\mathrm{gr}}(\mathbf{D}_{\lambda+\varpi_1-\alpha_{1,2}}^{4} ) +q^3\operatorname{ch}_{\mathrm{gr}}(\mathbf{D}_{\lambda+\varpi_1-\alpha_{1,4}}^{4} )
        \end{align*}
    In this example the terms in which the formula differs from the classical formula are of weight \(\lambda+\varpi_1-\alpha_{1,3}\) and \(\lambda+\varpi_1-\alpha_{1,5}\).
    The term of weight \(\lambda+\varpi_1-\alpha_{1,3}\)  does not appear as \(m_{3}=4\) while the term of weight \(\lambda+\varpi_1-\alpha_{1,5}\) does not appear as \(m_{4,5}=4\)
\end{enumerate} 
\end{example}

\subsection{}\label{section43} In this subsection we prove Theorem~\ref{mainthmres1}(3). The following set of characters is linearly independent 
\[\{\operatorname{ch}_{\mathfrak{h}}(\mathbf{D}^\ell_\lambda) : \lambda \in P^+\}\]
for the following reason.
Given a linear dependence relation \(\sum_{\lambda\in P^+} a_\lambda \operatorname{ch}_{\mathfrak{h}}(\mathbf{D}^\ell_\lambda)=0\) we choose \(\mu\in P^+\) to be one of the maximal elements in \(\{\lambda\in P^+ : a_\lambda\neq0\}\) with respect to $\geq$. If $a_{\lambda}\neq 0$, we have \(\mu \not\leq\lambda\) and thus $(\mathbf{D}^\ell_\lambda)_\mu =0$. Therefore the basis element \(e_\mu\) has non-zero coefficient only in \(\operatorname{ch}_{\mathfrak{h}}(\mathbf{D}^\ell_\mu)\). Thus \(a_\mu =0\) for the linear dependence relation to hold which is a contradiction. The next result gives a recurrence relation for the numerical multiplicities; the proof follows closely \cite[Section 6.2]{BCSV15} and \cite[Section 6.4]{BCK18}.
\begin{lem}
    \label{lem:numerical_mult_rec}
Let $\lambda,\nu\in P^+$ and $\ell\geq k$. The numerical multiplicities are subject to the following relation
    \begin{align*}
        \sum_{\chi\in \mathbf{R}_{\lambda,i}^k}\left[\mathbf{D}^k_{\lambda+\chi}:\mathbf{D}^\ell_\nu  \right]_{q=1} = \sum\limits_{\substack{\mu\in P^+\\ \nu-\mu\in \mathbf{R}_{\mu,i}^\ell}} \left[\mathbf{D}^k_\lambda:\mathbf{D}^\ell_{\mu}  \right]_{q=1}.
    \end{align*}
\end{lem}
\begin{proof}
    Note that we can write with Theorem~\ref{flagc}
    \begin{align}
    \label{eq:character-flag}
       \operatorname{ch}_{\mathfrak{h}}(\mathbf{D}^k_\lambda) = \sum\limits_{\mu\in P^+} \left[\mathbf{D}^k_\lambda:\mathbf{D}^\ell_\mu  \right]_{q=1} \operatorname{ch}_{\mathfrak{h}}(\mathbf{D}^\ell_\mu).
    \end{align}
    Multiplying the left-hand side with \(\operatorname{ch}_\mathfrak{h}V(\varpi_i)\) and applying the Pieri formula yields 
    \begin{align*}
        \operatorname{ch}_{\mathfrak{h}}(\mathbf{D}^k_\lambda) \operatorname{ch}_{\mathfrak{h}}V(\varpi_i)= \sum\limits_{\chi\in \mathbf{R}_{\lambda,i}^k} \operatorname{ch}_{\mathfrak{h}}(\mathbf{D}^k_{\lambda+\chi}),
    \end{align*}
    which we can write once more with Theorem~\ref{flagc} in terms of level \(\ell\) Demazure modules as
    \begin{align*}
        \operatorname{ch}_{\mathfrak{h}}(\mathbf{D}^k_\lambda) \operatorname{ch}_{\mathfrak{h}}V(\varpi_i)= \sum\limits_{\chi\in \mathbf{R}_{\lambda,i}^k} \sum\limits_{\nu\in P^+} \left[\mathbf{D}^k_{\lambda+\chi}:\mathbf{D}^\ell_\nu  \right]_{q=1} \operatorname{ch}_{\mathfrak{h}}(\mathbf{D}^\ell_{\nu}).        
    \end{align*}
    The application of the Pieri formula to the right-hand side of equation \ref{eq:character-flag} yields 
    \begin{align*}
\sum\limits_{\mu\in P^+} \left[\mathbf{D}^k_\lambda:\mathbf{D}^\ell_\mu  \right]_{q=1} \operatorname{ch}_{\mathfrak{h}}(\mathbf{D}^\ell_\mu)\operatorname{ch}_{\mathfrak{h}}V(\varpi_i) = \sum\limits_{\mu\in P^+} \sum\limits_{\chi\in \mathbf{R}_{\mu,i}^\ell}\left[\mathbf{D}^k_\lambda:\mathbf{D}^\ell_\mu  \right]_{q=1} \operatorname{ch}_{\mathfrak{h}}(\mathbf{D}^\ell_{\mu+\chi})
    \end{align*}
   Combining the above equations gives
    \begin{align*}
   \sum_{\chi\in \mathbf{R}_{\lambda,i}^k} \sum\limits_{\nu\in P^+} \left[\mathbf{D}^k_{\lambda+\chi}:\mathbf{D}^\ell_\nu  \right]_{q=1} \operatorname{ch}_{\mathfrak{h}}(\mathbf{D}^\ell_{\nu})= \sum\limits_{\mu\in P^+} \sum\limits_{\chi\in \mathbf{R}_{\mu,i}^\ell}\left[\mathbf{D}^k_\lambda:\mathbf{D}^\ell_\mu  \right]_{q=1} \operatorname{ch}_{\mathfrak{h}}(\mathbf{D}^\ell_{\mu+\chi})
    \end{align*}
    Extracting coefficient in front of $\operatorname{ch}_{\mathfrak{h}}(\mathbf{D}^\ell_{\nu})$ together with the linear independence of the characters yields the relation.
\end{proof}
\begin{rem}
Specializing to \(k=1\) in Lemma~\ref{lem:numerical_mult_rec} simplifies the left hand side of the equation to $\left[\mathbf{D}^1_{\lambda+\varpi_i}:\mathbf{D}^\ell_\nu  \right]_{q=1}$ since \(R_{\lambda,i}^1=\{\varpi_i\}\).
\end{rem}
The next lemma proves Theorem~\ref{mainthmres1}(3).
\begin{lem} Let $\lambda\in P^+$, $\ell\in \mathbb{N}$ and $i\in\{1,\dots,n\}$. 
    For the generating series the following recurrence relation holds
    $$\mathbf{A}^{1\to \ell}_{\lambda+\varpi_i}(x_1,\dots,x_n) =\sum\limits_{\substack{\mu\in P^+\\ \lambda+\varpi_i-\mu\in\mathbf{R}_{\mu,i}^\ell} } \ \mathbf{x}^{\mu-\lambda} \mathbf{A}^{1\rightarrow \ell}_{\mu}(x_1,\dots,x_n).$$

\begin{proof}
 
We will suppress the dependence of the generating series on the variables \(x_1,\dots,x_n\) for readability. We get with Lemma~\ref{lem:numerical_mult_rec}
    \begin{align*}
        \mathbf{A}^{1\to \ell}_{\lambda+\varpi_i} &= \sum_{(k_1,\dots,k_n)\in\mathbb{Z}_+^n}\left[\mathbf{D}_{\lambda+\varpi_i+\sum_{j=1}^nk_j\alpha_j}^{1}:\mathbf{D}_{\lambda+\varpi_i}^{\ell}\right]_{q=1} \ x_1^{k_1}\cdots x_n^{k_n}
        \\
        &=\sum_{(k_1,\dots,k_n)\in\mathbb{Z}_+^n}
        \sum\limits_{\substack{\mu\in P^+\\ \lambda+\varpi_i-\mu\in\mathbf{R}_{\mu,i}^\ell}} \left[\mathbf{D}^1_{\lambda+\sum_{j=1}^nk_j\alpha_j}:\mathbf{D}^\ell_{\mu}  \right]_{q=1}\ x_1^{k_1}\cdots x_n^{k_n}
    \end{align*}
  Now writing $\mu-\lambda=\sum_{j=1}^n k^{\mu}_j\alpha_j$ we get
\begin{align*}
        \mathbf{A}^{1\to \ell}_{\lambda+\varpi_i}&=\sum_{(k_1,\dots,k_n)\in\mathbb{Z}_+^n}
        \sum\limits_{\substack{\mu\in P^+\\ \lambda+\varpi_i-\mu\in\mathbf{R}_{\mu,i}^\ell}} \left[\mathbf{D}^1_{\mu+\sum_{j=1}^n(k_j-k^{\mu}_j)\alpha_j}:\mathbf{D}^\ell_{\mu}  \right]_{q=1}\ x_1^{k_1}\cdots x_n^{k_n}
    \end{align*}
    Exchanging the summations and shifting the indices \(k_j\) by \(k^{\mu}_j\) we arrive at
    \begin{align*}
        \mathbf{A}^{1\to \ell}_{\lambda+\varpi_i} &=\sum\limits_{\substack{\mu\in P^+\\ \lambda+\varpi_i-\mu\in\mathbf{R}_{\mu,i}^\ell}} \sum_{(k_1,\dots,k_n)\in\mathbb{Z}_+^n}\left[\mathbf{D}^1_{\mu+\sum_{j=1}^nk_j\alpha_j}:\mathbf{D}^\ell_{\mu}  \right]_{q=1}\ x_1^{k_1+k^{\mu}_1}\cdots x_n^{k_n+k^{\mu}_n}&\\&
        =\sum\limits_{\substack{\mu\in P^+\\ \lambda+\varpi_i-\mu\in\mathbf{R}_{\mu,i}^\ell} } \ x_1^{k^{\mu}_1}\cdots x_n^{k^{\mu}_n}\ \mathbf{A}^{1\rightarrow \ell}_{\mu}.
    \end{align*}
\end{proof}
\end{lem}
\begin{example}
    For \(n=2,\ \ell\geq3\) and  \( \lambda\in P^+\) the recursion simplifies to
    \begin{align*}
        \mathbf{A}^{1\to \ell}_{\lambda+\varpi_1}(x_1,x_2)  = \mathbf{A}^{1\to \ell}_{\lambda} (x_1,x_2)&+ (1-\delta_{m_1,\ell-2}) x_1 \mathbf{A}^{1\to \ell}_{\lambda+\alpha_1}(x_1,x_2)&\\& + (1-\delta_{m_2,\ell-1})(1-\delta_{m_{1,2},\ell-2})x_1x_2\mathbf{A}^{1\to \ell}_{\lambda+\alpha_1+\alpha_2}(x_1,x_2) 
         \end{align*}
        \begin{align*}
        \mathbf{A}^{1\to \ell}_{\lambda+\varpi_2}(x_1,x_2)  = \mathbf{A}^{1\to \ell}_{\lambda}(x_1,x_2) &+ (1-\delta_{m_2,\ell-2}) x_2 \mathbf{A}^{1\to \ell}_{\lambda+\alpha_2}(x_1,x_2)&\\& + (1-\delta_{m_1,\ell-1})(1-\delta_{m_{1,2},\ell-2})x_1x_2\mathbf{A}^{1\to \ell}_{\lambda+\alpha_1+\alpha_2} (x_1,x_2).
    \end{align*}
\end{example}
\subsection{}\label{section44} We discuss first the case $\mathfrak{sl}_3$ where the proof is purely combinatorial. We hope that this approach can be lifted later to higher rank leading possibly to combinatorial models for higher level Demazure crystals. The Pieri formula simplifies to
$$\mathbf{D}_{a\varpi_1+b\varpi_2}^{\ell}\otimes V(\varpi_1)=\mathbf{D}_{(a+1)\varpi_1+b\varpi_2}^{\ell}+\zeta_1 \mathbf{D}_{(a-1)\varpi_1+(b+1)\varpi_2}^{\ell}+\delta_1 \mathbf{D}_{a\varpi_1+(b-1)\varpi_2}^{\ell}$$
where 
$$\zeta_1=\begin{cases} 0,& \text{if $a=0\mod \ell $}\\
1,& \text{otherwise }
\end{cases},\ \ \delta_1=\begin{cases} 0,& \text{if $b=0\mod \ell $ \ or \ $a+b= 0 \mod \ell$}\\
1,& \text{otherwise }
\end{cases}$$
and 
$$\mathbf{D}_{a\varpi_1+b\varpi_2}^{\ell}\otimes V(\varpi_2)=\mathbf{D}_{a\varpi_1+(b+1)\varpi_2}^{\ell}+\zeta_2 \mathbf{D}_{(a+1)\varpi_1+(b-1)\varpi_2}^{\ell}+\delta_2 \mathbf{D}_{(a-1)\varpi_1+b\varpi_2}^{\ell}$$
where 
$$\zeta_2=\begin{cases} 0,& \text{if $a=0\mod \ell $}\\
1,& \text{otherwise }
\end{cases},\ \ \delta_2=\begin{cases} 0,& \text{if $b=0\mod \ell $ \ or \ $a+b= 0 \mod \ell$}\\
1,& \text{otherwise }
\end{cases}$$
By the Steinberg type decomposition formula \eqref{stde} we can assume without loss of generality that $a,b< \ell$ and $a+b\geq \ell$. The case $a+b<\ell$ reduces to the usual Pieri formula stated in Section~\ref{section13}.

In the $\mathfrak{sl}_3$ case it seems easier to directly compute the classical decomposition of Demazure modules. We follow the strategy in \cite{LS19a} using the combinatorial models developed in \cite{K13a,K16a}. From \cite{LS19a} we know that the classical highest weight vectors of $\mathbf{D}_{a\varpi_1+b\varpi_2}^{\ell}$ can be obtained from the classical highest weight vectors in the tensor product of KR crystals $B^{1,a}\otimes B^{2,b}$ which are connected to the unique element in $B^{1,a}\otimes B^{2,b}$ of weight $a\varpi_1+b\varpi_2$ by using only arrows labeled by $1$ and $2$ or Demazure arrows. Recall that an arrow $b\rightarrow b'$ is called a Demazure arrow if the label of the arrow is $0$ and $\widetilde{\epsilon}_0(b)\geq \ell$. Using the realization mentioned in Section~\ref{section36} and \cite[Lemma 4.1]{K16a}, the classical highest weight vectors in $B^{1,a}\otimes B^{2,b}$ are of the form 

\begin{equation}\label{123}\begin{tikzpicture}
\node at (0,0)   (A){$\vline\tableau{ 0 \\ c}
\quad  \otimes \quad  {\vline\tableau{  0&0  }}\quad \quad \quad 0\leq c\leq \min\{a,b\}$};
\end{tikzpicture}
\end{equation}

If additionally $a+b-c\geq \ell$ (this means that $\widetilde{\epsilon}_0$ of the above element is $\geq \ell$) we can act with the Kashiwara operator $\tilde{f}_0$ (exactly $c$-times) to reach the unique element in $B^{1,a}\otimes B^{2,b}$ of weight $a\varpi_1+b\varpi_2$; this is the element with only zero entries. Conversely if $a+b-c<\ell$ we can never reach this element which we can see as follows. If we would be able to reach this element, there must be an element
\begin{equation}\label{123b}\begin{tikzpicture}
\node at (0,0)   (A){$\vline\tableau{ u \\ y}
\quad  \otimes \quad  {\vline\tableau{z&w }}$};
\end{tikzpicture}
\end{equation}
in the classical connected component of \eqref{123} where we can act with a Demazure arrow. However, by the definition of the Kashiwara operators $\tilde{f}_i$ with $i\neq 0$ we must have $u+y\geq c$ and thus
$$\ell> a+b-c\geq \max\{b-z-w,a+b-2z-w-u-y\}$$
Since the maximum above is $\widetilde{\epsilon}_0$ of the element \eqref{123b} we are done. So we get
$$\mathbf{D}_{a\varpi_1+b\varpi_2}^{\ell}=\bigoplus_{\substack{0\leq c\leq \min\{a,b\}\\ a+b-c\geq \ell}} V( (a-c)\varpi_1+(b-c)\varpi_2)$$
Now the Pieri formulas are immediate consequences of the above decomposition and the classical Pieri formula.

\section{The $(\lambda,i)$-Poset}\label{section5}

\subsection{} Here we define an order on $R_i^+$ which will be needed later to define a filtration of the kernel of a certain map. If $\alpha,\beta\in R^+$ and $\mathrm{supp}(\alpha)\cap \mathrm{supp}(\beta)\neq \emptyset$ we denote by $\alpha\cup \beta$ the unique positive root whose support is $\mathrm{supp}(\alpha)\cup \mathrm{supp}(\beta)$ and by $\alpha\cap \beta$ the unique positive root whose support is $\mathrm{supp}(\alpha)\cap \mathrm{supp}(\beta)$. For example if $\alpha-\beta\in R^+$ we have
$$\alpha\cup \beta=\beta,\ \ \alpha\cap \beta=\alpha.$$
Moreover, for $\alpha,\beta\in R_i^+$ with $\alpha\neq \beta$ and $\beta-\alpha\notin R$ we can write 
\begin{equation}\label{dec11}\alpha\cup\beta-\alpha\hspace{0,03cm}\cap\hspace{0,03cm} \beta=\gamma_1+\gamma_2\end{equation} for unique positive roots $\gamma_1,\gamma_2\in R^+$. To see this, one has only to observe that
the supporting intervals can have only an overlap as follows (or vice versa; $\alpha=\alpha_{r',p'}$ and $\beta=\alpha_{r,p}$)

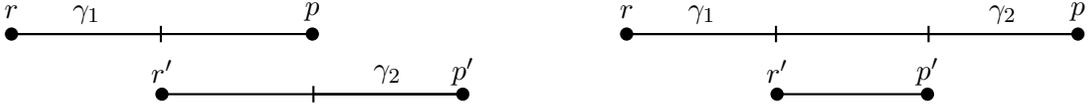
$$
\begin{tikzpicture}[thick]
\draw (0,0) -- (4, 0);
\node (A) at (0,0.3) [] {$r$};
\node (R) at (0,0) [] {};
\path [-|] (0,0) edge node [above] {$\gamma_1$} (2,0);
\path [|-] (4,-0.8) edge node [above] {$\gamma_2$} (6,-0.8);
\node (P) at (4,0) [] {};
\node (B) at (4,0.3) [circle] {$p$};
\fill (R) circle (2.5pt);
\fill (P) circle (2.5pt);
\draw (2,-0.8) -- (6, -0.8);
\node (C) at (2,-0.5) [circle] {$r'$};
\node (D) at (6,-0.5) [circle] {$p'$};
\node (RE) at (2,-0.8) [] {};
\node (PE) at (6,-0.8) [] {};
\fill (RE) circle (2.5pt);
\fill (PE) circle (2.5pt);
\end{tikzpicture}
\ \ \ \ \ \ \ \ \ \ \ \ 
\begin{tikzpicture}[thick]
\draw (-1,0) -- (5, 0);
\node (A) at (-1,0.3) [] {$r$};
\node (R) at (-1,0) [] {};
\path [-|] (-1,0) edge node [above] {$\gamma_1$} (1,0);
\path [|-] (3,0) edge node [above] {$\gamma_2$} (5,0);
\node (P) at (5,0) [] {};
\node (B) at (5,0.3) [circle] {$p$};
\fill (R) circle (2.5pt);
\fill (P) circle (2.5pt);
\draw (1,-0.8) -- (3, -0.8);
\node (C) at (1,-0.5) [circle] {$r'$};
\node (D) at (3,-0.5) [circle] {$p'$};
\node (RE) at (1,-0.8) [] {};
\node (PE) at (3,-0.8) [] {};
\fill (RE) circle (2.5pt);
\fill (PE) circle (2.5pt);
\end{tikzpicture}
$$
\textit{By convention we denote by $\gamma_1$ the left part of the string and by $\gamma_2$ the right part as indicated in the above picture.}

\subsection{} Recall the definition of $m_{\alpha},s_{\alpha}$ from Section~\ref{section13}. 
\begin{defn}\label{poset1} Assume that $\alpha,\beta\in R_{i}^+$. We say $\alpha \succeq \beta$ if $\alpha=\beta$ or either of the following conditions hold.
\begin{enumerate}[(i)]
    \item If $\alpha-\beta\in R$, then $$m_{\alpha\hspace{0,03cm}\cup\hspace{0,03cm} \beta}=m_{\alpha\hspace{0,03cm}\cap\hspace{0,03cm} \beta}+m_{\alpha\hspace{0,03cm}\cup\hspace{0,03cm} \beta-\alpha\hspace{0,03cm}\cap\hspace{0,03cm} \beta}-\ell \delta_{\alpha,\alpha\hspace{0,03cm}\cap\hspace{0,03cm} \beta}$$
    
    \item If $\alpha-\beta\notin R$ and $\alpha\cap\beta\in\{\alpha,\beta\}$, then 
    $$m_{\alpha\hspace{0,03cm}\cup\hspace{0,03cm} \beta}=m_{\alpha\hspace{0,03cm}\cap\hspace{0,03cm} \beta}+m_{\gamma_1}+m_{\gamma_2}-2\ell\delta_{\alpha,\alpha\hspace{0,03cm}\cap \hspace{0,03cm}\beta}$$
    where $\gamma_1,\gamma_2$ are as in \eqref{dec11}. 
    
    \vspace{0,1cm}
    
    \item  If $\alpha-\beta\notin R$ and $\alpha\cap\beta\notin\{\alpha,\beta\}$, then 
    $$m_{\alpha\hspace{0,03cm}\cup\hspace{0,03cm} \beta}=m_{\alpha}+m_{\alpha\hspace{0,03cm}\cup\hspace{0,03cm} \beta-\alpha}-\ell=m_{\beta}+m_{\alpha\hspace{0,03cm}\cup\hspace{0,03cm} \beta-\beta}$$
\end{enumerate}
\end{defn}
\begin{rem}\label{rem1}
\begin{enumerate}[(1)]
\item  It is worthwhile to mention that in our applications we will only compare roots using relation (i). Relations (ii) and (iii) are needed only to ensure transitivity of \(\succeq\).
\item It is straightforward to see that $\alpha\succeq \beta$ implies $m_{\alpha}\geq m_{\beta}$.
\end{enumerate}
\end{rem}
\begin{lem}\label{poslambdai}
We have that $(R_{i}^+,\succeq)$ is a partially ordered set which we call the $(\lambda,i)$-poset.
\begin{proof}
\underline{\textit{Reflexivity:}} This is clear by definition.

\underline{\textit{Antisymmetry:}} Let $\alpha,\beta\in R_{i}^+$ such that $\alpha\succeq \beta$, $\beta\succeq \alpha$ and $\alpha\neq \beta$. If $\alpha-\beta\in R$
we get 
$$m_{\alpha\hspace{0,03cm}\cup\hspace{0,03cm} \beta}=m_{\alpha\hspace{0,03cm}\cap\hspace{0,03cm} \beta}+m_{\alpha\hspace{0,03cm}\cup\hspace{0,03cm} \beta-\alpha\hspace{0,03cm}\cap\hspace{0,03cm} \beta}-\ell \delta_{\alpha,\alpha\hspace{0,03cm}\cap\hspace{0,03cm} \beta}=m_{\alpha\hspace{0,03cm}\cap\hspace{0,03cm} \beta}+m_{\alpha\hspace{0,03cm}\cup\hspace{0,03cm} \beta-\alpha\hspace{0,03cm}\cap\hspace{0,03cm} \beta}-\ell \delta_{\beta,\alpha\hspace{0,03cm}\cap\hspace{0,03cm} \beta}$$
which is a contradiction to $\delta_{\alpha,\alpha\hspace{0,03cm}\cap\hspace{0,03cm} \beta}\neq \delta_{\beta,\alpha\hspace{0,03cm}\cap\hspace{0,03cm} \beta}$. Similarly we get a contradiction if we are in the situation of (ii) or (iii). 

\underline{\textit{Transitivity:}} The proof of the transitivity considers many cases and is quite a long calculation. We postpone the proof to the appendix (see Section~\ref{appendix}). 
\end{proof}
\end{lem}

\section{Fusion product with fundamental modules}\label{section6}
\subsection{} In this section we aim to determine the structure of the fusion product of an arbitrary Demazure module with a fundamental module. Recall that
$$\mathbf{D}_{\lambda}^{\ell}* V(\varpi_i)=\bigoplus_{d\geq 0} V^d/V^{d-1},\ \ V^{-1}=0,\ V^d=\sum_{s=0}^{d}\mathbf{U}_s(v_1\otimes v_2)$$ with distinct fusion parameters $(z_1,z_2)\in \mathbb{C}\times \mathbb{C}$. We remark that the case $\ell=1$ is trivial and the difficulties appear for higher level prime Demazure modules. We keep the notation established in Section~\ref{mainresultsection}.

\begin{lem}\label{smap}
We have a surjective map $$\mathbf{M}^{\ell}_{\lambda,i}\rightarrow\mathbf{D}_{\lambda}^{\ell}*V(\varpi_i)\rightarrow 0$$
\begin{proof}
We have to check step-wise that the elements listed in \eqref{req1}-\eqref{req5} act trivially on the cyclic generator $v_1*v_2$ of the fusion product $\mathbf{D}_{\lambda}^{\ell}*V(\varpi_i)$. Recall that $V(\varpi_i)=\Lambda^i\mathbb{C}^{n+1}$ with highest weight vector $v_2=e_1\wedge\cdots \wedge e_{i}$. If $s_{\boldsymbol{\alpha}}=s_{\kappa(\boldsymbol{\alpha})}$, we have in $V_{s_{\boldsymbol{\alpha}}}/V_{s_{\boldsymbol{\alpha}-1}}$
\begin{align*}X_{\boldsymbol{\alpha}}(v_1\otimes v_2)&=(x^-_{r_{1},p_1}\otimes t^{s_{r_{1},p_{1}}})\cdots (x^-_{r_{k},p_k}\otimes t^{s_{r_{k},p_{k}}})(v_1\otimes v_2)&\\&
=(x^-_{r_{1},p_1}\otimes (t-z_1)^{s_{r_{1},p_{1}}})\cdots (x^-_{r_{k},p_k}\otimes (t-z_1)^{s_{r_{k},p_{k}}})(v_1\otimes v_2)&\\&
= (z_2-z_1)^{s_{\boldsymbol{\alpha}}}\left (v_1\otimes x^-_{r_{1},p_1}\cdots x^-_{r_{k},p_k}v_2\right)&\\&
= \mathrm{sgn}(\kappa)(z_2-z_1)^{s_{\sigma(\boldsymbol{\alpha})}}\left (v_1\otimes x^-_{r_{\kappa(1)},p_1}\cdots x^-_{r_{\kappa(k)},p_k}v_2\right)&\\&=\mathrm{sgn}(\kappa)X_{\kappa(\boldsymbol{\alpha})}(v_1\otimes v_2)\end{align*}
If $s_{\boldsymbol{\alpha}}\neq s_{\boldsymbol{\alpha}}^{\mathrm{min}}$ there exists $\rho\in \Sigma_k$ such that $s_{\boldsymbol{\alpha}}>s_{\rho(\boldsymbol{\alpha})}$. In particular, the same calculation as above shows in $V_{s_{\boldsymbol{\alpha}}}/V_{s_{\boldsymbol{\alpha}-1}}$:
$$X_{\boldsymbol{\alpha}}(v_1\otimes v_2)=\mathrm{sgn}(\rho)(z_2-z_1)^{s_{\boldsymbol{\alpha}}-s_{\rho(\boldsymbol{\alpha})}}X_{\rho(\boldsymbol{\alpha})}(v_1\otimes v_2)=0$$
This proves that \eqref{req1} holds in the fusion product. The relations in \eqref{req2} and \eqref{req3} are immediate from the defining relations of $\mathbf{D}_{\lambda}^{\ell}$ and $V(\varpi_i)$. For example, to see \eqref{req3} we simply replace the variable $t$ by $(t-z_1)$ and act only on the second factor in the fusion product. Noting that $x_{\alpha}^-x_{\beta}^-v_2=0$ shows this part. Next we note that $\alpha-\gamma\in R^+$ and $m_{\alpha-\gamma}+m_{\gamma}>\ell$ implies
\begin{equation}\label{rfgac}0=\left(x^-_{\alpha-\gamma}\otimes t^{s_{\alpha-\gamma}}\right)\left(x^-_{\gamma}\otimes t^{s_{\gamma}-1}\right)^{m_{\gamma}+1}v_1=\left(x^-_{\gamma}\otimes t^{s_{\gamma}-1}\right)^{m_{\gamma}}\left(x_{\alpha}^-\otimes t^{s_{\alpha}-1}\right)v_1\end{equation} 
Thus,  if $\varpi_i(h_{\gamma})=0$, we get modulo terms of lower degree and \eqref{rfgac}
\begin{align*}&\left(x^-_{\gamma}\otimes t^{s_{\gamma}-1}\right)^{m_{\gamma}}\left(x_{\alpha}^-\otimes t^{s_{\alpha}}\right)(v_1\otimes v_2)&\\&=\left(x^-_{\gamma}\otimes t^{s_{\gamma}-1}\right)^{m_{\gamma}}\left(x_{\alpha}^-\otimes t^{s_{\alpha}}\right)(v_1\otimes v_2)-(z_2-z_1)\left(x^-_{\gamma}\otimes t^{s_{\gamma}-1}\right)^{m_{\gamma}}\left(x_{\alpha}^-\otimes t^{s_{\alpha}-1}\right)(v_1\otimes v_2)&\\&
=(z_2-z_1)^{s_{\alpha}}\left(
\left(x^-_{\gamma}\otimes t^{s_{\gamma}-1}\right)^{m_{\gamma}}v_1\otimes x_{\alpha}^-v_2\right)-(z_2-z_1)^{s_{\alpha}}\left(
\left(x^-_{\gamma}\otimes t^{s_{\gamma}-1}\right)^{m_{\gamma}}v_1\otimes x_{\alpha}^-v_2\right)&\\&
\hspace{3cm}-(z_2-z_1)\left(\left(x^-_{\gamma}\otimes t^{s_{\gamma}-1}\right)^{m_{\gamma}}\left(x_{\alpha}^-\otimes t^{s_{\alpha}-1}\right)v_1\otimes v_2\right)=0\end{align*}
and \eqref{req4} is obtained in this case. For \eqref{req4} and $\varpi_i(h_{\gamma})=1$ we observe that we have modulo terms of lower degree
$$\left(x^-_{\gamma}\otimes t^{s_{\gamma}-1}\right)^{m_{\gamma}+1}\left(x_{\alpha}^-\otimes t^{s_{\alpha}}\right)(v_1\otimes v_2)=(z_2-z_1)^{s_{\alpha}}\left(\left(x^-_{\gamma}\otimes t^{s_{\gamma}-1}\right)^{m_{\gamma}+1}v_1\otimes x_{\alpha}^-v_2\right)=0$$
since $x_{\gamma}^-x_{\alpha}^-v_2=0$ and we are done. Now we consider \eqref{req5} with $\varpi_i(h_{\gamma})$=1; the case $\varpi_i(h_{\gamma})=0$ proceeds similarly as the proof of relation \eqref{req4} and we omit the details.  Note that the roots are in the following position:

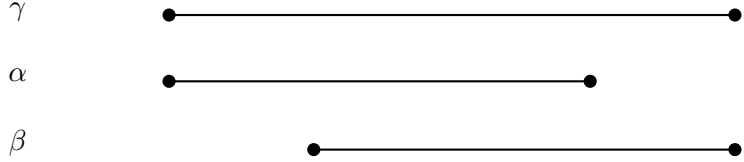
\begin{center}
\begin{tikzpicture}[baseline=1em]
  \matrix (m) [matrix of math nodes,row sep=1em,column sep=3em,minimum width=2em] {
    \gamma& \mbox{} &  \mbox{} & \mbox{} &  \mbox{} & \mbox{} &  \mbox{} \\
    \alpha& \mbox{} &  \mbox{} & \mbox{} &  \mbox{} & \mbox{} &  \mbox{} \\
  \beta & \mbox{} & \mbox{} & \mbox{} &  \mbox{} & \mbox{} &  \mbox{}   \\};
  \path[Circle-Circle,thick]
  (m-1-2.center) edge  (m-1-6.center)
  (m-2-2.center) edge  (m-2-5.center)
  (m-3-3.center) edge  (m-3-6.center)
  ;
\end{tikzpicture}
\end{center}

We get modulo terms of lower degree
\begin{align}\label{rfgac1}\left(x^-_{\gamma}\otimes t^{s_{\gamma}-1}\right)^{m_{\gamma}}&\left(x^-_{\alpha}\otimes t^{s_{\alpha}}\right)\left(x^-_{\beta}\otimes t^{s_{\beta}}\right)(v_1\otimes v_2)&\\&\notag=(z_2-z_1)^{s_{\alpha}+s_{\beta}}\left(\left(x^-_{\gamma}\otimes t^{s_{\gamma}-1}\right)^{m_{\gamma}}v_1\otimes x_{\alpha}^-x_{\beta}^-v_2\right)\end{align}
If $m_{\gamma-\alpha}+m_{\alpha}\leq \ell$ or $m_{\gamma-\beta}+m_{\beta}\leq \ell$ we obtain with \eqref{req1} that 
\[\left(x^-_{\alpha}\otimes t^{s_{\alpha}}\right)\left(x^-_{\beta}\otimes t^{s_{\beta}}\right)(v_1\otimes v_2)\approx\left(x^-_{\alpha\hspace{0.03cm}\cap\hspace{0.03cm} \beta}\otimes t^{s_{\alpha\hspace{0.03cm}\cap \hspace{0.03cm}\beta}}\right)\left(x^-_{\gamma}\otimes t^{s_{\gamma}}\right)(v_1\otimes v_2)\]
are proportional and we are done since
$$\left(x_{\gamma}^-\otimes t^{s_{\gamma}-1}\right)^{m_{\gamma}}\left(x_{\gamma}^-\otimes t^{s_{\gamma}}\right)(v_1*v_2)=0$$
Similarly the elements are also proportional if either 
$m_{\alpha\hspace{0.03cm}\cap \hspace{0.03cm}\beta}+m_{\gamma-\alpha}>\ell$ or $m_{\alpha\hspace{0.03cm}\cap \hspace{0.03cm}\beta}+m_{\gamma-\beta}>\ell$. So we can assume in the rest of the proof that
$$m_{\gamma-\alpha}+m_{\alpha}>\ell,\ \ m_{\gamma-\beta}+m_{\beta}> \ell$$
$$m_{\alpha\hspace{0.03cm}\cap \hspace{0.03cm}\beta}+m_{\gamma-\alpha}\leq\ell,\ \ m_{\alpha\hspace{0.03cm}\cap \hspace{0.03cm}\beta}+m_{\gamma-\beta}\leq \ell.$$
Thus, modulo terms of lower degree we get 
\begin{align*}&\left(x^-_{\gamma}\otimes t^{s_{\gamma}-1}\right)^{m_{\gamma}}\left(x^-_{\alpha}\otimes t^{s_{\alpha}}\right)\left(x^-_{\beta}\otimes t^{s_{\beta}}\right)(v_1\otimes v_2)&\\&=\left(x^-_{\gamma}\otimes t^{s_{\gamma}-1}\right)^{m_{\gamma}}\left(x^-_{\alpha}\otimes t^{s_{\alpha}}\right)\left(x^-_{\gamma-\alpha}\otimes t^{s_{\gamma-\alpha}-1}\right)\left(x^-_{\alpha\hspace{0.03cm}\cap \hspace{0.03cm}\beta}\otimes t^{s_{\alpha\hspace{0.03cm}\cap \hspace{0.03cm}\beta}}\right)(v_1\otimes v_2)&\\&
\hspace{3cm}-\left(x^-_{\gamma}\otimes t^{s_{\gamma}-1}\right)^{m_{\gamma}}\left(x^-_{\alpha}\otimes t^{s_{\alpha}}\right)\left(x^-_{\alpha\hspace{0.03cm}\cap \hspace{0.03cm}\beta}\otimes t^{s_{\alpha\hspace{0.03cm}\cap \hspace{0.03cm}\beta}}\right)\left(x^-_{\gamma-\alpha}\otimes t^{s_{\gamma-\alpha}-1}\right)(v_1\otimes v_2)&\\&
=-\left(x^-_{\gamma}\otimes t^{s_{\gamma}-1}\right)^{m_{\gamma}+1}\left(x^-_{\alpha\hspace{0.03cm}\cap \hspace{0.03cm}\beta}\otimes t^{s_{\alpha\hspace{0.03cm}\cap \hspace{0.03cm}\beta}}\right)(v_1\otimes v_2)&\\&
\hspace{3cm}-\left(x^-_{\gamma}\otimes t^{s_{\gamma}-1}\right)^{m_{\gamma}}\left(x^-_{\alpha}\otimes t^{s_{\alpha}}\right)\left(x^-_{\alpha\hspace{0.03cm}\cap \hspace{0.03cm}\beta}\otimes t^{s_{\alpha\hspace{0.03cm}\cap \hspace{0.03cm}\beta}}\right)\left(x^-_{\gamma-\alpha}\otimes t^{s_{\gamma-\alpha}-1}\right)(v_1\otimes v_2)&\\&
=-(z_2-z_1)^{s_{\alpha\hspace{0.03cm}\cap \hspace{0.03cm}\beta}+s_{\gamma}-1}\left(\left(x^-_{\gamma}\otimes t^{s_{\gamma}-1}\right)^{m_{\gamma}}v_1\otimes  x_{\gamma}^-x_{\alpha\hspace{0.03cm}\cap \hspace{0.03cm}\beta}^-v_2\right)&\\&
\hspace{3cm}-(z_2-z_1)^{s_{\alpha\hspace{0.03cm}\cap \hspace{0.03cm}\beta}}\left(x^-_{\gamma}\otimes t^{s_{\gamma}-1}\right)^{m_{\gamma}}\left(x^-_{\alpha}\otimes t^{s_{\alpha}}\right)\left(\left(x^-_{\gamma-\alpha}\otimes t^{s_{\gamma-\alpha}-1}\right)v_1\otimes x_{\alpha\hspace{0.03cm}\cap \hspace{0.03cm}\beta}^{-}v_2\right)&\\&
=-(z_2-z_1)^{s_{\alpha\hspace{0.03cm}\cap \hspace{0.03cm}\beta}+s_{\gamma}-1}\left(\left(x^-_{\gamma}\otimes t^{s_{\gamma}-1}\right)^{m_{\gamma}}v_1\otimes  x_{\gamma}^-x_{\alpha\hspace{0.03cm}\cap \hspace{0.03cm}\beta}^-v_2\right)&\\&
\hspace{3cm}+(z_2-z_1)^{s_{\alpha\hspace{0.03cm}\cap \hspace{0.03cm}\beta}}\left(x^-_{\gamma}\otimes t^{s_{\gamma}-1}\right)^{m_{\gamma}}\left(\left(x^-_{\gamma}\otimes t^{s_{\gamma}-1}\right)v_1\otimes x_{\alpha\hspace{0.03cm}\cap \hspace{0.03cm}\beta}^{-}v_2\right)
&\\&=-(z_2-z_1)^{s_{\alpha\hspace{0.03cm}\cap \hspace{0.03cm}\beta}+s_{\gamma}-1}\left(\left(x^-_{\gamma}\otimes t^{s_{\gamma}-1}\right)^{m_{\gamma}}v_1\otimes  x_{\gamma}^-x_{\alpha\hspace{0.03cm}\cap \hspace{0.03cm}\beta}^-v_2\right)&\\&
\hspace{3cm}+(z_2-z_1)^{s_{\alpha\hspace{0.03cm}\cap \hspace{0.03cm}\beta}+s_{\gamma}-1}\left(\left(x^-_{\gamma}\otimes t^{s_{\gamma}-1}\right)^{m_{\gamma}}v_1\otimes x_{\gamma}^-x_{\alpha\hspace{0.03cm}\cap \hspace{0.03cm}\beta}^{-}v_2\right)=0.
\end{align*}
\end{proof}
\end{lem}
Setting $$\mathfrak{g}^+[t]=(\lie b\otimes 1) \oplus (\lie g\otimes t \mathbb{C}[t]),\ \ \ \widetilde{\mathbf{M}}^{\ell}_{\lambda,i}:= \mathbf{U}(\mathfrak{g}^+[t])m\subseteq \mathbf{M}^{\ell}_{\lambda,i}$$ we can derive that $\widetilde{\mathbf{M}}^{\ell}_{\lambda,i}$ is a quotient of a module which has been introduced in \cite[Definition 3.2]{KV21} (in the notation of \cite{KV21} it is the module $M''_{\lambda+\varpi_i,\mathbf{p}}$). The only consequence we want to emphasize is that all relations in $M''_{\lambda+\varpi_i,\mathbf{p}}$ also hold in $\widetilde{\mathbf{M}}^{\ell}_{\lambda,i}$ and therefore in  $\mathbf{M}^{\ell}_{\lambda,i}$. We describe these relations next. 
For $r\in\mathbb{N}$ and $s,k\in\mathbb{Z}_+$, let \begin{equation*}_k\mathbf{S}(r,s)=\left\{(b_p)_{p\ge 0}: b_p\in\bz_+, \ \  \sum_{p\ge 0} b_p=r,\ \ \sum_{p\ge 0} pb_p=s,\ \ b_j=0\ \ j<k\right\}\end{equation*}
\begin{align*} {}^{t}_{k}\mathbf{x}_{\alpha}(r,s)&=\sum_{\mathbf{b}\in _k\mathbf{S}(r,s)}(x\otimes t^{k+1})^{(b_k)}\cdots (x\otimes t^{s+1})^{(b_s)}
\end{align*}
where for any integer $p$ and any $y\in\lie g[t]$,  we set $y^{(p)}=y^p/p!$. The following is an immediate consequence of the discussion above and \cite[Proposition 3.5 and Theorem 2]{KV21}.
\begin{prop}\label{prop1a}
For all $\alpha\in R^+$ and  $r\in\mathbb{N}$ we have $${}^{t}_{k}\mathbf{x}^-_\alpha(r,s)m=0,\ \forall s, k\in \mathbb{Z}_+: \ s+r\ge 1+ rk+\max\{0,\lambda(h_{\alpha})-\ell (k+1)\}+\varpi_i(h_{\alpha}).$$
\qed
\begin{example}\label{ex11}
If $s_{\alpha}\geq 1$ and $\varpi_i(h_{\alpha})=1$ we obtain the relation 
$$(x^-_{\alpha}\otimes t^{s_{\alpha}-1})^{m_{\alpha}}(x^-_{\alpha}\otimes t^{s_{\alpha}})m=0$$
in $\mathbf{M}^{\ell}_{\lambda,i}$. For $s_{\alpha}=1$ this follows from $\mathfrak{sl}_2$-theory. Otherwise this follows from Proposition~\ref{prop1a} if we set $s=(s_{\alpha}-2)m_{\alpha}+s_{\alpha}-1$, $r=m_{\alpha}+1$ and $k=s_{\alpha}-2$. Then 
$$s+r=s_{\alpha}m_{\alpha}+s_{\alpha}-m_{\alpha}\ge 1+ rk+\max\{0,m_{\alpha}\}+\varpi_i(h_{\alpha})=s_{\alpha}m_{\alpha}-m_{\alpha}+s_{\alpha}$$
and hence 
$${}^{t}_{k}\mathbf{x}^-_\alpha(r,s)m=
(x^-_{\alpha}\otimes t^{s_{\alpha}-1})^{(m_{\alpha})}(x^-_{\alpha}\otimes t^{s_{\alpha}})m=0.$$
\end{example}
\hfill\qed
\end{prop}
\subsection{} \label{section42}
From the defining relations of $\mathbf{M}^{\ell}_{\lambda,i}$ stated in \eqref{req1}-\eqref{req5} it is straightforward to check that we have a short exact sequence
$$0\rightarrow \mathrm{Ker}(\psi)\rightarrow \mathbf{M}^{\ell}_{\lambda,i}\xrightarrow{\psi} \mathbf{D}^{\ell}_{\lambda+\varpi_i}\rightarrow 0$$
where the kernel $K:=\mathrm{Ker}(\psi)$ is generated by the elements 
$$\left\{(x^-_{\alpha}\otimes t^{s_{\alpha}})m: \varpi_i(h_{\alpha})=1,\ \ m_{\alpha}<\ell\right\}.$$
However the number of generators of the kernel can be reduced by the following observation.  Assume that $\alpha_{r,p}\in R_i^+$ and $m_{j,p}=\ell$ for some $j\in\{i+1,\dots,p\}$.
Then we have $s_{r,p}=s_{r,j-1}+s_{j,p}$ and thus
$$(x_{r,p}^-\otimes t^{s_{r,p}})m\in \mathbf{U}\cdot (x_{{r,j-1}}^-\otimes t^{s_{r,j-1}})m$$
Similarly, we get $$(x_{r,p}^-\otimes t^{s_{r,p}})m\in \mathbf{U}\cdot (x_{{j+1,p}}^-\otimes t^{s_{j+1,p}})m$$ provided that $m_{r,j}=\ell$ for some $j\in\{r,\dots,i-1\}$.
So we can derive the following from the above calculations. 
\begin{prop}
\label{prop1ref}
The Kernel $K$ of the map $\mathbf{M}^{\ell}_{\lambda,i}\twoheadrightarrow \mathbf{D}^{\ell}_{\lambda+\varpi_i}$ is generated by the elements 
\begin{equation}\label{generatorsK1}\left\{(x^-_{\alpha}\otimes t^{s_{\alpha}})m: \alpha\in R^+\cap \widetilde{\mathbf{R}}_{\lambda,i}^{\ell} \right\}.\end{equation}
\qed
\end{prop}
\subsection{}\label{section63} The symmetric group $\Sigma_k$ acts on $\left(R_{i}^+\right)^{\times k}$ by permuting  the entries. We extend the $(\lambda,i)$-poset structure on $R_{i}^+$ to a total order and consider the induced order on the set $$\displaystyle\bigcup_{k\geq 1}^{}\left(R_{i}^+\right)^{\times k}/\Sigma_k$$
defined as follows. Given two tuples 
$$\boldsymbol{\beta}=(\beta_1,\dots,\beta_{k_1})\in \left(R_{i}^+\right)^{\times k_1},\  \beta_1\preceq \cdots \preceq \beta_{k_1},\ \  \boldsymbol{\beta}'=(\beta'_1,\dots,\beta'_{k_2})\in \left(R_{i}^+\right)^{\times k_2},\  \beta'_1\preceq \cdots \preceq \beta'_{k_2}$$ we say that $\boldsymbol{\beta}\prec \boldsymbol{\beta}'$ if and only if there exists $j\in\{1,\dots,\min\{k_1,k_2\}\}$ with 
$$ \beta_1=\beta_1',\dots,\beta_j=\beta_j',\ \beta_{j+1}\prec \beta_{j+1}'$$
or $\beta_1=\beta_1',\dots, \beta_{k_2}=\beta_{k_2}'$ and $k_1>k_2$.
Now we want to identify any non-zero element $\mu\in \widetilde{\mathbf{R}}_{\lambda,i}^{\ell}$ as a tuple of positive roots as follows. First we remark that each non-zero $\mu\in \widetilde{\mathbf{R}}_{\lambda,i}^{\ell}$ can be written uniquely as 
$$\mu=\alpha_{r_{1},p_1}+\cdots+\alpha_{r_{k},p_k},\ \ k=(\varpi_i,\mu),\ \ \mathrm{supp}(\alpha_{r_1,p_1})\subsetneq \cdots \subsetneq \mathrm{supp}(\alpha_{r_k,p_k}).$$
Consider the set
$$N_{\mu}=\left\{(\alpha_{r_{\kappa(1)},p_1},\dots,\alpha_{r_{\kappa(k)},p_k}): \kappa\in \tilde{\Sigma}_k\right\}$$
where $\tilde{\Sigma}_k\subseteq \Sigma_k$ denotes the subset of elements where the following sum 
$$s_{r_{\kappa(1)},p_1}+\cdots +s_{r_{\kappa(k)},p_k}$$
takes the minimal value. 
\begin{defn}\label{orbitdec9} Given a non-zero element $\mu\in \widetilde{\mathbf{R}}_{\lambda,i}^{\ell}$ we define the \textit{orbit decomposition} of $\mu$ to be the (up to permutation) unique minimal tuple
$\mathbf{o}(\mu)=(\mu_1,\dots,\mu_k)$ in $N_{\mu}$ with respect to the above defined order $\preceq$. We extend this definition by $\mathbf{o}(0)=\emptyset$.
\end{defn}

Our aim is to define a filtration on $K$ and study the associated graded space. Note that  $\preceq$ also defines a total order on the set of monomials
$$X_{\boldsymbol{\alpha}},\ \ \boldsymbol{\alpha}\in \left(R_{i}^+\right)^{\times k}$$
in the obvious way. We introduce 
$$\mathrm{HW}_i:=\left\{X_{\bold{o}(\mu)}m: \mu\in \widetilde{\mathbf{R}}_{\lambda,i}^{\ell}\backslash\{0\}\right\}$$
which shall represent the highest weight vectors in the Pieri formula corresponding to $K$ up to some filtration which we define now. 
We define a filtration of the kernel $K$ as follows. We order the representatives of the elements in 
$\mathrm{HW}_i$ as
\begin{equation}\label{orddef}\bold{o}(\mu)_1\prec \cdots \prec \bold{o}(\mu)_s\end{equation}
and define
$$0\subseteq K_1\subseteq \cdots \subseteq K_{s-1}\subseteq K_s=K,\ \ K_j=\sum_{r=1}^j \mathbf{U}\cdot  X_{\boldsymbol{o}(\mu)_{r}}m$$
Hence each $K_j/K_{j-1}$ is a cyclic $\mathbf{U}$-module with cyclic generator $X_{\bold{o}(\mu)_j}m\in \mathrm{HW}_i$.
\subsection{}
Given a tuple $\boldsymbol{\alpha}=(\alpha_{r_1,s_1},\dots,\alpha_{r_k,s_k})\in \left(R_i^+\right)^k$ we denote by $X_{\boldsymbol{\alpha}}^{i_1,\dots,i_d}$
the element obtained from $X_{\boldsymbol{\alpha}}$ by erasing the vectors
$$(x^-_{r_{i_1},p_{i_1}}\otimes t^{s_{r_{i_1},p_{i_1}}}),\dots, (x^-_{r_{i_d},p_{i_d}}\otimes t^{s_{r_{i_d},p_{i_d}}})$$
The following lemma is crucial. 
\begin{lem}\label{nplus}
Let $\mathbf{o}(\mu)=(\mu_1,\dots,\mu_k)$ be the orbit decomposition of $\mu\in \widetilde{\mathbf{R}}_{\lambda,i}^{\ell}\backslash\{0\}$, and $\beta\in R^+$ such that $\varpi_i-(\mu-\beta)\in W(\varpi_i)$ and $\eta:=\mu_w-\beta\in R^+$ for some $w\in\{1,\dots,k\}$. We suppose also that the tuple
$$\boldsymbol{\mu}_{\beta}=\begin{cases}
(\mu_1,\dots,\mu_w-\beta,\dots,\mu_{k}),& \text{ if $\varpi_i(h_{\beta})=0$}\\
(\mu_{1},\dots,\mu_p+\eta,\dots,\mu_{k})\backslash\{\mu_w\}, & \text{ if $\varpi_i(h_{\beta})=1$} 
\end{cases}$$
satisfies $s_{\boldsymbol{\mu}_{\beta}}=s^{\mathrm{min}}_{\boldsymbol{\mu}_{\beta}}$ where $p\in\{1,\dots,k\}$ is the unique index with $\mu_p+\eta\in R^+$.
Moreover, we assume that
$$m_{\beta}+m_{\eta}\leq \ell ,\ \ m_{\eta}+m_{\mu_{p}}>\ell.$$ 
Then $\mu-\beta\in \widetilde{\mathbf{R}}_{\lambda,i}^{\ell}\backslash\{0\}$ or there exists $\rho\in \widetilde{\mathbf{R}}_{\lambda,i}^{\ell}\backslash\{0\}$ with the following property

$$X_{\boldsymbol{\mu}_{\beta}}m\in \mathbf{U}\cdot X_{\mathbf{o}(\rho)}m,\ \ X_{\mathbf{o}(\rho)}\prec X_{\mathbf{o}(\mu)}.$$
    
\end{lem}
\begin{proof}
Assume that $\mu-\beta\notin \widetilde{\mathbf{R}}_{\lambda,i}^{\ell}$, i.e.  
there exists $\beta_1\in R^+$ such that $m_{\beta_1}=\ell$ and \begin{equation}\label{p11}\varpi_i-\mu+\beta+\beta_1\in W(\varpi_i)\end{equation}
and we choose $\beta_1$ of maximal height.  
So we have the situation $$\varpi_i-\mu,\varpi_i-\mu+\beta, \varpi_i-\mu+\beta+\beta_1\in W(\varpi_i),\ \varpi_i-\mu+\beta_1\notin W(\varpi_i).$$
Writing $\varpi_i-\mu=\epsilon_{j_1}+\cdots+\epsilon_{j_i}$, $\beta=\epsilon_r-\epsilon_s$ and $\beta_1=\epsilon_{r_1}-\epsilon_{s_1}$ we have
\begin{itemize}
    \item $s_1\notin \{j_1,\dots,j_i\}$ or $r_1\in\{j_1,\dots,j_i\}$
    \vspace{0,2cm}
    
    \item $s_1\in \{r,j_1,\dots,j_i\}\backslash\{s\}$,\ $r_1\notin\{r,j_1,\dots,j_i\}\backslash\{s\}$
\end{itemize}
which implies $r_1=s$ or $r=s_1$ and hence $\beta+\beta_1\in R^+$. 
The above calculations imply that we must have locally one of the following situations, where the cases listed below can also appear simultaneously

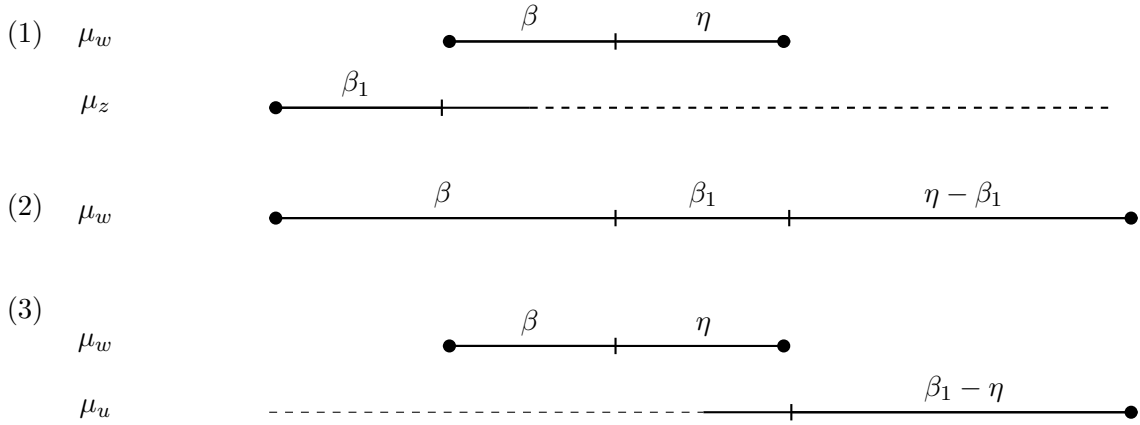
\begin{enumerate}
\item  
\begin{tikzpicture}[baseline=1em]
  \matrix (m) [matrix of math nodes,row sep=1em,column sep=4em,minimum width=2em] {
    \mu_w& \mbox{} &  \mbox{} & \mbox{} &  \mbox{} & \mbox{} &  \mbox{} \\
  \mu_z & \mbox{} & \mbox{} & \mbox{} &  \mbox{} & \mbox{} &  \mbox{}   \\};
  
  \node (a) at ($(m-2-3)!0.5!(m-2-4)$) {};
  
  \path[Circle-Circle,thick]
  (m-1-3.center) edge  (m-1-5.center)
  ;
  \path[Circle-,thick]
  (m-2-2.center) edge (a.center) 
  ;
  \path[-,thick, dashed]
  (a.center) edge (m-2-7) 
  ;
  
  \path[-|,thick]
  (m-1-3.center) edge node [above] {\(\beta\)} (m-1-4.center) 
  (m-2-2.center) edge node [above] {\(\beta_1\)} (m-2-3.center); 
  \path[-,thick]
  (m-1-4.center) edge node [above] {\(\eta\)} (m-1-5.center); 
\end{tikzpicture}
\\ \\
  \item  \begin{tikzpicture}[baseline=0em]
  \matrix (m) [matrix of math nodes,row sep=1em,column sep=4em,minimum width=2em] {
  \mu_w & \mbox{} & \mbox{} & \mbox{} &  \mbox{} & \mbox{} &  \mbox{}   \\};
  \path[Circle-Circle,thick]
  (m-1-2.center) edge node [above] {\(\)} (m-1-7.center);
  \path[-|,thick]
  (m-1-2.center) edge node [above] {\(\beta\)} (m-1-4.center) 
  (m-1-4.center) edge node [above] {\(\beta_1\)} (m-1-5.center); 
  \path[-,thick]
  (m-1-5.center) edge node [above] {\(\eta-\beta_1\)} (m-1-7.center); 
\end{tikzpicture}
\\ \\ 
  \item  \begin{tikzpicture}[baseline=2em]
  \matrix (m) [matrix of math nodes,row sep=1em,column sep=4em,minimum width=2em] {
    \mu_w & \mbox{} &  \mbox{} & \mbox{} &  \mbox{} & \mbox{} &  \mbox{} \\
    \mu_u & \mbox{} &  \mbox{} & \mbox{} &  \mbox{} & \mbox{} &  \mbox{} \\
 }
 ;
 \node (a) at ($(m-2-4)!0.5!(m-2-5)$) {};
 \node (b) at ($(m-1-4)!0.5!(m-1-5)$) [above] {\(\eta\)};
  \path[Circle-Circle,thick]
  (m-1-3.center) edge node [above] {\(\)} (m-1-5.center) 
  ;
 \path[dashed]
 (m-2-2.center) edge (a.center);
  \path[-Circle,thick]
  (a.center) edge node [above] {\(\)} (m-2-7.center);
  \path[-|,thick]
  (m-1-3.center) edge node [above] {\(\beta\)} (m-1-4.center);
  \path[|-,thick]
  (m-2-5.center) edge node [above] {\(\beta_1-\eta\)} (m-2-7.center);
\end{tikzpicture}
\end{enumerate}

We consider in the rest of the proof the case when (1) and (3) appears simultaneously and we omit the details for the remaining cases; the cases where only one of the above situations hold is the most easiest. We shall remark that the proof in all cases is similar and the definition of $\rho$ is uniform. To distinguish the roots we denote the root in picture (3) by $\beta_2$ instead of $\beta_1$ and again we assume that it is of maximal height.\medskip

\textit{Case 1:} Assume that $\varpi_i(h_{\beta})=1$. Recall that $\mu_p$ is the unique root which does not commute with $\eta$ and assume for the rest of the case that $Y:=\mu_p-\beta\in R^+$. The case $\beta-\mu_p\in R^+$ follows exactly the same idea and we omit the details. Note that we have $m_{Y}+m_{\beta}\leq \ell$ since otherwise
$$m_Y+m_{\mu_w}=m_Y+m_{\beta}+m_{\eta}=m_{\mu_p}+m_{\eta}+\ell=m_{\mu_p+\eta}+2\ell$$
which would be a contradiction. So we get
\begin{equation}\label{eq1r2}m_{\mu_p}+m_{\eta}>\ell,\ \ m_{Y}+m_{\beta}\leq \ell\end{equation}
\textit{Case 1.1:} In this case we assume that $A:=\beta_1-Y\in R^+$. So we have the following situation with $m_{\beta_1}=m_{\beta_2}=\ell$. 

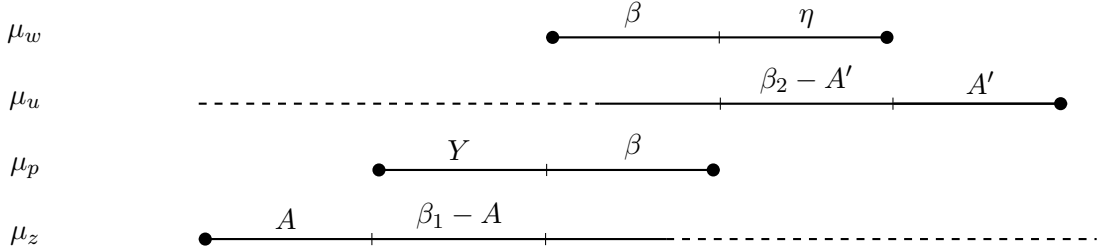
\begin{tikzpicture}[baseline=2em]
  \matrix (m) [matrix of math nodes,row sep=1em,column sep=4em,minimum width=2em] {
    \mu_w & \mbox{} &  \mbox{} & \mbox{} &  \mbox{} & \mbox{} &  \mbox{} \\
    \mu_u & \mbox{} &  \mbox{} & \mbox{} &  \mbox{} & \mbox{} &  \mbox{} \\
  \mu_p  & \mbox{} & \mbox{} & \mbox{} &  \mbox{} & \mbox{} &  \mbox{}   \\
  \mu_z& \mbox{} & \mbox{} & \mbox{} &  \mbox{} & \mbox{} &  \mbox{}   \\
 }
 ;
 \node (a) at ($(m-2-4)!0.25!(m-2-5)$) {};
 \node (b) at ($(m-4-4)!0.75!(m-4-5)$) {};
 \node (c) at ($(m-1-5)!0.5!(m-1-6)$) [above] {\(\eta\)};
 \node at ($(m-3-3)!0.5!(m-3-4)$) [above] {\(Y\)};
 \node at ($(m-2-6)!0.5!(m-2-7)$) [above] {\(A'\)};
 \path[Circle-,thick]
  (m-4-2.center) edge node [above] {\(\)} (b)
 ;
 \path[-Circle,thick]
 (a) edge node [above] {\(\)} (m-2-7.center)
 ;
 \path[Circle-Circle,thick]
  (m-1-4.center) edge node [above] {\(\)} (m-1-6.center) 
  (m-3-3.center) edge node [above] {\(\)} (m-3-5.center)
 ;
  \path[-|]
  (m-1-4.center) edge node [above] {\(\beta\)} (m-1-5.center)
  (m-4-2.center) edge node [above] {\(A\)} (m-4-3.center)
  (m-4-3.center) edge node [above] {\(\beta_1-A\)} (m-4-4.center)
  ;
  \path[|-,thin]
  (m-3-4.center) edge node [above] {\(\beta\)} (m-3-5.center)
  ;
  \path[|-|,thin]
  (m-2-5.center) edge node [above] {\(\beta_2-A'\)} (m-2-6.center)
  ;
  \path[-,thick]
  (m-2-6.center) edge node [above] {\(\)} (m-2-7.center)
  ;
  \path[dashed,thick]
  (m-2-2.center) edge node [above] {\(\)} (a.east)
  (b.west) edge node [above] {\(\)} (m-4-7.east)
  ;
\end{tikzpicture}

Moreover,
$$m_A+m_{\mu_p+\eta}=m_A+m_{\mu_p}+m_{\eta}-\ell =m_A+m_Y+m_{\beta}+m_{\eta}-\ell=m_{\beta}+m_{\eta}\leq \ell$$
where the first and second equation follows from \eqref{eq1r2}, the last equation follows from $m_{Y}+m_A=m_{\beta_1}=\ell$ (note that $m_Y=\ell$ is impossible since otherwise $m_{\mu_p}=m_{\beta}$ contradicting $m_{\mu_p}+m_{\eta}>\ell$). 
So we must have \begin{equation}\label{eqs}m_{A}+m_{\mu_z-A}\leq \ell \Longleftrightarrow \mu_{z}-A\preceq \mu_z\end{equation}
since otherwise 
$$s_{\mu_p+\eta}+s_{\mu_z}=s_{\mu_p+\eta}+s_{\mu_z-A}+s_A=s_{\mu_p+\eta+A}+s_{\mu_z-A}+1>s_{\mu_p+\eta+A}+s_{\mu_z-A}$$
which contradicts
$s_{\boldsymbol{\mu}_{\beta}}=s^{\mathrm{min}}_{\boldsymbol{\mu}_{\beta}}$.
This gives with defining relation \eqref{req1}
\begin{align} X_{\boldsymbol{\mu}_{\beta}} m& \notag=X^{p,w}_{\boldsymbol{o}(\mu)}(x^-_{\mu_p+\eta}\otimes t^{s_{\mu_p+\eta}})m
&\\&\label{ff441}=X^{p,w,z}_{\boldsymbol{o}(\mu)}(x^-_{\mu_z-A}\otimes t^{s_{\mu_z-A}})(x^-_{\mu_p+\eta+A}\otimes t^{s_{\mu_p+\eta+A}})m
\end{align}
Now we write $$(x^-_{\mu_p+\eta+A}\otimes t^{s_{\mu_p+\eta+A}})=[x^-_{\beta_1}\otimes t^{s_{\beta_1}},x^-_{\mu_w}\otimes t^{s_{\mu_w}}]$$
and obtain from \eqref{ff441} that
\begin{equation}\label{uztg} X_{\boldsymbol{\mu}_{\beta}} m=(x^-_{\beta_1}\otimes t^{s_{\beta_1}})X^{p,w,z}_{\boldsymbol{o}(\mu)}(x^-_{\mu_z-A}\otimes t^{s_{\mu_z-A}})(x^-_{\mu_w}\otimes t^{s_{\mu_w}})m\end{equation}
Now from \eqref{eqs} we already know $\mu_z-A\preceq \mu_z$ and we aim to show in the rest of this part that $\mu_z-A\preceq \mu_p$. \medskip 

If $Y_2:=\mu_p-(\mu_z-A)\in R^+$ (we set $Y_1:=\mu_z-\beta_1$) we must have $m_Y+m_{Y_1}>\ell$. To see this we consider with the converse assumption
\begin{align*}m_{Y_1}+m_{Y_2}&=m_{\mu_z}+m_{Y_2}=m_A+m_{\mu_z-A}+m_{Y_2}&\\&=m_A+m_{Y+Y_1}+m_{Y_2}= m_A+m_{Y}+m_{Y_1}+m_{Y_2}\end{align*}
where the second equation follows from \eqref{eqs}, the third from $\mu_u-A=Y+Y_1$ and the last one from the converse assumption $m_{Y}+m_{Y_1}\leq \ell $. This is a contradiction and hence we have in fact $m_Y+m_{Y_1}>\ell$. This implies
$$m_{\mu_p}+m_{Y_1}+m_{Y_2}-\ell\geq m_Y+m_{Y_1}+m_{Y_2}-\ell=m_{Y+Y_1}+m_{Y_2}\in\{m_{\mu_p},m_{\mu_p}+\ell\}$$
and thus $m_{Y+Y_1}+m_{Y_2}\leq \ell$ which is exactly $\mu_z-A\preceq \mu_p$. \medskip 

Now if $(\mu_z-A)-\mu_p\in R^+$ we assume by contradiction $\mu_z-A\succeq \mu_p$. This is impossible since 
$$m_{\mu_p+A}+m_{\mu_z-A-\mu_p}=m_{\mu_p}+m_A-\ell+m_{\mu_z-A-\mu_p}=m_{A}+m_{\mu_z-A}-\ell=m_{\mu_z}-\ell$$
where the first equation follows from $m_{\mu_p}=m_{\beta}+m_Y$, the second from $\mu_z-A\succeq \mu_p$ and the last equation from \eqref{eqs}. So we must have
\begin{equation}\label{kleinergleich}\mu_z-A\preceq \mu_p,\ \mu_z-A\preceq \mu_z.\end{equation}
Now we proceed as follows. Since $\mu_p$ has been removed in \eqref{uztg} we obtain that $\beta_2$ commutes with all other roots in the product. This gives
\begin{equation}\label{uztg1} X_{\boldsymbol{\mu}_{\beta}}m=(x^-_{\beta_1}\otimes t^{s_{\beta_1}})(x^-_{\beta_2}\otimes t^{s_{\beta_2}})X^{p,u,z}_{\boldsymbol{o}(\mu)}(x^-_{\mu_z-A}\otimes t^{s_{\mu_z-A}})(x^-_{\mu_u-\beta_2}\otimes t^{s_{\mu_u-\beta_2}})m\end{equation}
If $\mu_z-A\preceq \mu_u$ we would have together with \eqref{uztg1} and $\rho=\mu-\beta-\beta_1-\beta_2$:
$$X_{\boldsymbol{\mu}_{\beta}}m\in \mathbf{U}\cdot X_{\mathbf{o}(\rho)}m,\ \ X_{\mathbf{o}(\rho)}\prec X_{\mathbf{o}(\mu)}$$
and $\rho\in \widetilde{\mathbf{R}}_{\lambda,i}^{\ell}\backslash\{0\}$ follows from the height maximality of $\beta_1$ and $\beta_2$ respectively. So we will assume that 
\begin{equation}\label{kleinergleich22}\mu_p,\mu_z\succeq \mu_z-A\succeq \mu_u\end{equation}
The idea is to repeat the above arguments with the roots $\beta$ and $\beta_2$ instead of $\beta$ and $\beta_1$. As a first step, similarly as above, we can show
\begin{equation}\label{eqs1}m_{A'}+m_{\mu_u-A'}\leq \ell \Longleftrightarrow \mu_{u}-A'\preceq \mu_u\end{equation}
We obtain as above with \eqref{req1}
\begin{align} X_{\boldsymbol{\mu}_{\beta}}m& \notag=X^{p,w}_{\boldsymbol{o}(\mu)}(x^-_{\mu_p+\eta}\otimes t^{s_{\mu_p+\eta}})m
&\\&\notag=X^{p,w,u}_{\boldsymbol{o}(\mu)}(x^-_{\mu_u-A'}\otimes t^{s_{\mu_u-A'}})(x^-_{\mu_p+\eta+A'}\otimes t^{s_{\mu_p+\eta+A'}})m
&\\&\notag=(x^-_{\beta_2}\otimes t^{s_{\beta_2}})X^{w,u}_{\boldsymbol{o}(\mu)}(x^-_{\mu_u-A'}\otimes t^{s_{\mu_u-A'}})m
&\\&\label{ff44122}=(x^-_{\beta_1}\otimes t^{s_{\beta_1}})(x^-_{\beta_2}\otimes t^{s_{\beta_2}})X^{w,u,z}_{\boldsymbol{o}(\mu)}(x^-_{\mu_z-\beta_1}\otimes t^{s_{\mu_z-\beta_1}})(x^-_{\mu_u-A'}\otimes t^{s_{\mu_u-A'}})m
\end{align}

From \eqref{eqs1} we have $\mu_u-A'\preceq \mu_u$ and we claim in what follows that $\mu_u-A'\preceq \mu_w$. \medskip
If $Y'_2:=\gamma_w-(\mu_u-A')\in R^+$ we obtain as above with the assumption $m_{\eta}+m_{Y'_1}\leq \ell$ a contradiction (we set $Y_1':=\mu_u-A'-\tau$):
$$m_{Y_1'}+m_{Y_2'}=m_{Y_2'}+m_{\mu_u}=m_{Y_1'}+m_{\mu_u-A'}+m_{A'}=m_{Y_2'}+m_{Y_1'+\eta}+m_{A'}=m_{Y_2'}+m_{Y_1'}+\ell.$$
Hence $m_{\eta}+m_{Y'_1}>\ell$ which gives 
$$m_{\mu_w}+m_{Y_1'}+m_{Y_2'}-\ell\geq m_{\eta}+m_{Y_1'}+m_{Y_2'}-\ell =m_{\eta+Y'_1}+m_{Y_2'}\in\{m_{\mu_w},m_{\mu_w}+\ell\}$$
Therefore $m_{\eta+Y'_1}+m_{Y_2'}\leq \ell$ and $\mu_u-A'\preceq \mu_w$ is obtained. \medskip 

If $(\mu_u-A')-\gamma_w\in R^+$ we assume by contradiction $\gamma_u-A'\succeq \gamma_w$. This gives
$$m_{\mu_w+A'}+m_{\mu_u-A'-\mu_w}=m_{\mu_w}+m_{A'}-\ell+m_{\mu_u-A'-\mu_w}=m_{A'}+m_{\mu_u-A'}-\ell=m_{\mu_u}-\ell$$
which is impossible. So summarizing we get with \eqref{kleinergleich22}
$$\mu_p,\mu_z\succeq \mu_z-A\succeq \mu_u\succeq \mu_u-A',\ \mu_u-A'\preceq \mu_w$$
which finally gives with \eqref{ff44122}
$$X_{\boldsymbol{\mu}_{\beta}}m\in \mathbf{U}\cdot X_{\mathbf{o}(\rho)}m,\ \ X_{\mathbf{o}(\rho)}\prec X_{\mathbf{o}(\mu)}$$
where $\rho=\mu-\beta-\beta_1-\beta_2$ and $\rho\in \widetilde{\mathbf{R}}_{\lambda,i}^{\ell}\backslash\{0\}$ follows again from the height maximality of $\beta_1$ and $\beta_2$ respectively. 

\textit{Remark: Note that we haven't used the assumption of Case 1.1 in order to obtain \eqref{ff44122} and the inequalities $\mu_u-A'\preceq \mu_w$ and $\mu_u-A'\preceq \mu_u$. We emphasize this here, since we want to use \eqref{ff44122} also in Case 1.2.}

\textit{Case 1.2:} In this case we assume that $A'':=Y-\beta_1\in R^+$. So we have the following situation with $m_{\beta_1}=m_{\beta_2}=\ell$. 

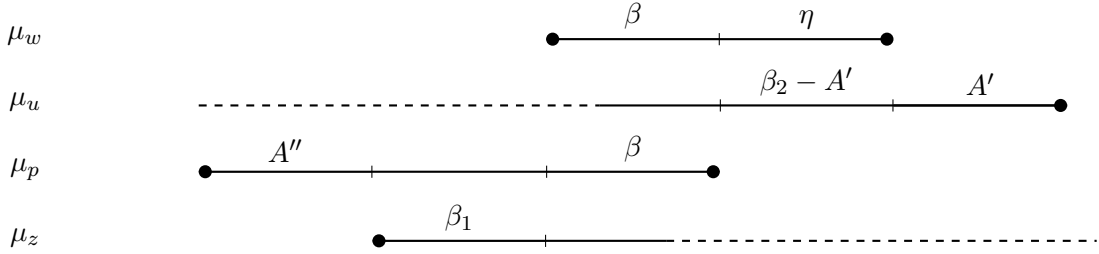
\begin{tikzpicture}[baseline=2em]
  \matrix (m) [matrix of math nodes,row sep=1em,column sep=4em,minimum width=2em] {
    \mu_w & \mbox{} &  \mbox{} & \mbox{} &  \mbox{} & \mbox{} &  \mbox{} \\
    \mu_u & \mbox{} &  \mbox{} & \mbox{} &  \mbox{} & \mbox{} &  \mbox{} \\
  \mu_p  & \mbox{} & \mbox{} & \mbox{} &  \mbox{} & \mbox{} &  \mbox{}   \\
  \mu_z& \mbox{} & \mbox{} & \mbox{} &  \mbox{} & \mbox{} &  \mbox{}   \\
 }
 ;
 \node (a) at ($(m-2-4)!0.25!(m-2-5)$) {};
 \node (b) at ($(m-4-4)!0.75!(m-4-5)$) {};
 \node (c) at ($(m-1-5)!0.5!(m-1-6)$) [above] {\(\eta\)};
  \node at ($(m-2-6)!0.5!(m-2-7)$) [above] {\(A'\)};
\path[|-|,thin]
  (m-2-5.center) edge node [above] {\(\beta_2-A'\)} (m-2-6.center)
  ;
 \path[Circle-,thick]
  (m-4-3.center) edge node [above] {\(\)} (b)
 ;
 \path[-Circle,thick]
 (a) edge node [above] {\(\)} (m-2-7.center)
 ;
 \path[Circle-Circle,thick]
  (m-1-4.center) edge node [above] {\(\)} (m-1-6.center) 
  (m-3-2.center) edge node [above] {\(\)} (m-3-5.center)
 ;
  \path[-|]
  (m-1-4.center) edge node [above] {\(\beta\)} (m-1-5.center)
  (m-3-2.center) edge node [above] {\(A''\)} (m-3-3.center)
  (m-4-3.center) edge node [above] {\(\beta_1\)} (m-4-4.center);
  \path[|-,thin]
  (m-3-4.center) edge node [above] {\(\beta\)} (m-3-5.center)
  ;
  \path[-,thick]
  (m-2-6.center) edge node [above] {\(\)} (m-2-7.center);
  
  \path[dashed,thick]
  (m-2-2.center) edge node [above] {\(\)} (a.east)
  (b.west) edge node [above] {\(\)} (m-4-7.east)
  ;
\end{tikzpicture}

We first note that $m_{A''}+m_{\beta+\beta_1}\leq \ell $ since otherwise
$$m_{A''}+m_{\beta+\beta_1}+m_{\eta}=m_{\mu_p}+\ell+m_{\eta}=m_{\mu_p+\eta}+2\ell$$
which is absurd. This gives in particular
$$m_{A''}+m_{\gamma_w+\beta_1}=m_{A''}+m_{\beta+\beta_1}+m_{\eta}=m_{\mu_p}+m_{\eta}>\ell$$
Thus we must also have with $s_{\boldsymbol{\mu}_\beta}=s^{\mathrm{min}}_{\boldsymbol{\mu}_\beta}$
\begin{equation}\label{hhggttq}m_{\mu_z}+m_{A''}>\ell \iff \mu_z+A''\preceq \mu_z\end{equation}
Now with \eqref{req1} we get 
\begin{align} X_{\boldsymbol{\mu}_{\beta}}m& \notag=X^{p,w}_{\boldsymbol{o}(\mu)}(x^-_{\mu_p+\eta}\otimes t^{s_{\mu_p+\eta}})m
&\\&\notag=X^{p,w,z}_{\boldsymbol{o}(\mu)}(x^-_{\mu_z+A''}\otimes t^{s_{\mu_z+A''}})(x^-_{\mu_w+\beta_1}\otimes t^{s_{\mu_w+\beta_1}})m
&\\&\label{44dfq}=(x^-_{\beta_1}\otimes t^{s_{\beta_1}})X^{p,z}_{\boldsymbol{o}(\mu)}(x^-_{\mu_z+A''}\otimes t^{s_{\mu_z+A''}})m
\end{align}
From \eqref{hhggttq} we know $\mu_z+A''\preceq \mu_z$ and we aim to show that we also have $\mu_z+A''\preceq \mu_p$. \medskip 

If $Y_2:=\mu_p-(\mu_z+A'')\in R^+$ (we set $Y_1:=\mu_z-\beta_1$) we must have $m_{\mu_p-Y_2}+m_{Y_2}\leq \ell$, since otherwise
\begin{align*}m_{Y_1}+m_{Y_2}&=m_{\mu_z}+m_{Y_2}=m_{\mu_z+A''}+\ell-m_{A''}+m_{Y_2}&\\&=m_{\mu_p-Y_2}+\ell-m_{A''}+m_{Y_2}=m_{\mu_p}-m_{A''}+2\ell=m_{\beta}+2\ell\end{align*}
which is a contradiction, where the second equation follows from \eqref{hhggttq}, the third from $\mu_z+A''=\mu_p-Y_2$ and the last one from $m_{\mu_p}=m_{A''}+m_{\beta}$. Thus $\mu_z+A''\preceq \mu_p$. \medskip 

If $(\mu_z+A'')-\mu_p\in R^+$ we assume by contradiction $\mu_z+A''\succeq \mu_p$. But this is impossible since 
$$m_{\mu_p-A''}+m_{\mu_z+A''-\mu_p}=m_{\mu_p}-m_{A''}+m_{\mu_z+A''-\mu_p}=m_{\mu_z+A''}-m_{A''}=m_{\mu_z}-\ell$$
where the first equation follows from $m_{\mu_p}=m_{\mu_p-A''}+m_{A''}$, the second from $\gamma_z+A''\succeq \gamma_p$ and the last equation from \eqref{hhggttq}. So we must have
\begin{equation}\label{kleinergleich1251}\mu_z+A''\preceq \mu_p,\ \mu_z+A''\preceq \mu_z.\end{equation}
Now we proceed as follows; the idea is exactly the same as before. In a first step we remove $\beta_2$ from $\mu_u$ and get with \eqref{44dfq}
\begin{align}\notag X_{\boldsymbol{\mu}_{\beta}}m =(x^-_{\beta_1}\otimes t^{s_{\beta_2}})(x^-_{\beta_1}\otimes t^{s_{\beta_2}})X^{p,z,u}_{\boldsymbol{o}(\mu)}(x^-_{\mu_z+A''}\otimes t^{s_{\mu_z+A''}})(x^-_{\mu_u-\beta_2}\otimes t^{s_{\mu_u-\beta_2}})m
\end{align}
If $\mu_z+A''\preceq \mu_u$ we would have as before
$$X_{\boldsymbol{\mu}_{\beta}}m\in \mathbf{U}\cdot X_{\mathbf{o}(\rho)}m,\ \ X_{\mathbf{o}(\rho)}\prec X_{\mathbf{o}(\mu)}$$
where $\rho=\mu-\beta-\beta_1-\beta_2$ and $\rho\in \widetilde{\mathbf{R}}_{\lambda,i}^{\ell}\backslash\{0\}$ follows again from the height maximality of $\beta_1$ and $\beta_2$ respectively. 
So we will assume that 
\begin{equation}\label{kleinergleich2277ü}\mu_p,\mu_z\succeq \mu_z+A''\succeq \mu_u\end{equation}
The idea is again the same and we repeat the above arguments with the roots $\beta$ and $\beta_2$ instead of $\beta$ and $\beta_1$. Recall from \eqref{ff44122} (which was independent of the assumption in Case 1.1)
\begin{align} X_{\boldsymbol{\mu}_{\beta}}m= (x^-_{\beta_1}\otimes t^{s_{\beta_1}})(x^-_{\beta_2}\otimes t^{s_{\beta_2}})X^{w,u,z}_{\boldsymbol{o}(\mu)}(x^-_{\mu_z-\beta_1}\otimes t^{s_{\mu_z-\beta_1}})(x^-_{\mu_u-A'}\otimes t^{s_{\mu_u-A'}})m
\end{align}
with $\mu_u-A'\preceq \mu_u$ and $\mu_u-A'\preceq \mu_w$. However with \eqref{kleinergleich2277ü}
$$\mu_p,\mu_z\succeq \mu_z+A''\succeq \mu_u\succeq \mu_u-A',\ \ \mu_u-A'\preceq \mu_w$$
and we obtain once more the desired property.

\textit{Case 2:} We assume that $\varpi_i(h_{\beta})=0$.

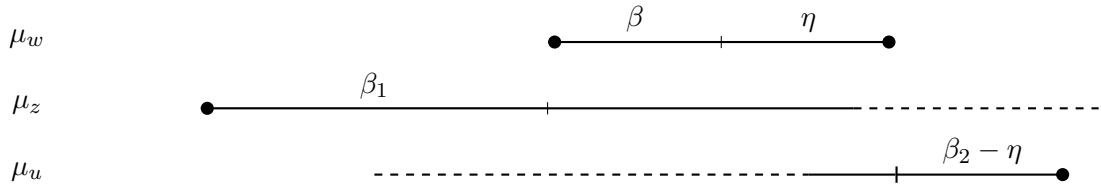
\begin{tikzpicture}[baseline=2em]
  \matrix (m) [matrix of math nodes,row sep=1em,column sep=4em,minimum width=2em] {
    \mu_w & \mbox{} &  \mbox{} & \mbox{} &  \mbox{} & \mbox{} &  \mbox{} \\
    \mu_z& \mbox{} & \mbox{} & \mbox{} &  \mbox{} & \mbox{} &  \mbox{}   \\
    \mu_u& \mbox{} & \mbox{} & \mbox{} &  \mbox{} & \mbox{} &  \mbox{}   \\
 }
 ;
 \node (b) at ($(m-2-5)!0.7!(m-2-6)$) {};
 \node (c) at ($(m-1-5)!0.5!(m-1-6)$) [above] {\(\eta\)};
 \node (r) at ($(m-3-6)!0.5!(m-3-7)$) [] {};
 \node (d) at ($(m-3-5)!0.5!(m-3-6)$) {}; 
 \path[Circle-,thick]
  (m-2-2.center) edge node [above] {\(\)} (b.east)
 ;
 \path[|-Circle,thick]
 (m-3-6.center) edge node [above] {\(\beta_2-\eta\)} (m-3-7.center)
 ;
 \path[-,thick]
 (d.center) edge node [above] {\(\)} (m-3-6.center)
 ;
 \path[Circle-Circle,thick]
  (m-1-4.center) edge node [above] {\(\)} (m-1-6.center) 
  ;
  \path[-|]
  (m-1-4.center) edge node [above] {\(\beta\)} (m-1-5.center)
  (m-2-2.center) edge node [above] {\(\beta_1\)} (m-2-4.center)
  ;
  \path[dashed,thick]
  (b.east) edge node [above] {\(\)} (m-2-7.east)
  (m-3-3.center) edge node [above] {\(\)} (d.center)
  ;
\end{tikzpicture}

Since $m_{\beta+\beta_1}+m_{\eta}\leq \ell$ we must have by the degree minimality $s_{\boldsymbol{\mu}_\beta}=s^{\mathrm{min}}_{\boldsymbol{\mu}_\beta}$ also 
\begin{equation}\label{55t}m_{\mu_z-\beta-\beta_1}+m_{\beta+\beta_1}\leq \ell.\end{equation} 
\textit{Case 2.1:} We consider first the case \(A:=\mu_z-\beta-\beta_1-\beta_2 \in R^+\).  
So \eqref{req1} gives
\begin{align*}
	X_{\boldsymbol{\mu}_\beta}m&= (x^-_{\beta_1} \otimes t^{s_{\beta_1}}) X^{z}_{\boldsymbol{o}(\mu)}(x^-_{\beta_2+A} \otimes t^{s_{\beta_2+A}})m\\ &= (x^-_{\beta_1} \otimes t^{s_{\beta_1}}) (x^-_{\beta_2} \otimes t^{s_{\beta_2}}) X^{z,u}_{\boldsymbol{o}(\mu)}(x^-_{\mu_u+A} \otimes t^{s_{\mu_u}+s_A})m,
\end{align*}
where the second equation follows from $m_{\beta_2}=\ell$ and $(x^-_{\beta_2} \otimes t^{s_{\beta_2}})(x^-_{\mu_u} \otimes t^{s_{\mu_u}})m=0$.
Now we continue as follows. The above element vanishes if $s_{\mu_u}+s_A>s_{\mu_u+A}$, so that we can assume 
$$m_{\mu_u}+ m_A > \ell \iff \mu_u+A \preceq \mu_u.$$
If \(\mu_u+A \preceq \mu_z\) we are done as usual with the maximal height property of $\beta_1$ and $\beta_2$. 
So suppose additionally that \(\mu_u+A \succeq \mu_z\) and consider first the subcase 
$$A':=\mu_u -\beta-\beta_1-\beta_2 \in R^+.$$
In this subcase \(\mu_u+A \preceq \mu_z\) means
$$m_{A'}+m_{\mu_z}\leq \ell \Rightarrow m_{A'}+m_{\gamma_u-A'} \leq \ell$$
where the implication follows again from the degree minimality $s_{\boldsymbol{\mu}_\beta}=s^{\mathrm{min}}_{\boldsymbol{\mu}_\beta}$. So we have once more
\begin{align*}
	X_{\boldsymbol{\mu}_\beta}m&= X^{z,u,w}_{\boldsymbol{o}(\mu)}(x^-_{\eta} \otimes t^{s_{\eta}})(x^-_{\mu_z+A'} \otimes t^{s_{\mu_z+A'}})(x^-_{\beta+\beta_1+\beta_2} \otimes t^{s_{\beta+\beta_1+\beta_2}})m&\\&
 = X^{z,u,w}_{\boldsymbol{o}(\mu)}(x^-_{\gamma_w+\beta_1} \otimes t^{s_{\eta}+s_{\beta+\beta_1}})(x^-_{\mu_z+A'} \otimes t^{s_{\mu_z+A'}})(x^-_{\beta_2} \otimes t^{s_{\beta_2}})m=0
\end{align*}
since 
\begin{equation}\label{nhg}s_{\eta}+s_{\beta+\beta_1}=s_{\eta}+s_{\beta}+s_{\beta_1}=s_{\gamma_w+\beta_1}+1.\end{equation}
The idea for the subcase \(-A'\in R^+\) is exactly the same. Namely, we move $-A'$ from $\mu_z$ to $\mu_u$ and remain with a root $\beta+\beta_1+\beta_2$. We omit the details. 

\textit{Case 2.2:} Now we consider the case that \(B=\beta+\beta_1+\beta_2-\mu_z\in R^+\). 
We proceed by first noting that \(m_{B}+m_{\mu_u-B}\leq \ell\), otherwise if \(m_{B}+m_{\mu_u-B}> \ell\) we would also have \(m_{B}+m_{\mu_z}>\ell\) by the degree minimality and we would obtain once more with \eqref{req1}
\begin{align*}
	X_{\boldsymbol{\mu}_\beta}m&= X^{u,w,z}_{\boldsymbol{o}(\mu)}(x^-_{\eta} \otimes t^{s_{\eta}})(x^-_{\mu_u-B} \otimes t^{s_{\mu_u-B}})(x^-_{\beta+\beta_1+\beta_2} \otimes t^{s_{\beta+\beta_1+\beta_2}})m&\\&
 =X^{u,w,z}_{\boldsymbol{o}(\mu)} (x^-_{\gamma_w+\beta_1} \otimes t^{s_{\eta}+s_{\beta+\beta_1}})(x^-_{\mu_u-B} \otimes t^{s_{\mu_u-B}})(x^-_{\beta_2} \otimes t^{s_{\beta_2}})m
\end{align*}
which vanishes for the same reason as above (see \eqref{nhg}). So we assume in the rest of the proof that
$m_{B}+m_{\mu_u-B}\leq \ell$ which means $\mu_u-B \preceq \mu_u$. We obtain 
\begin{align*}
	X_{\boldsymbol{\mu}_\beta}m&=  (x^-_{\beta_1} \otimes t^{s_{\beta_1}})X^{z}_{\boldsymbol{o}(\mu)}(x^-_{\mu_z-\beta_1-\beta} \otimes t^{s_{\mu_z-\beta-\beta_1}})m \\
	&=(x^-_{\beta_1} \otimes t^{s_{\beta_1}})X^{u,z}_{\boldsymbol{o}(\mu)}(x^-_{\mu_u-B} \otimes t^{s_{\mu_u-B}}) (x^-_{\mu_z-\beta_1-\beta+B} \otimes t^{s_{\mu_z-\beta-\beta_1+B}})m \\
	&=(x^-_{\beta_1} \otimes t^{s_{\beta_1}})(x^-_{\beta_2} \otimes t^{s_{\beta_2}})X^{u,z}_{\boldsymbol{o}(\mu)}(x^-_{\mu_u-B} \otimes t^{s_{\mu_u-B}})m 
\end{align*}
where the first equality is obtained by splitting off $\beta_1$ to the front (this root commutes with all other roots) and passing $\beta$ from $\mu_z-\beta_1$ to $\eta$ (this is possible since \eqref{55t} holds). The second equality follows from \eqref{req1} by noting that
$$m_{B}+m_{\mu_u-B}\leq \ell,\ \ m_{B}+m_{\mu_z-\beta-\beta_1}=m_{\beta_2}\leq \ell.$$
The third equality is implied by \(\mu_z-\beta_1-\beta+B = \beta_2\).\medskip

It is left to show that \(\mu_u-B \preceq \mu_z\). If \(A'=\mu_u-\beta-\beta_1-\beta_2\in R^+\) this is equivalent to \(m_{\mu_z}+m_{A'} > \ell \). If this is false, i.e. \(m_{\mu_z}+m_{A'} \leq \ell \) we can move $A'$ to $\mu_z$ with \eqref{req1} (note that we also have $m_{\mu_u-A'}+m_{A'}\leq \ell$ by the degree minimality). After this step we have
\begin{align*}X_{\boldsymbol{\mu}_\beta}m&=  (x^-_{\beta_1} \otimes t^{s_{\beta_1}})X^{u,w,z}_{\boldsymbol{o}(\mu)}(x^-_{\eta} \otimes t^{s_{\eta}}) (x^-_{\mu_z+A'} \otimes t^{s_{\mu_z+A'}})(x^-_{\beta+\beta_1+\beta_2} \otimes t^{s_{\beta+\beta_1+\beta_2}})m&\\&
=(x^-_{\beta_1} \otimes t^{s_{\beta_1}})X^{u,w,z}_{\boldsymbol{o}(\mu)}(x^-_{\beta_1+\mu_w} \otimes t^{s_{\eta}+s_{\beta+\beta_1}}) (x^-_{\mu_z+A'} \otimes t^{s_{\mu_z+A'}})(x^-_{\beta_2} \otimes t^{s_{\beta_2}})m
\end{align*}
Together with \eqref{nhg} we would get again $X_{\boldsymbol{\mu}_\beta}m=0$. \medskip
If \(-A'\in R^+\) the statement is equivalent to \(m_{\mu_z+A'}+ m_{-A'} \leq \ell\). By exactly the same reason as above the converse assumption (moving $-A'$ to $\mu_u$ etc.) would lead to a contradiction. This finishes the proof. 
\end{proof}

\begin{rem}\label{remnplus} Given $\mu\in \widetilde{\mathbf{R}}_{\lambda,i}^{\ell}\backslash\{0\}$ and $\beta\in R^+$, we either have $(x_{\beta}^+\otimes 1)X_{\mathbf{o}(\mu)}m=0$ or the assumptions of Lemma~\ref{nplus} hold. To see this, let $(x_{\beta}^+\otimes 1)X_{\mathbf{o}(\mu)}m\neq 0$, i.e. $\varpi_i-(\mu-\beta)\in W(\varpi_i)$ and there exists $w\in\{1,\dots,k\}$ with $\mu_w-\beta\in R^+$. Since $(x_{\mu_w-\beta}^-\otimes t^{s_{\mu_w-\beta}+\varpi_i(h_{\mu_w-\beta})})m=0$ we also have
$$s_{\mu_w}=s_{\mu_w-\beta}=s_{\mu_w-\beta}+s_{\beta}-1\Rightarrow s_{\beta}=1,\ m_{\beta}+m_{\mu_w-\beta}\leq \ell.$$
Moreover, if $\varpi_i(h_{\beta})=1$, there exists another index $p\in\{1,\dots,k\}$ such that $\mu_w-\beta+\mu_p\in R^+$. So $(x_{\beta}^+\otimes 1)X_{\mathbf{o}(\mu)}m$ is proportional to
$$Y(x_{\mu_w-\beta+\mu_p}^-\otimes t^{s_{\mu_w-\beta}+s_{\mu_p}})m$$
which forces $m_{\mu_p}+m_{\mu_w-\beta}>\ell.$  The property $s_{\boldsymbol{\mu}_{\beta}}=s_{\boldsymbol{\mu}_{\beta}}^{\mathrm{min}}$ is a direct consequence of \eqref{req1}. So all assumptions are checked. 
\end{rem}

\subsection{}The following proposition will give an upper bound for the characters involved in the Pieri formula. 
\begin{prop}\label{mainprophel} Let $\mathbf{o}(\mu)=(\mu_1,\dots,\mu_k)$ be the orbit decomposition of $\mu\in \widetilde{\mathbf{R}}_{\lambda,i}^{\ell}\backslash\{0\}$. Let $j\in\{1,\dots,s\}$ be the index such that $\mathbf{o}(\mu)=\mathbf{o}(\mu)_j$ (recall the notation from \eqref{orddef}). Then we have a surjective map of graded $\mathbf{U}$-modules
$$\mathbf{D}^{\ell}_{\lambda+\varpi_i-\mu}\rightarrow \tau_{s_{\mathbf{o}(\mu)}} K_{j}/K_{j-1}\rightarrow 0.$$
with cyclic generator $X_{\mathbf{o}(\mu)}m.$

\begin{proof}We will step wise show that the defining relations of the Demazure module are satisfied by the cyclic generator. The relations 
$$(\mathfrak{n}^+\otimes 1)X_{\mathbf{o}(\mu)}m=0=(h\otimes t^{k+1})X_{\mathbf{o}(\mu)}m=0,\ \ k\geq 0$$
hold where the second part is clear and the first part is the statement of Lemma~\ref{nplus} together with Remark~\ref{remnplus}. Now we consider the remaining relations.

\textit{Case 1}: Let $\gamma\in R^+$ such that $(\varpi_i-\mu)(h_{\gamma})=1$.  
We have to show that the following identities hold
$$(x^-_{\gamma}\otimes t^{s_{\gamma}-1})^{m_{\gamma}+2}X_{\mathbf{o}(\mu)}m=0,\ \ (x^-_{\gamma}\otimes t^{s_{\gamma}})X_{\mathbf{o}(\mu)}m=0,\ \ m_{\gamma}<\ell$$
$$(x^-_{\gamma}\otimes t^{s_{\gamma}})^{2}X_{\mathbf{o}(\mu)}m=0,\ \ (x^-_{\gamma}\otimes t^{s_{\gamma}+1})X_{\mathbf{o}(\mu)}m=0,\ \ m_{\gamma}=\ell$$
where the only non trivial relation is 
\begin{equation}\label{nontg4}(x^-_{\gamma}\otimes t^{s_{\gamma}})X_{\mathbf{o}(\mu)}m=0,\ \ m_{\gamma}<\ell\end{equation}
by the following argument. If $\varpi_i(h_{\gamma})=1$, then it commutes with $X_{\mathbf{o}(\mu)}$ and the three other relations are coming from the defining relations of $\mathbf{M}^{\ell}_{\lambda,i}$. Otherwise $\varpi_i(h_{\gamma})=0$ and $\gamma$ does not commute with $X_{\mathbf{o}(\mu)}$ (otherwise the statement is again clear). But then we have a unique root $\mu_u$ which does not commute with $\gamma$ since $\mu_u\in \widetilde{\mathbf{R}}_{\lambda,i}^{\ell}\backslash\{0\}$. Now
$$s_{\mu_u}+s_{\gamma}+1>s_{\mu_u+\gamma}$$
implies $(x^-_{\gamma}\otimes t^{s_{\gamma}+1})X_{\mathbf{o}(\mu)}m=0$.
Moreover, since $(x^-_{\gamma}\otimes t^{s_{\gamma}})m=0$
\begin{align*}(x^-_{\gamma}\otimes t^{s_{\gamma}})^{2}X_{\mathbf{o}(\mu)}m =(x^-_{\gamma}\otimes t^{s_{\gamma}})(x^-_{\mu_u+\gamma}\otimes t^{s_{\mu_u}+s_{\gamma}})X^u_{\mathbf{o}(\mu)}m=0.
\end{align*}
Similarly 
\begin{align*}&(x^-_{\gamma}\otimes t^{s_{\gamma}-1})^{m_{\gamma}+2}X_{\mathbf{o}(\mu)}m&\\&=(x^-_{\gamma}\otimes t^{s_{\gamma}-1})^{m_{\gamma}+1}X_{\mathbf{o}(\mu)}(x^-_{\gamma}\otimes t^{s_{\gamma}-1})m+(x^-_{\gamma}\otimes t^{s_{\gamma}-1})^{m_{\gamma}+1}(x^-_{\mu_u+\gamma}\otimes t^{s_{\mu_u}+s_{\gamma}-1})X^u_{\mathbf{o}(\mu)}m&\\&
=(x^-_{\gamma}\otimes t^{s_{\gamma}-1})^{m_{\gamma}+1}X_{\mathbf{o}(\mu)}(x^-_{\gamma}\otimes t^{s_{\gamma}-1})m&\\&
=(x^-_{\gamma}\otimes t^{s_{\gamma}-1})^{m_{\gamma}}X_{\mathbf{o}(\mu)}(x^-_{\gamma}\otimes t^{s_{\gamma}-1})^2m+(x^-_{\mu_u+\gamma}\otimes t^{s_{\mu_u}+s_{\gamma}-1})X^u_{\mathbf{o}(\mu)}(x^-_{\gamma}\otimes t^{s_{\gamma}-1})^{m_{\gamma}+1}m&\\&
=(x^-_{\gamma}\otimes t^{s_{\gamma}-1})^{m_{\gamma}}X_{\mathbf{o}(\mu)}(x^-_{\gamma}\otimes t^{s_{\gamma}-1})^2m=\cdots=X_{\mathbf{o}(\mu)}(x^-_{\gamma}\otimes t^{s_{\gamma}-1})^{m_{\gamma}+2}m=0
\end{align*}
\medskip 
So we have to show the non-trivial relation \eqref{nontg4} only.

\textit{Case 1.1}: Suppose that $\varpi_i(h_{\gamma})=0$. If $\gamma$ commutes with $X_{\mathbf{o}(\mu)}$ we are done. So let $\mu_w$ the unique root such that $\gamma+\mu_w\in R^+$. Note that 
\begin{equation}\label{bbvx50}(x^-_{\gamma}\otimes t^{s_{\gamma}})X_{\mathbf{o}(\mu)}m=X^w_{\mathbf{o}(\mu)}(x^-_{\mu_w+\gamma}\otimes t^{s_{\mu_w}+s_{\gamma}})m\end{equation}
and we can assume that $s_{\mu_w}+s_{\gamma}=s_{\mu_w+\gamma}$ (otherwise the above element vanishes). This means $m_{\mu_w}+m_{\gamma}>\ell$ and thus $\mu_w+\gamma\preceq \gamma_w$. If the right hand side of \eqref{bbvx50} corresponds to an element in $\widetilde{\mathbf{R}}_{\lambda,i}^{\ell}\backslash\{0\}$ we are done. Otherwise there exists a positive root $\beta_1$ such that $m_{\beta_1}=\ell$ and $\varpi_i-\mu+\beta_1\notin W(\varpi_i)$, $\varpi_i-\mu-\gamma+\beta_1\in W(\varpi_i)$ (as in the proof of Lemma~\ref{nplus}). Again we have the following local situation

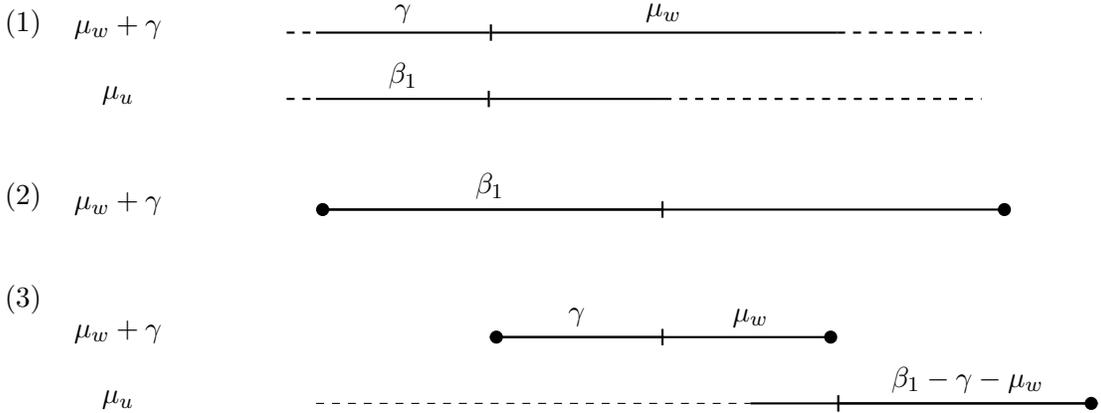
\begin{enumerate}
\item  
\begin{tikzpicture}[baseline=1em]
  \matrix (m) [matrix of math nodes,row sep=1em,column sep=4em,minimum width=2em] {
    \mu_w+\gamma& \mbox{} &  \mbox{} & \mbox{} &  \mbox{} & \mbox{} &  \mbox{} \\
  \mu_u & \mbox{} & \mbox{} & \mbox{} &  \mbox{} & \mbox{} &  \mbox{}   \\}
  ;
  
  \node (a) at ($(m-2-3)!0.5!(m-2-2)$) {}
  ;
  \node (b) at ($(m-1-3)!0.5!(m-1-2)$) {}
  ;
  \node at (b) [above] {\(\gamma\)} 
  ;
  \node at (a) [above] {\(\beta_1\)} 
  ;
  
  \path[|-,thick]
  (m-1-3.center) edge node [above] {\(\mu_w\)} (m-1-5.center) 
  ;
  \path[-|,thick]
  (m-2-2.center) edge (m-2-3.center)
  ;
  \path[-,thick]
  (m-2-3.center) edge (m-2-4.center) 
  ;
  \path[-,thick]
  (m-1-2.center) edge (m-1-3.center) 
  ;
  \path[-,thick, dashed]
  (m-2-2.west) edge node [above] {} (m-2-2.center)
  (m-2-4.center) edge node [above] {} (m-2-6.west) 
  (m-1-5.center) edge node [above] {} (m-1-6.west)
  (m-1-2.west) edge node [above] {} (m-1-2.center)
  ;
  \end{tikzpicture}
\\
\\
  \item  
  \begin{tikzpicture}[baseline=0em]
  \matrix (m) [matrix of math nodes,row sep=1em,column sep=4em,minimum width=2em] {
  \mu_w+\gamma & \mbox{} & \mbox{} & \mbox{} &  \mbox{} & \mbox{} &  \mbox{}   \\}
  ;
  \path[Circle-Circle,thick]
  (m-1-2.center) edge node [above] {\(\)} (m-1-6.center)
  ;
  \path[-|,thick]
  (m-1-2.center) edge node [above] {\(\beta_1\)} (m-1-4.center)
  ;
  \end{tikzpicture}
\\ \\ 
  \item  
  \begin{tikzpicture}[baseline=2em]
  \matrix (m) [matrix of math nodes,row sep=1em,column sep=4em,minimum width=2em] {
    \mu_w +\gamma& \mbox{} &  \mbox{} & \mbox{} &  \mbox{} & \mbox{} &  \mbox{} \\
    \mu_u & \mbox{} &  \mbox{} & \mbox{} &  \mbox{} & \mbox{} &  \mbox{} \\
  }
  ;
  \node (a) at ($(m-2-4)!0.5!(m-2-5)$) {}
  ;
  \node (b) at ($(m-1-4)!0.5!(m-1-5)$) [above] {\(\mu_w\)}
  ;
  \node (c) at ($(m-2-6)!0.5!(m-2-7)$) {}
  ;
  \path[Circle-Circle,thick]
  (m-1-3.center) edge node [above] {\(\)} (m-1-5.center) 
  ;
  \path[dashed]
  (m-2-2.center) edge (a.center)
  ;
  \path[-Circle,thick]
  (a.center) edge node [above] {\(\)} (c.center)
  ;
  \path[-|,thick]
  (m-1-3.center) edge node [above] {\(\gamma\)} (m-1-4.center)
  ;
  \path[|-,thick]
  (m-2-5.center) edge node [above] {\(\beta_1-\gamma-\mu_w\)} (c.center)
  ;
\end{tikzpicture}
\end{enumerate}
As in the proof of Lemma~\ref{nplus} we restrict ourselves when case (1) and (3) appears simultaneously which means the following constellation of roots

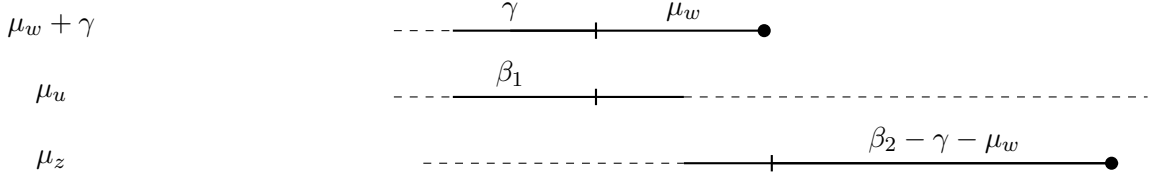
\begin{tikzpicture}[baseline=2em]
  \matrix (m) [matrix of math nodes,row sep=1em,column sep=4em,minimum width=2em] {
    \mu_w +\gamma & \mbox{} &  \mbox{} & \mbox{} &  \mbox{} & \mbox{} &  \mbox{} \\
    \mu_u & \mbox{} &  \mbox{} & \mbox{} &  \mbox{} & \mbox{} &  \mbox{} \\
    \mu_z & \mbox{} &  \mbox{} & \mbox{} &  \mbox{} & \mbox{} &  \mbox{} \\
  }
  ;
  \node (a) at ($(m-3-4)!0.5!(m-3-5)$) {}
  ;
  \node (b) at ($(m-1-4)!0.5!(m-1-5)$) [above] {\(\mu_w\)}
  ;
  \node (c) at ($(m-2-3)!0.5!(m-2-4)$) {}
  ;
  \node (d) at ($(m-2-4)!0.5!(m-2-5)$) {}
  ;
  \node at (c) [above] {\(\beta_1\)} 
  ;
  \node (e) at ($(m-1-3)!0.5!(m-1-4)$) {}
  ;
  \node at (e) [above] {\(\gamma\)}
  ;
  \node (f) at ($(m-3-3)!0.5!(m-3-4)$) {}
  ;
  \path[-Circle,thick]
  (e.center) edge node [above] {\(\)} (m-1-5.center) 
  ;
  \path[dashed]
  (m-3-3.center) edge (a.center)
  (m-1-3.west) edge (m-1-3.east)
  (m-2-3.west) edge node [above] {} (m-2-3.east)
  (d.center) edge node [above] {} (m-2-7.east)
  ;
  \path[-Circle,thick]
  (a.center) edge node [above] {\(\)} (m-3-7.center)
  ;
  \path[-|,thick]
  (m-1-3.east) edge node [above] {} (m-1-4.center)
  ;
  \path[|-,thick]
  (m-3-5.center) edge node [above] {\(\beta_2-\gamma-\mu_w\)} (m-3-7.center)
  ;
  \path[-|,thick]
  (m-2-3.east) edge (m-2-4.center)
  ;
  \path[-,thick]
  (m-2-4.center) edge (d.center) 
  ;
\end{tikzpicture}
\\
We can assume that $\beta_1$ and $\beta_2$ are of maximal height.  Define \(B=\beta_2-\gamma-\mu_w\) and note that we must have \(m_{\mu_w}+m_{B} > \ell\) since otherwise (note that $m_{\gamma+\mu_w}=m_{\gamma}+m_{\mu_w}-\ell\neq \ell$)
$$2\ell=m_{\gamma+\mu_w}+m_B+\ell=m_{\gamma}+m_{\mu_w}+m_{B}=m_{\gamma}+m_{\mu_w+B}.$$

\textit{Case 1.1.1}: Suppose first that \(m_{B}+ m_{\mu_z-B} > \ell\). 
This implies with \eqref{req1}
\begin{align*}
	X_{\boldsymbol{o}(\mu)}m=
	X^{w,z}_{\boldsymbol{o}(\mu)}(x^-_{\mu_w+B}\otimes t^{s_{\gamma_w +B}}) (x^-_{\mu_z-B}\otimes t^{s_{\mu_z-B}})m.
\end{align*}
Thus
\begin{align*}
	(x^-_{\gamma}\otimes t^{s_{\gamma}})X_{\boldsymbol{o}(\mu)}m =
	X^{w,z}_{\boldsymbol{o}(\mu)} (x^-_{\gamma+\mu_w+B}\otimes t^{s_{\mu_w +B}+s_\gamma}) (x^-_{\mu_z-B}\otimes t^{s_{\mu_z-B}})m.
\end{align*}
Since $m_{\gamma}<\ell$ and  \(\mu_w+\gamma+B=\beta_2\) we obtain $m_{\mu_w+B}+m_{\gamma}\leq \ell$ and therefore \(s_{\mu_w +B}+s_\gamma = s_{\beta_2} +1\). This forces that the above element is zero. \medskip

\textit{Case 1.1.2}: Here we assume
$$m_B+m_{\mu_z-B} \leq \ell\iff \mu_z\succeq \mu_z-B $$
We have with $m_{\gamma+\mu_w}+m_B\leq \ell$ and \eqref{req1}
\begin{align}
(x^-_{\gamma}\otimes t^{s_{\gamma}})X_{\boldsymbol{o}(\mu)}m &=\notag X^w_{\mathbf{o}(\mu)}(x^-_{\mu_w+\gamma}\otimes t^{s_{\mu_w+\gamma}})m&\\&\notag=
	X^{w,z}_{\boldsymbol{o}(\mu)}(x^-_{\beta_2}\otimes t^{s_{\beta_2}}) (x^-_{\mu_z-B}\otimes t^{s_{\mu_z-B}})m\\
	&=(x^-_{\beta_2}\otimes t^{s_{\beta_2}})X^{w,z}_{\boldsymbol{o}(\mu)}(x^-_{\mu_z-B}\otimes t^{s_{\mu_z-B}})m
  \label{eq:nminus-beta2}
	\\
	&\notag=(x^-_{\beta_2}\otimes t^{s_{\beta_2}}) (x^-_{\beta_1}\otimes t^{s_{\beta_1}})X^{w,u,z}_{\boldsymbol{o}(\mu)}(x^-_{\mu_u-\beta_1}\otimes t^{s_{\mu_u-\beta_1}})(x^-_{\mu_z-B}\otimes t^{s_{\mu_z-B}})m
 \end{align}
First we will argue why we also have $\mu_z-B\preceq \mu_w$ and we consider the cases \(\mu_z-B-\mu_w\in R^+\) or \(-(\mu_z-B-\mu_w)\in R^+\) separately. If the former holds, the converse assumption $\mu_z-B\succeq \mu_w$ would end in a contradiction:
$$m_{\mu_z}=m_{\mu_z-B}+m_{B}=m_{\mu_z-B-\mu_w}+m_{\mu_w}+m_{B}>\ell$$
If the latter holds, again the converse assumption $\mu_z-B\succeq \mu_w$ would end in a contradiction:
$$m_{\mu_w+B}+2\ell=m_{\mu_w}+m_{B}+\ell=m_{\mu_w-\mu_z+B}+m_{\mu_z-B}+m_{B}=m_{\mu_w-\mu_z+B}+m_{\mu_z}.$$
Hence we also have  \(\mu_w\succeq \mu_z-B\). If \(\mu_u \succeq \mu_z-B\) we are done (recall the maximal height choice).\medskip 

If instead \(\mu_u \preceq \mu_z-B\) we will differentiate the cases \(Y:=\mu_z-B-\mu_w-\beta_1 \in R^+\) and \(-Y \in R^+\).
In the former case our assumption \(\mu_u \preceq \mu_z-B\) gives \(m_{\mu_u}+m_{Y} \leq \ell\) and in the latter case \(m_{\mu_u+Y} +m_{-Y} > \ell \). However, in both cases we obtain \(\mu_u\succeq \mu_u+Y\). Summarizing, we have 
\begin{equation}\label{244as}\mu_w,\mu_z \succeq \mu_z-B \succeq \mu_u \succeq \mu_u+Y \end{equation}
As $s_{\boldsymbol{o}(\mu)}=s_{\boldsymbol{o}(\mu)}^{\mathrm{min}}$ we deduce further \(m_{\mu_z-B-Y}+m_{Y} \leq \ell\) (if $Y\in R^+$) and \(m_{\mu_z-B} +m_{-Y} > \ell \) (if $-Y\in R^+$) respectively.
We continue from equation \eqref{eq:nminus-beta2} and get
\begin{align*}
	(x^-_{\gamma}\otimes t^{s_{\gamma}})&X_{\boldsymbol{o}(\mu)}m =(x^-_{\beta_2}\otimes t^{s_{\beta_2}}) X^{w,z}_{\boldsymbol{o}(\mu)}(x^-_{\mu_z-B}\otimes t^{s_{\mu_z-B}})m\\
		&= (x^-_{\beta_2}\otimes t^{s_{\beta_2}}) X^{w,u,z}_{\boldsymbol{o}(\mu)}(x^-_{\mu_u+Y}\otimes t^{s_{\mu_u+Y}})(x^-_{\mu_z-B-Y}\otimes t^{s_{\mu_z-B-Y}})\\
  &= (x^-_{\beta_2}\otimes t^{s_{\beta_2}}) (x^-_{\beta_1}\otimes t^{s_{\beta_1}}) X^{w,u,z}_{\boldsymbol{o}(\mu)}(x^-_{\mu_u+Y}\otimes t^{s_{\mu_u+Y}})(x^-_{\mu_z-B-Y-\beta_1}\otimes t^{s_{\mu_z-B-Y-\beta_1}})
  .
\end{align*}
The second equality is obtained from the equations \(m_{\mu_u}+m_{Y} \leq \ell\) and \(m_{\mu_z-B-Y}+m_{Y} \leq \ell\) if \(Y\in R^+\) and the equations \(m_{\mu_u+Y} +m_{-Y}> \ell\) and \( m_{\mu_z-B} +m_{-Y}  > \ell \) if \(-Y \in R^+\).
The third equality is due to \(m_{\beta_1} =\ell\). It holds that \(\mu_z-B-Y-\beta_1=\mu_w\) and we are done with \eqref{244as}.

\textit{Case 1.2}: 
Let \(\varpi_i(h_{\gamma})=1\). 
Thus \(\mu(h_\gamma)=0\) which implies \(\mu_u(h_\gamma)=0\) for all \(u\) and we are done if \(\mu+\gamma \in \widetilde{\mathbf{R}}_{\lambda,i}^{\ell}\backslash\{0\} \). 
Otherwise there exists \(\beta_1\in R^+\) with \(m_{\beta_1}=\ell\) such that \(\varpi_i-\mu+\beta_1 \not\in W(\varpi_i)\) and \(\varpi_i-\mu-\gamma+\beta_1 \in W(\varpi_i)\) and we choose \(\beta_1\) of maximal height.
We assume first that \(\varpi_i(h_{\beta_1})=1\).
The roots must be arranged as follows
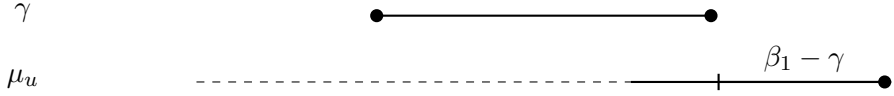
\begin{center}
\begin{tikzpicture}[baseline=2em]
  \matrix (m) [matrix of math nodes,row sep=1em,column sep=4em,minimum width=2em] {
    \gamma & \mbox{} &  \mbox{} & \mbox{} &  \mbox{} & \mbox{} &  \mbox{} \\
    \mu_u & \mbox{} &  \mbox{} & \mbox{} &  \mbox{} & \mbox{} &  \mbox{} \\
  }
 ;
 \node (a) at ($(m-2-4)!0.5!(m-2-5)$) {};
 \path[Circle-Circle,thick]
  (m-1-3.center) edge node [above] {\(\)} (m-1-5.center) 
  ;
 \path[dashed]
 (m-2-2.center) edge (a.center);
  \path[-Circle,thick]
  (a.center) edge node [above] {\(\)} (m-2-6.center);
 \path[|-,thick]
  (m-2-5.center) edge node [above] {\(\beta_1-\gamma\)} (m-2-6.center);
\end{tikzpicture}.
\end{center}
Note that \(\mu_u\) does not have to be unique, but there are at most two such roots. However, we assume in the rest of this case that $\mu_u$ is unique and omit the details otherwise; the strategy is similar to Case 1.1.2.\medskip 
Note that \(m_{\beta_1-\gamma}\neq \ell\) as \(\mu\in \widetilde{\mathbf{R}}_{\lambda,i}^{\ell}\backslash\{0\} \) and hence
\(m_{\beta_1-\gamma}+m_{\gamma}= \ell\).
As $s_{\boldsymbol{o}(\mu)}=s_{\boldsymbol{o}(\mu)}^{\mathrm{min}}$ we also have \(m_{\beta_1-\gamma}+m_{\mu_u-(\beta_1-\gamma)}\leq \ell \). Equivalently \(\mu_u-(\beta_1-\gamma) \preceq \mu_u\) and \eqref{req1} gives
\begin{align}
\label{eq:nminus-operation-fall12.1}
    (x^-_{\gamma}\otimes t^{s_{\gamma}})X_{\boldsymbol{o}(\mu)}m &=(x^-_{\beta_1}\otimes t^{s_{\beta_1}})X^{u}_{\boldsymbol{o}(\mu)} (x^-_{\mu_u-(\beta_1-\gamma)}\otimes t^{s_{\mu_u-(\beta_1-\gamma)}})m
    .
\end{align}
The right hand side is an element of \(\widetilde{\mathbf{R}}_{\lambda,i}^{\ell}\backslash\{0\}\) as \(\mu_u\) is unique and \(\beta_1\) is of maximal height. Also \ref{eq:nminus-operation-fall12.1} vanishes in the quotient since \(\mu_u-(\beta_1-\gamma)\preceq\mu_u\).

If \(\varpi_i(h_{\beta_1})=0\) we know from Case 1.1 that the term \((x^-_{\beta_1}\otimes t^{s_{\beta_1}})X_{\boldsymbol{o}(\mu)}m\) vanishes in the quotient and thus 
\begin{align*}
    (x^-_{\gamma}\otimes t^{s_{\gamma}})X_{\boldsymbol{o}(\mu)}m &=(x^-_{\beta_1}\otimes t^{s_{\beta_1}})(x^-_{\gamma-\beta_1}\otimes t^{s_{\gamma-\beta_1}})X_{\boldsymbol{o}(\mu)}m.
\end{align*}
Now setting $\gamma'=\gamma-\beta_1$ we will use the above arguments to finish this case also (note that $(\varpi_i-\mu)(h_{\gamma'})=1$ and $\varpi_i(h_{\gamma'})=1$). So either \(\mu+\gamma' \in \widetilde{\mathbf{R}}_{\lambda,i}^{\ell}\backslash\{0\}\) or there exists as above \(\beta_2\) with the aforementioned properties. If $\varpi_i(h_{\beta_2})=1$ we argue as above to show that   
$(x^-_{\gamma'}\otimes t^{s_{\gamma'}})X_{\boldsymbol{o}(\mu)}m$ vanishes in the quotient. Hence let $\varpi_i(h_{\beta_2})=0$ and we can write 
\begin{align*}
    (x^-_{\gamma'}\otimes t^{s_{\gamma'}})X_{\boldsymbol{o}(\mu)}m &=(x^-_{\beta_2}\otimes t^{s_{\beta_2}})(x^-_{\gamma'r-\beta_2}\otimes t^{s_{\gamma'-\beta_2}})X_{\boldsymbol{o}(\mu)}m.
\end{align*}
Repeating the above strategy gives the claim. 

\textit{Case 2}: Now suppose that $(\varpi_i-\mu)(h_{\gamma})=-1$.
Then
$$(\lambda+\varpi_i-\mu)(h_{\gamma})=(s_{\gamma}-1)\ell+m_{\gamma}-1$$
and we have to prove the following relations
$$(x^-_{\gamma}\otimes t^{s_{\gamma}})X_{\mathbf{o}(\mu)}m=0,\ \ m_{\gamma}>1$$
$$(x^-_{\gamma}\otimes t^{s_{\gamma}-1})^{m_{\gamma}}X_{\mathbf{o}(\mu)}m=0.$$

\textit{Case 2.1}: If $\varpi_i(h_{\gamma})=0$, then we have $\mu(h_{\gamma})=1$ and there exists a root $\mu_p$ with $(\mu_p,\gamma)=1$, i.e. $\mu_z-\gamma\in R^+$. Moreover, we note that $(\mu_w,\gamma)=0$ with all $w\neq p$ since $\mu\in \widetilde{\mathbf{R}}_{\lambda,i}^{\ell}\backslash\{0\}$. In particular $\gamma$ commutes with $X_{\mathbf{o}(\mu)}$ and the first relation is trivial and the second relations is implied by \eqref{req4} if $m_{\gamma}+m_{\mu_p-\gamma}>\ell$.
So we assume that $m_{\gamma}+m_{\mu_p-\gamma}\leq \ell$. In this case we obtain 
\begin{align*}(x^-_{\gamma}\otimes t^{s_{\gamma}-1})^{m_{\gamma}}X_{\mathbf{o}(\mu)}m&=X^p_{\mathbf{o}(\mu)}(x^-_{\gamma}\otimes t^{s_{\gamma}-1})^{m_{\gamma}}(x^-_{\mu_p}\otimes t^{s_{\gamma_p}})m&\\&=X^p_{\mathbf{o}(\mu)}(x^-_{\gamma}\otimes t^{s_{\gamma}-1})^{m_{\gamma}}\left[(x^-_{\gamma}\otimes t^{s_{\gamma}-1}),(x^-_{\mu_p-\gamma}\otimes t^{s_{\mu_p-\gamma}})\right]m&\\&=X^p_{\mathbf{o}(\mu)}(x^-_{\gamma}\otimes t^{s_{\gamma}-1})^{m_{\gamma}+1}(x^-_{\mu_p-\gamma}\otimes t^{s_{\mu_p-\gamma}})m
&\\& \hspace{0,5cm}-X^p_{\mathbf{o}(\mu)}(x^-_{\gamma}\otimes t^{s_{\gamma}-1})^{m_{\gamma}}(x^-_{\mu_p-\gamma}\otimes t^{s_{\mu_p-\gamma}})(x^-_{\gamma}\otimes t^{s_{\gamma}-1})m&\\&=X^p_{\mathbf{o}(\mu)}(x^-_{\gamma}\otimes t^{s_{\gamma}-1})^{m_{\gamma}+1}(x^-_{\mu_p-\gamma}\otimes t^{s_{\gamma_p-\gamma}})m&\\&\hspace{0,5cm}- m_{\gamma}(x^-_{\gamma}\otimes t^{s_{\gamma}-1})^{m_{\gamma}}X_{\mathbf{o}(\mu)}m
\end{align*}
and thus 
$$(x^-_{\gamma}\otimes t^{s_{\gamma}-1})^{m_{\gamma}}X_{\mathbf{o}(\mu)}m\in \mathbf{U}\cdot X^p_{\mathbf{o}(\mu)}(x^-_{\mu_p-\gamma}\otimes t^{s_{\mu_p-\gamma}})m$$
and $\mu_p-\gamma\preceq \mu_p$. So either the above element is contained in $\mathrm{HW}_i\cup\{0\}$ or we continue as in Lemma~\ref{nplus} (all assumptions hold) to write it as a linear combination of smaller elements in $\mathrm{HW}_i$. 

\textit{Case 2.2}: Let $\varpi_i(h_{\gamma})=1$, then we must have $\gamma=\mu_w$ for some $w$ in which case both equations are obvious (see also Example~\ref{ex11}) or we have two different roots $\mu_w$ and $\mu_p$ such that $\gamma-\mu_w\in R$ and $\gamma-\mu_p\in R$. The first relation is now immediate from \eqref{req3}. For the second relation we first consider the case when $\mu_w-\gamma,\mu_p-\gamma\in R^+$. We
can assume also that 
$$m_{\mu_p}+m_{\mu_w-\gamma}>\ell,\ \ m_{\gamma}+m_{\mu_w-\gamma}\leq \ell$$
$$m_{\mu_w}+m_{\mu_p-\gamma}>\ell,\ \ m_{\gamma}+m_{\mu_p-\gamma}\leq \ell$$
because otherwise $X_{\mathbf{o}(\mu)}m$ is proportional to an element which involves the factor $(x_{\gamma}^{-}\otimes t^{s_{\gamma}})$ (see relation \eqref{req1}) and the claim is immediate (see again Example~\ref{ex11}).
We obtain with the defining relations of $M_{i,\lambda}$:
\begin{align*}0&=(x_{\mu_p-\gamma}^-\otimes t^{s_{\mu_p-\gamma}})(x_{\mu_w-\gamma}^-\otimes t^{s_{\mu_w-\gamma}})(x^-_{\gamma}\otimes t^{s_{\gamma}-1})^{m_{\gamma}+2}m&\\&=
(x_{\mu_p-\gamma}^-\otimes t^{s_{\mu_p-\gamma}})(x_{\mu_w}^-\otimes t^{s_{\mu_w}})(x^-_{\gamma}\otimes t^{s_{\gamma}-1})^{m_{\gamma}+1}m&\\&
=(x^-_{\gamma}\otimes t^{s_{\gamma}-1})^{m_{\gamma}+1}(x_{\mu_w+\mu_p-\gamma}^-\otimes t^{s_{\mu_w+\mu_p-\gamma}})m+(x^-_{\gamma}\otimes t^{s_{\gamma}-1})^{m_{\gamma}}(x^-_{\mu_w}\otimes t^{s_{\gamma_w}})(x^-_{\mu_p}\otimes t^{s_{\mu_p}})m\end{align*}
with $$\mu_w+\mu_p-\gamma\leq \mu_p,\ \ \mu_w+\mu_p-\gamma\leq \mu_w$$
So either $X_{\mathbf{o}(\mu)}^{p,w}(x_{\mu_w+\mu_p-\gamma}^-\otimes t^{s_{\mu_w+\mu_p-\gamma}})m$ lies in $\mathrm{HW}_i\cup\{0\}$ or we continue as in Lemma~\ref{nplus} to write it as a linear combination of smaller elements. 

The second case considers $\gamma-\mu_p,\ \mu_w-\gamma\in R^+$. Again the element $X_{\mathbf{o}(\mu)}m$ is proportional to an element which involves an element $(x_{\gamma}^{-}\otimes t^{s_{\gamma}})$ (in this case we are done) or we can assume all the inequalities below 
$$m_{\mu_p+\mu_w-\gamma}+m_{\gamma-\mu_p}\leq \ell$$
$$m_{\mu_p}+m_{\mu_w-\gamma}>\ell,\ \ m_{\gamma}+m_{\mu_w-\gamma}\leq \ell$$
We get with \eqref{req4}
\begin{align*}0&=(x_{\mu_w-\gamma}^-\otimes t^{s_{\mu_w-\gamma}})(x^-_{\gamma}\otimes t^{s_{\gamma}-1})^{m_{\gamma}+1}(x_{\mu_p}^-\otimes t^{s_{\mu_p}})m&\\&
=(x^-_{\gamma}\otimes t^{s_{\gamma}-1})^{m_{\gamma}+1}(x_{\mu_w+\mu_p-\gamma}^-\otimes t^{s_{\mu_w+\mu_p-\gamma}})m+(x^-_{\gamma}\otimes t^{s_{\gamma}-1})^{m_{\gamma}}(x^-_{\mu_w}\otimes t^{s_{\mu_w}})(x^-_{\mu_p}\otimes t^{s_{\mu_z}})m\end{align*}
with $$\mu_w+\mu_p-\gamma\leq \mu_p,\ \ \mu_w+\mu_p-\gamma\leq \mu_w$$
Again we have that $X_{\mathbf{o}(\mu)}^{p,w}(x_{\mu_w+\mu_p-\gamma}^-\otimes t^{s_{\mu_w+\mu_p-\gamma}})m$ lies in $\mathrm{HW}_i\cup\{0\}$ or we continue as in Lemma~\ref{nplus} to write it as a linear combination of smaller elements. 

The third case considers $\gamma-\mu_w,\gamma-\mu_z\in R^+$. But this case follows immediately from \eqref{req5}.

\textit{Case 3}: 
Now suppose that $(\varpi_i-\mu)(h_{\gamma})=0$.
Then we have to prove the following relations
$$(x^-_{\gamma}\otimes t^{s_{\gamma}})X_{\mathbf{o}(\mu)}m=0$$
$$(x^-_{\gamma}\otimes t^{s_{\gamma}-1})^{m_{\gamma}+1}X_{\mathbf{o}(\mu)}m=0,\ \ m_{\gamma}<\ell$$

\textit{Case 3.1}: Let $\varpi_i(h_{\gamma})=0$ and note that there is nothing to show if $x^-_{\gamma}$ commutes with $X_{\mathbf{o}(\mu)}$. So assume that $\gamma$ does not commute with the unique root $\mu_w$. Since $\mu_j(h_{\gamma})=0$ and $(\mu_w,\gamma)=-1$ there must be a another root $\mu_z$ with $z\neq w$ and $(\mu_z,\gamma)=1$. This means $\mu_z-\gamma\in R^+$. So in particular with \eqref{req2} and \eqref{req3} we get
$$(x^-_{\mu_w+\gamma}\otimes t^{s_{\mu_w}+s_{\gamma}})(x^-_{\mu_z}\otimes t^{s_{\mu_z}})m=0$$
and the first relation follows immediately. A similar argument shows the second relation unless $s_{\mu_w}+s_{\gamma}=s_{\gamma+\mu_w}$, i.e. $m_{\mu_w}+m_{\gamma}>\ell$. 
But than we have 
$$(x^-_{\gamma_w}\otimes t^{s_{\gamma_w}})(x^-_{\gamma}\otimes t^{s_{\gamma}-1})^{m_{\gamma}+1}v_1=(x^-_{\gamma_w+\gamma}\otimes t^{s_{\gamma_w+\gamma}-1})(x^-_{\gamma}\otimes t^{s_{\gamma}-1})^{m_{\gamma}}v_1=0$$ and hence
with \eqref{req2} (for the first equation) 
$$0=(x^-_{\mu_w}\otimes t^{s_{\mu_w}})(x^-_{\gamma}\otimes t^{s_{\gamma}-1})^{m_{\gamma}+1}(x^-_{\mu_z}\otimes t^{s_{\mu_z}})X_{\mathbf{o}(\mu)}^{z,w}m=(x^-_{\gamma}\otimes t^{s_{\gamma}-1})^{m_{\gamma}+1}X_{\mathbf{o}(\mu)}m$$
$$+(x^-_{\gamma_w+\gamma}\otimes t^{s_{\gamma_w+\gamma}-1})(x^-_{\gamma}\otimes t^{s_{\gamma}-1})^{m_{\gamma}}(x^-_{\gamma_z}\otimes t^{s_{\gamma_z}})X_{\mathbf{o}(\mu)}^{z,w}m$$
and the claim follows from \eqref{req5}.

\textit{Case 3.2}: If $\varpi_i(h_{\gamma})=1$, i.e. $\mu(h_{\gamma})=1$ (in particular $\gamma\notin\{\mu_1,\dots,\mu_k\}$) we must have a unique root $\mu_w$ such that $(\mu_w,\gamma)=1$ since always $(\mu_z,\gamma)\geq 0$. So $\gamma-\mu_w\in R$ and we have with \eqref{req3} and \eqref{req5}
$$(x^-_{\gamma}\otimes t^{s_{\gamma}})(x^-_{\gamma_w}\otimes t^{s_{\mu_w}})m=0,\ \ (x^-_{\gamma}\otimes t^{s_{\gamma}-1})^{m_{\gamma}+1}(x^-_{\gamma_w}\otimes t^{s_{\mu_w}})m=0$$
This shows both relations in this case.
\end{proof}
\end{prop}
\subsection{}\label{section66} Now we finish the proof of parts (1) and (2) of Theorem~\ref{mainthmres1}. As a consequence of the discussions in Section~\ref{section6} we get the following character estimate
 \begin{flalign*}
 &&\mathrm{ch}_{\mathrm{gr}}\left(\mathbf{D}^{\ell}_{\lambda}*V(\varpi_i)\right)&\leq \mathrm{ch}_{\mathrm{gr}}(\mathbf{M}^{\ell}_{\lambda,i}) && \text{(by Lemma~\ref{smap})} \\
    && &= \mathrm{ch}_{\mathrm{gr}}(\mathbf{D}_{\lambda+\varpi_i}^{\ell})+\mathrm{ch}_{\mathrm{gr}}(K) && \text{(by Section~\ref{section42})} \\
    && &\leq \mathrm{ch}_{\mathrm{gr}}(\mathbf{D}_{\lambda+\varpi_i}^{\ell}) +\sum_{\mu\in \widetilde{\mathbf{R}}_{\lambda,i}^{\ell}\backslash\{0\}} \mathrm{ch}_{\mathrm{gr}}(\mathbf{D}_{\lambda+\varpi_i-\mu}^{\ell})\ q^{s_{\mathbf{o}(\mu)}} && \text{(by Proposition~\ref{mainprophel})}  \\
    && &
    =\sum_{\mu\in \widetilde{\mathbf{R}}_{\lambda,i}^{\ell}}\mathrm{ch}_{\mathrm{gr}}(\mathbf{D}_{\lambda+\varpi_i-\mu}^{\ell}) \ q^{s_{\mathbf{o}(\mu)}} 
  \end{flalign*}
Hence the proof is finished by a dimension estimate reversing the above estimate.  Writing \(\lambda=\lambda_0+\ell\lambda_1\) for some \(\lambda_0,\lambda_1\in P^+\) satisfying \(\lambda_0(h_\alpha) \leq \ell\) for all \(\alpha\in \Pi\) we get
\begin{align*}
 \dim\left(\mathbf{D}^{\ell}_{\lambda}*V(\varpi_i)\right)&=\dim(\mathbf{D}_{\lambda}^{\ell})\cdot \dim V(\varpi_i)&\\&=\dim(\mathbf{D}_{\ell\lambda_1}^{\ell})\cdot \dim(\mathbf{D}_{\lambda_0}^{\ell})\cdot \dim V(\varpi_i)\ \ \ \ \ (\text{by \eqref{stde}}) &\\&\geq \sum_{\mu\in \widetilde{\mathbf{R}}_{\lambda_0,i}^{\ell}}\dim(\mathbf{D}_{\ell\lambda_1}^{\ell})\cdot \dim(\mathbf{D}_{\lambda_0+\varpi_i-\mu}^{\ell}) \ \ \ \  (\text{by  Theorem~\ref{crysesti1})}&\\&
 =\sum_{\mu\in \widetilde{\mathbf{R}}_{\lambda,i}^{\ell}}\dim(\mathbf{D}_{\ell\lambda_1}^{\ell})\cdot \dim(\mathbf{D}_{\lambda_0+\varpi_i-\mu}^{\ell})\ \ \ \ \ (\text{since } \widetilde{\mathbf{R}}_{\lambda_0,i}^{\ell}=\widetilde{\mathbf{R}}_{\lambda,i}^{\ell})&\\&
 =\sum_{\mu\in \widetilde{\mathbf{R}}_{\lambda,i}^{\ell}}\dim(\mathbf{D}_{\lambda+\varpi_i-\mu}^{\ell}) \ \ \ \ \ \ \ \ \ \ \ \ \ \ \ \ \ \ \ \ \ \  (\text{by \eqref{stde}}) 
 \end{align*}
Thus $\mathrm{ch}_{\mathrm{gr}}\left(\mathbf{D}^{\ell}_{\lambda}*V(\varpi_i)\right)=\mathrm{ch}_{\mathrm{gr}}(\mathbf{M}^{\ell}_{\lambda,i})$ giving the statement of Theorem~\ref{mainthmres1}(2) and all maps in Proposition~\ref{mainprophel} are isomorphisms which completes the proof of Theorem~\ref{mainthmres1}(1).

\section{Appendix}\label{appendix}
We complete the proof of Lemma~\ref{poslambdai} and show in the rest of the article the transitivity of $\succeq$ by several case considerations. Let $\alpha,\beta,\alpha'\in R_i^+$ with $\alpha\succeq \beta\succeq \alpha'$. 

\textit{Case 1}: In this case we assume that $\alpha-\alpha'\in R$. 

\textit{Case 1.1}: Assume that $(\alpha\cap \alpha')-\beta\notin R$ and $\mathrm{supp}(\alpha\cap \alpha')\subseteq \mathrm{supp}(\beta)$. Then the only possible constellation of the roots is one of the following (or reflected):

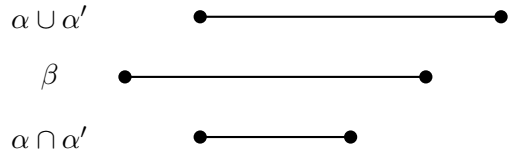
\begin{enumerate}
\item $$\begin{tikzpicture}[thick]
\draw (0,0) -- (4, 0);
\node (A) at (-2,0) [] {$\alpha\cup\alpha'$};
\node (R) at (0,0) [] {};
\node (P) at (4,0) [] {};
\node (E) at (0,-1.6) [] {};
\node (G) at (2, -1.6) [] {};
\fill (R) circle (2.5pt);
\fill (P) circle (2.5pt);
\draw (-1,-0.8) -- (3, -0.8);
\node (C) at (-2,-0.8) [circle] {$\beta$};
\node (D) at (-2,-1.6) [circle] {$\alpha\cap \alpha'$};
\node (RE) at (-1,-0.8) [] {};
\node (PE) at (3,-0.8) [] {};
\fill (RE) circle (2.5pt);
\fill (PE) circle (2.5pt);
\fill (E) circle (2.5pt);
\fill (G) circle (2.5pt);
\draw (0,-1.6) -- (2, -1.6);
\end{tikzpicture}$$
\item 
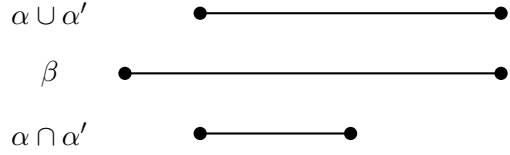
$$\begin{tikzpicture}[thick]
\draw (0,0) -- (4, 0);
\node (A) at (-2,0) [] {$\alpha\cup\alpha'$};
\node (R) at (0,0) [] {};
\node (P) at (4,0) [] {};
\node (E) at (0,-1.6) [] {};
\node (G) at (2, -1.6) [] {};
\fill (R) circle (2.5pt);
\fill (P) circle (2.5pt);
\draw (-1,-0.8) -- (4, -0.8);
\node (C) at (-2,-0.8) [circle] {$\beta$};
\node (D) at (-2,-1.6) [circle] {$\alpha\cap\alpha'$};
\node (RE) at (-1,-0.8) [] {};
\node (PE) at (4,-0.8) [] {};
\fill (RE) circle (2.5pt);
\fill (PE) circle (2.5pt);
\fill (E) circle (2.5pt);
\fill (G) circle (2.5pt);
\draw (0,-1.6) -- (2, -1.6);
\end{tikzpicture}$$
\item 
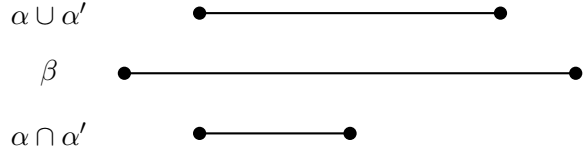
$$
\begin{tikzpicture}[thick]
\draw (0,0) -- (4, 0);
\node (A) at (-2,0) [] {$\alpha\cup\alpha'$};
\node (R) at (0,0) [] {};
\node (P) at (4,0) [] {};
\node (E) at (0,-1.6) [] {};
\node (G) at (2, -1.6) [] {};
\fill (R) circle (2.5pt);
\fill (P) circle (2.5pt);
\draw (-1,-0.8) -- (5, -0.8);
\node (C) at (-2,-0.8) [circle] {$\beta$};
\node (D) at (-2,-1.6) [circle] {$\alpha\cap\alpha'$};
\node (RE) at (-1,-0.8) [] {};
\node (PE) at (5,-0.8) [] {};
\fill (RE) circle (2.5pt);
\fill (PE) circle (2.5pt);
\fill (E) circle (2.5pt);
\fill (G) circle (2.5pt);
\draw (0,-1.6) -- (2, -1.6);
\end{tikzpicture}$$
\end{enumerate}

In the first case we have 
 $$m_{(\alpha\hspace{0,03cm}\cup\hspace{0,03cm}\alpha')\hspace{0,03cm}\cup\hspace{0,03cm} \beta}=m_{\alpha\hspace{0,03cm}\cup \hspace{0,03cm}\alpha'}+m_{(\alpha\hspace{0,03cm}\cup\hspace{0,03cm}\alpha')\hspace{0,05cm} \cup\hspace{0,03cm} \beta-(\alpha\hspace{0,03cm}\cup\hspace{0,03cm} \alpha')}-\ell\delta_{\alpha,\alpha\hspace{0,03cm}\cup\hspace{0,03cm}\alpha'}=m_{\beta}+m_{(\alpha\hspace{0,03cm}\cup\hspace{0,03cm}\alpha')\hspace{0,03cm}\cup \hspace{0,03cm} \beta-\beta}-\ell\delta_{\alpha',\alpha\hspace{0,03cm}\cup\hspace{0,03cm}\alpha'}$$
and 
\begin{equation}\label{111}m_{\beta}=m_{\alpha\hspace{0,03cm}\cap\hspace{0,03cm}\alpha'}+m_{\gamma_1}+m_{\gamma_2}-2\ell\delta_{\alpha,\alpha\hspace{0,03cm}\cap\hspace{0,03cm}\alpha'}.\end{equation} Now substituting the second equation in the first one gives 
$$m_{\alpha\hspace{0,03cm}\cup \hspace{0,03cm}\alpha'}+m_{(\alpha\hspace{0,03cm}\cup\hspace{0,03cm}\alpha')\hspace{0,05cm} \cup\hspace{0,03cm} \beta-(\alpha\hspace{0,03cm}\cup\alpha')}=m_{\alpha\hspace{0,03cm}\cap\hspace{0,03cm}\alpha'}+m_{\gamma_1}+m_{\gamma_2}-2\ell\delta_{\alpha,\alpha\hspace{0,03cm}\cup\hspace{0,03cm}\alpha'}+m_{(\alpha\hspace{0,03cm}\cup\hspace{0,03cm}\alpha')\hspace{0,03cm}\cup \hspace{0,03cm} \beta-\beta}$$
Since 
$$\gamma_1=(\alpha\cup\alpha')\cup \beta-(\alpha\cup\alpha'),\ \ \alpha\cup \alpha'-\alpha\cap \alpha'=\big((\alpha\cup\alpha')\cup \beta-\beta\big)+\gamma_2$$ we must have
$m_{\alpha\hspace{0,03cm}\cup \hspace{0,03cm}\alpha'}=m_{\alpha\hspace{0,03cm}\cap \hspace{0,03cm}\alpha'}+m_{\alpha\hspace{0,03cm}\cup \hspace{0,03cm}\alpha'-\alpha\hspace{0,03cm}\cap \hspace{0,03cm}\alpha'}-\ell\delta_{\alpha,\alpha\hspace{0,03cm}\cap\hspace{0,03cm}\alpha'} $
implying that $\alpha\succeq \alpha'$.\par
In the second case we have \eqref{111} and $m_{\beta}=m_{\alpha\hspace{0,03cm}\cup \hspace{0,03cm}\alpha'}+m_{\beta-(\alpha\hspace{0,03cm}\cup \hspace{0,03cm}\alpha')}-\ell \delta_{\alpha,\alpha\hspace{0,03cm}\cup \hspace{0,03cm}\alpha'}$ which gives (note that $(\alpha\cap \alpha')+\gamma_2=\alpha\cup\alpha'$ and $\beta-(\alpha\hspace{0,03cm}\cup \hspace{0,03cm}\alpha')=\gamma_1$)
$$m_{\alpha\hspace{0,03cm}\cap\hspace{0,03cm}\alpha'}+m_{\gamma_1}+m_{\gamma_2}-2\ell\delta_{\alpha,\alpha\hspace{0,03cm}\cap\hspace{0,03cm}\alpha'}=m_{\alpha\hspace{0,03cm}\cup \hspace{0,03cm}\alpha'}+m_{\beta-(\alpha\hspace{0,03cm}\cup \hspace{0,03cm}\alpha')}-\ell \delta_{\alpha,\alpha\hspace{0,03cm}\cup \hspace{0,03cm}\alpha'}$$
Therefore
$$m_{\alpha\hspace{0,03cm}\cap\hspace{0,03cm}\alpha'}+m_{\alpha\hspace{0,03cm}\cup\hspace{0,03cm}\alpha'-\alpha\hspace{0,03cm}\cap\hspace{0,03cm}\alpha'}-2\ell\delta_{\alpha,\alpha\hspace{0,03cm}\cap\hspace{0,03cm}\alpha'}=m_{\alpha\hspace{0,03cm}\cup \hspace{0,03cm}\alpha'}-\ell \delta_{\alpha,\alpha\hspace{0,03cm}\cup \hspace{0,03cm}\alpha'}$$
and we end in a contradiction.\par
In the third case we have \eqref{111} and $$m_{\beta}=m_{\alpha\hspace{0,03cm}\cup\hspace{0,03cm}\alpha'}+m_{\gamma_1}+m_{\gamma_2-((\alpha\hspace{0,03cm}\cap\hspace{0,03cm}\alpha')-(\alpha\hspace{0,03cm}\cup\hspace{0,03cm}\alpha'))}-2\ell \delta_{\alpha,\alpha\hspace{0,03cm}\cup\hspace{0,03cm}\alpha'}$$
and the same substitution as above shows once more that this case is impossible. 

\textit{Case 1.2}: Suppose that $(\alpha\cap\alpha')-\beta\notin R$ and $\mathrm{supp}(\beta)\subseteq \mathrm{supp}(\alpha\cap \alpha')$. This case gives
$$m_{\alpha\hspace{0,03cm}\cap\hspace{0,03cm}\alpha'}=m_{\gamma_1}+m_{\beta}+m_{\gamma_2}-2\ell \delta_{\alpha',\alpha\hspace{0,03cm}\cap\hspace{0,03cm}\alpha'},\ \ m_{\alpha\hspace{0,03cm}\cup\hspace{0,03cm}\alpha'}=m_{\gamma_1}+m_{\beta}+m_{\tilde{\gamma}_2}-2\ell\delta_{\alpha',\alpha\hspace{0,03cm}\cup\hspace{0,03cm}\alpha'}$$
where $$(\alpha\cap \alpha')-\beta=\gamma_1+\gamma_2,\ \ (\alpha\cup \alpha')-\beta=\tilde{\gamma}_1+\tilde{\gamma}_2$$
and we assume without loss of generality that $\tilde{\gamma}_1=\gamma_1$.
Hence $$m_{\alpha\hspace{0,03cm}\cup\hspace{0,03cm}\alpha'}=m_{\alpha\hspace{0,03cm}\cap\hspace{0,03cm}\alpha'}+m_{(\alpha\hspace{0,03cm}\cup\hspace{0,03cm}\alpha')-(\alpha\hspace{0,03cm}\cup\hspace{0,03cm}\alpha')}-\ell \delta_{\alpha,\alpha\hspace{0,03cm}\cap\hspace{0,03cm}\alpha'}$$ because otherwise  
$$m_{\alpha\hspace{0,03cm}\cup\hspace{0,03cm}\alpha'}=m_{\gamma_1}+m_{\beta}+m_{\gamma_2}+m_{(\alpha\hspace{0,03cm}\cup\hspace{0,03cm}\alpha')-(\alpha\hspace{0,03cm}\cup\hspace{0,03cm}\alpha')}-3\ell \delta_{\alpha',\alpha\hspace{0,03cm}\cap\hspace{0,03cm}\alpha'}=m_{\gamma_1}+m_{\beta}+m_{\tilde{\gamma}_2}-2\ell\delta_{\alpha',\alpha\hspace{0,03cm}\cup\hspace{0,03cm}\alpha'}$$
$$\Rightarrow m_{\gamma_2}+m_{(\alpha\hspace{0,03cm}\cup\hspace{0,03cm}\alpha')-(\alpha\hspace{0,03cm}\cup\hspace{0,03cm}\alpha')}-3\ell \delta_{\alpha',\alpha\hspace{0,03cm}\cap\hspace{0,03cm}\alpha'}=m_{\tilde{\gamma}_2}-2\ell\delta_{\alpha',\alpha\hspace{0,03cm}\cup\hspace{0,03cm}\alpha'}$$
which is a contradiction. 

\textit{Case 1.3}: Here we assume that $(\alpha\cap\alpha')\succeq \beta$ and $\beta \succeq (\alpha\cap\alpha')$ respectively is induced from Definition~\ref{poset1} (iii). We get
$$m_{(\alpha\hspace{0,03cm}\cap\hspace{0,03cm}\alpha')\hspace{0,03cm} \cup \hspace{0,03cm} \beta}=m_{(\alpha\hspace{0,03cm}\cap\hspace{0,03cm}\alpha')}+m_{(\alpha\hspace{0,03cm}\cap\hspace{0,03cm}\alpha')\hspace{0,03cm}\cup \hspace{0,03cm}\beta-(\alpha\hspace{0,03cm}\cap\hspace{0,03cm}\alpha')}-\ell \delta_{\alpha,\alpha\hspace{0,03cm}\cap\hspace{0,03cm}\alpha'}=m_{\beta}+m_{((\alpha\hspace{0,03cm}\cap\hspace{0,03cm}\alpha')\hspace{0,03cm}\cup\hspace{0,03cm} \beta)-\beta}-\ell \delta_{\alpha',\alpha\hspace{0,03cm}\cap\hspace{0,03cm}\alpha'}$$
If $(\alpha\cup\alpha')\succeq \beta$ and $\beta \succeq (\alpha\cup\alpha')$ respectively is also induced from Definition~\ref{poset1} (iii) we get 
$$m_{(\alpha\hspace{0,03cm}\cup\hspace{0,03cm}\alpha')\hspace{0,03cm}\cup \hspace{0,03cm}\beta}=m_{(\alpha\hspace{0,03cm}\cup\hspace{0,03cm}\alpha')}+m_{(\alpha\hspace{0,03cm}\cup\hspace{0,03cm}\alpha')\hspace{0,03cm}\cup \hspace{0,03cm}\beta-(\alpha\hspace{0,03cm}\cup\hspace{0,03cm}\alpha')}-\ell \delta_{\alpha,\alpha\hspace{0,03cm}\cup\hspace{0,03cm}\alpha'}=m_{\beta}+m_{((\alpha\hspace{0,03cm}\cup\hspace{0,03cm}\alpha')\hspace{0,03cm}\cup\hspace{0,03cm} \beta)-\beta}-\ell \delta_{\alpha',\alpha\hspace{0,03cm}\cup\hspace{0,03cm}\alpha'}$$
and hence by solving both equations for $m_{\beta}$:
$$m_{(\alpha\hspace{0,03cm}\cup\hspace{0,03cm}\alpha')}+m_{(\alpha\hspace{0,03cm}\cup\hspace{0,03cm}\alpha')\hspace{0,03cm}\cup \hspace{0,03cm}\beta-(\alpha\hspace{0,03cm}\cup\hspace{0,03cm}\alpha')}-2\ell \delta_{\alpha,\alpha\hspace{0,03cm}\cup\hspace{0,03cm}\alpha'}-m_{((\alpha\hspace{0,03cm}\cup\hspace{0,03cm}\alpha')\hspace{0,03cm}\cup\hspace{0,03cm} \beta)-\beta}+2\ell \delta_{\alpha',\alpha\hspace{0,03cm}\cup\hspace{0,03cm}\alpha'}
$$$$=m_{(\alpha\hspace{0,03cm}\cap\hspace{0,03cm}\alpha')}+m_{(\alpha\hspace{0,03cm}\cap\hspace{0,03cm}\alpha')\hspace{0,03cm}\cup \hspace{0,03cm}\beta-(\alpha\hspace{0,03cm}\cap\hspace{0,03cm}\alpha')}-m_{((\alpha\hspace{0,03cm}\cap\hspace{0,03cm}\alpha')\hspace{0,03cm}\cup\hspace{0,03cm} \beta)-\beta}$$
So if we assume that 
$$m_{\alpha\hspace{0,03cm}\cup\hspace{0,03cm}\alpha'}=m_{\alpha\hspace{0,03cm}\cap\hspace{0,03cm}\alpha'}+m_{(\alpha\hspace{0,03cm}\cup\hspace{0,03cm}\alpha')-(\alpha\hspace{0,03cm}\cup\hspace{0,03cm}\alpha')}-\ell \delta_{\alpha',\alpha\hspace{0,03cm}\cap\hspace{0,03cm}\alpha'}$$
then we get 
$$m_{(\alpha\hspace{0,03cm}\cup\hspace{0,03cm}\alpha')-(\alpha\hspace{0,03cm}\cup\hspace{0,03cm}\alpha')}+m_{(\alpha\hspace{0,03cm}\cup\hspace{0,03cm}\alpha')\hspace{0,03cm}\cup \hspace{0,03cm}\beta-(\alpha\hspace{0,03cm}\cup\hspace{0,03cm}\alpha')}-3\ell \delta_{\alpha,\alpha\hspace{0,03cm}\cup\hspace{0,03cm}\alpha'}-m_{((\alpha\hspace{0,03cm}\cup\hspace{0,03cm}\alpha')\hspace{0,03cm}\cup\hspace{0,03cm} \beta)-\beta}+2\ell \delta_{\alpha',\alpha\hspace{0,03cm}\cup\hspace{0,03cm}\alpha'}
$$$$=m_{(\alpha\hspace{0,03cm}\cap\hspace{0,03cm}\alpha')\hspace{0,03cm}\cup \hspace{0,03cm}\beta-(\alpha\hspace{0,03cm}\cap\hspace{0,03cm}\alpha')}-m_{((\alpha\hspace{0,03cm}\cap\hspace{0,03cm}\alpha')\hspace{0,03cm}\cup\hspace{0,03cm} \beta)-\beta}$$
which is impossible if $\alpha'=\alpha\cup \alpha'$ and if $\alpha=\alpha\cup \alpha'$ we have
\begin{align*}m_{(\alpha\hspace{0,03cm}\cup\hspace{0,03cm}\alpha')-(\alpha\hspace{0,03cm}\cup\hspace{0,03cm}\alpha')}&+m_{(\alpha\hspace{0,03cm}\cup\hspace{0,03cm}\alpha')\hspace{0,03cm}\cup \hspace{0,03cm}\beta-(\alpha\hspace{0,03cm}\cup\hspace{0,03cm}\alpha')}&\\&=m_{((\alpha\hspace{0,03cm}\cup\hspace{0,03cm}\alpha')\hspace{0,03cm}\cup\hspace{0,03cm} \beta)-\beta}+m_{(\alpha\hspace{0,03cm}\cap\hspace{0,03cm}\alpha')\hspace{0,03cm}\cup \hspace{0,03cm}\beta-(\alpha\hspace{0,03cm}\cap\hspace{0,03cm}\alpha')}-m_{((\alpha\hspace{0,03cm}\cap\hspace{0,03cm}\alpha')\hspace{0,03cm}\cup\hspace{0,03cm} \beta)-\beta}+3\ell\end{align*}
which is again a contradiction. If $(\alpha\cup \alpha')-\beta\in R$, then the roots are as follows 
$$
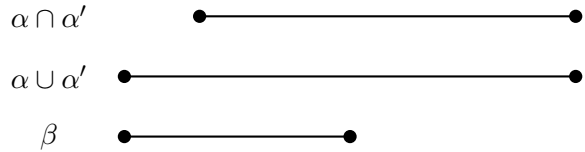
\begin{tikzpicture}[thick]
\draw (0,0) -- (5, 0);
\node (A) at (-2,0) [] {$\alpha\cap\alpha'$};
\node (R) at (0,0) [] {};
\node (P) at (5,0) [] {};
\node (E) at (-1,-1.6) [] {};
\node (G) at (2, -1.6) [] {};
\fill (R) circle (2.5pt);
\fill (P) circle (2.5pt);
\draw (-1,-0.8) -- (5, -0.8);
\node (C) at (-2,-0.8) [circle] {$\alpha\cup \alpha'$};
\node (D) at (-2,-1.6) [circle] {$\beta$};
\node (RE) at (-1,-0.8) [] {};
\node (PE) at (5,-0.8) [] {};
\fill (RE) circle (2.5pt);
\fill (PE) circle (2.5pt);
\fill (E) circle (2.5pt);
\fill (G) circle (2.5pt);
\draw (-1,-1.6) -- (2, -1.6);
\end{tikzpicture}
$$
and we have (note that $(\alpha\hspace{0,03cm}\cup\hspace{0,03cm}\alpha')\hspace{0,03cm}\cap\hspace{0,03cm} \beta=\beta)$
$$m_{\alpha\hspace{0,03cm}\cup\hspace{0,03cm}\alpha'}=m_{\beta}+m_{\alpha\hspace{0,03cm}\cup\hspace{0,03cm}\alpha'-\beta}-\ell \delta_{\alpha',\alpha\hspace{0,03cm}\cup\hspace{0,03cm}\alpha'}.$$
Together with 
$$m_{(\alpha\hspace{0,03cm}\cap\hspace{0,03cm}\alpha')\cup \beta}=m_{(\alpha\hspace{0,03cm}\cap\hspace{0,03cm}\alpha')}+m_{(\alpha\hspace{0,03cm}\cap\hspace{0,03cm}\alpha')\hspace{0,03cm}\cup \hspace{0,03cm}\beta-(\alpha\hspace{0,03cm}\cap\hspace{0,03cm}\alpha')}-\ell \delta_{\alpha,\alpha\hspace{0,03cm}\cap\hspace{0,03cm}\alpha'}=m_{\beta}+m_{((\alpha\hspace{0,03cm}\cap\hspace{0,03cm}\alpha')\hspace{0,03cm}\cup\hspace{0,03cm} \beta)-\beta}-\ell \delta_{\alpha',\alpha\hspace{0,03cm}\cap\hspace{0,03cm}\alpha'}$$
$$=m_{\alpha\hspace{0,03cm}\cup\hspace{0,03cm}\alpha'}-m_{\alpha\hspace{0,03cm}\cup\hspace{0,03cm}\alpha'-\beta}+\ell \delta_{\alpha',\alpha\hspace{0,03cm}\cup\hspace{0,03cm}\alpha'}+m_{((\alpha\hspace{0,03cm}\cap\hspace{0,03cm}\alpha')\hspace{0,03cm}\cup\hspace{0,03cm} \beta)-\beta}-\ell \delta_{\alpha',\alpha\hspace{0,03cm}\cap\hspace{0,03cm}\alpha'}$$
we obtain (note that $(\alpha\hspace{0,03cm}\cup\hspace{0,03cm}\alpha')-\beta=((\alpha\hspace{0,03cm}\cap\hspace{0,03cm}\alpha')\hspace{0,03cm}\cup\hspace{0,03cm} \beta)-\beta$)
$$m_{(\alpha\hspace{0,03cm}\cap\hspace{0,03cm}\alpha')}+m_{(\alpha\hspace{0,03cm}\cup\hspace{0,03cm}\alpha')-(\alpha\hspace{0,03cm}\cap\hspace{0,03cm}\alpha')}-2\ell \delta_{\alpha,\alpha\hspace{0,03cm}\cap\hspace{0,03cm}\alpha'}=m_{\alpha\hspace{0,03cm}\cup\hspace{0,03cm}\alpha'}-\ell \delta_{\alpha',\alpha\hspace{0,03cm}\cap\hspace{0,03cm}\alpha'}$$
which is again a contradiction. So the only possible constellation is 
$$
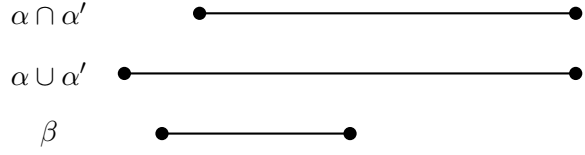
\begin{tikzpicture}[thick]
\draw (0,0) -- (5, 0);
\node (A) at (-2,0) [] {$\alpha\cap\alpha'$};
\node (R) at (0,0) [] {};
\node (P) at (5,0) [] {};
\node (E) at (-0.5,-1.6) [] {};
\node (G) at (2, -1.6) [] {};
\fill (R) circle (2.5pt);
\fill (P) circle (2.5pt);
\draw (-1,-0.8) -- (5, -0.8);
\node (C) at (-2,-0.8) [circle] {$\alpha\cup \alpha'$};
\node (D) at (-2,-1.6) [circle] {$\beta$};
\node (RE) at (-1,-0.8) [] {};
\node (PE) at (5,-0.8) [] {};
\fill (RE) circle (2.5pt);
\fill (PE) circle (2.5pt);
\fill (E) circle (2.5pt);
\fill (G) circle (2.5pt);
\draw (-0.5,-1.6) -- (2, -1.6);
\end{tikzpicture}
$$
This means
$$m_{(\alpha\hspace{0,03cm}\cup\hspace{0,03cm}\alpha')}=m_{\beta}+m_{\gamma_1}+m_{\gamma_2}-2\ell \delta_{\alpha',\alpha\hspace{0,03cm}\cup\hspace{0,03cm}\alpha'}$$
and 
$$m_{(\alpha\hspace{0,03cm}\cap\hspace{0,03cm}\alpha'\hspace{0,03cm})\cup \hspace{0,03cm}\beta}=m_{(\alpha\hspace{0,03cm}\cap\hspace{0,03cm}\alpha')}+m_{(\alpha\hspace{0,03cm}\cap\hspace{0,03cm}\alpha')\hspace{0,03cm}\cup \hspace{0,03cm}\beta-(\alpha\hspace{0,03cm}\cap\hspace{0,03cm}\alpha')}-\ell \delta_{\alpha,\alpha\hspace{0,03cm}\cap\hspace{0,03cm}\alpha'}=m_{\beta}+m_{\gamma_2}-\ell \delta_{\alpha',\alpha\hspace{0,03cm}\cap\hspace{0,03cm}\alpha'}
$$
$$=m_{(\alpha\hspace{0,03cm}\cup\hspace{0,03cm}\alpha')}-m_{\gamma_1}+2\ell \delta_{\alpha',\alpha\hspace{0,03cm}\cup\hspace{0,03cm}\alpha'}-\ell \delta_{\alpha',\alpha\hspace{0,03cm}\cap\hspace{0,03cm}\alpha'}$$
which is once more a contradiction.

\textit{Case 1.4}: We assume in this case that $(\alpha\hspace{0,03cm}\cap\hspace{0,03cm}\alpha')-\beta\in R$. We obtain 
\begin{equation}\label{4r3}m_{(\alpha\hspace{0,03cm}\cap\hspace{0,03cm}\alpha')\hspace{0,03cm}\cup\hspace{0,03cm} \beta}=m_{(\alpha\hspace{0,03cm}\cap\hspace{0,03cm}\alpha')\hspace{0,03cm}\cap \hspace{0,03cm}\beta}+m_{(\alpha\hspace{0,03cm}\cap\hspace{0,03cm}\alpha')\hspace{0,03cm}\cup \hspace{0,03cm} \beta-(\alpha\hspace{0,03cm}\cap\hspace{0,03cm}\alpha')\hspace{0,03cm}\cap \hspace{0,03cm} \beta}-\ell\delta_{X,(\alpha\hspace{0,03cm}\cap\hspace{0,03cm}\alpha')\hspace{0,03cm}\cap \hspace{0,03cm}\hspace{0,03cm} \beta}\end{equation}
where $X=\alpha$ if $\alpha=\alpha\cap \alpha'$ and $X=\beta$ otherwise. If additionally  $(\alpha\hspace{0,03cm}\cup\hspace{0,03cm}\alpha')-\beta\in R$ we also have \begin{equation}\label{4r4}m_{(\alpha\hspace{0,03cm}\cup\hspace{0,03cm}\alpha')\hspace{0,03cm}\cup \hspace{0,03cm}\beta}=m_{(\alpha\hspace{0,03cm}\cup\hspace{0,03cm}\alpha')\hspace{0,03cm}\cap \hspace{0,03cm}\beta}+m_{(\alpha\hspace{0,03cm}\cup\hspace{0,03cm}\alpha')\hspace{0,03cm}\cup \hspace{0,03cm}\beta-(\alpha\hspace{0,03cm}\cup\hspace{0,03cm}\alpha')\hspace{0,03cm}\cap\hspace{0,03cm} \beta}-\ell\delta_{Y,(\alpha\hspace{0,03cm}\cup\hspace{0,03cm}\alpha')\hspace{0,03cm}\cap \hspace{0,03cm}\beta}\end{equation}
where $Y=\alpha$ if $\alpha=\alpha\cup \alpha'$ and $Y=\beta$ otherwise. Note that we have either of the following situations:
\begin{enumerate}
\item   $$
\begin{tikzpicture}[thick]
\draw (-1,0) -- (1, 0);
\node (A) at (-2,0) [] {$\beta$};
\node (R) at (-1,0) [] {};
\node (P) at (1,0) [] {};
\node (E) at (-1,-1.6) [] {};
\node (G) at (2, -1.6) [] {};
\fill (R) circle (2.5pt);
\fill (P) circle (2.5pt);
\draw (-1,-0.8) -- (5, -0.8);
\node (C) at (-2,-0.8) [circle] {$\alpha\cup\alpha'$};
\node (D) at (-2,-1.6) [circle] {$\alpha\cap \alpha'$};
\node (RE) at (-1,-0.8) [] {};
\node (PE) at (5,-0.8) [] {};
\fill (RE) circle (2.5pt);
\fill (PE) circle (2.5pt);
\fill (E) circle (2.5pt);
\fill (G) circle (2.5pt);
\draw (-1,-1.6) -- (2, -1.6);
\end{tikzpicture}$$
    \item $$
\begin{tikzpicture}[thick]
\draw (-1,0) -- (4, 0);
\node (A) at (-2,0) [] {$\beta$};
\node (R) at (-1,0) [] {};
\node (P) at (4,0) [] {};
\node (E) at (-1,-1.6) [] {};
\node (G) at (2, -1.6) [] {};
\fill (R) circle (2.5pt);
\fill (P) circle (2.5pt);
\draw (-1,-0.8) -- (5, -0.8);
\node (C) at (-2,-0.8) [circle] {$\alpha\cup\alpha'$};
\node (D) at (-2,-1.6) [circle] {$\alpha\cap \alpha'$};
\node (RE) at (-1,-0.8) [] {};
\node (PE) at (5,-0.8) [] {};
\fill (RE) circle (2.5pt);
\fill (PE) circle (2.5pt);
\fill (E) circle (2.5pt);
\fill (G) circle (2.5pt);
\draw (-1,-1.6) -- (2, -1.6);
\end{tikzpicture}$$

 \item   $$
\begin{tikzpicture}[thick]
\draw (-1,0) -- (6, 0);
\node (A) at (-2,0) [] {$\beta$};
\node (R) at (-1,0) [] {};
\node (P) at (6,0) [] {};
\node (E) at (-1,-1.6) [] {};
\node (G) at (2, -1.6) [] {};
\fill (R) circle (2.5pt);
\fill (P) circle (2.5pt);
\draw (-1,-0.8) -- (5, -0.8);
\node (C) at (-2,-0.8) [circle] {$\alpha\cup\alpha'$};
\node (D) at (-2,-1.6) [circle] {$\alpha\cap \alpha'$};
\node (RE) at (-1,-0.8) [] {};
\node (PE) at (5,-0.8) [] {};
\fill (RE) circle (2.5pt);
\fill (PE) circle (2.5pt);
\fill (E) circle (2.5pt);
\fill (G) circle (2.5pt);
\draw (-1,-1.6) -- (2, -1.6);
\end{tikzpicture}$$
\end{enumerate}
In the first case \eqref{4r3} and \eqref{4r4} transform to 

$$m_{\alpha\hspace{0,03cm}\cap\hspace{0,03cm}\alpha'}=m_{\beta}+m_{\alpha\hspace{0,03cm}\cap\hspace{0,03cm}\alpha'- \beta}-\ell\delta_{\alpha,\alpha\hspace{0,03cm}\cup\hspace{0,03cm}\alpha'}$$
and 
$$m_{\alpha\hspace{0,03cm}\cup\hspace{0,03cm}\alpha'}=m_{\beta}+m_{\alpha\hspace{0,03cm}\cup\hspace{0,03cm}\alpha'-\beta}-\ell\delta_{\alpha,\alpha\hspace{0,03cm}\cap\hspace{0,03cm}\alpha'}$$
Thus we get the desired property
$$m_{\alpha\hspace{0,03cm}\cup\hspace{0,03cm}\alpha'}=m_{\alpha\hspace{0,03cm}\cap\hspace{0,03cm}\alpha'}+m_{\alpha\hspace{0,03cm}\cup\hspace{0,03cm}\alpha'-\beta}-\ell\delta_{\alpha,\alpha\hspace{0,03cm}\cap\hspace{0,03cm}\alpha'}-m_{\alpha\hspace{0,03cm}\cap\hspace{0,03cm}\alpha'- \beta}+\ell\delta_{\alpha,\alpha\hspace{0,03cm}\cup\hspace{0,03cm}\alpha'}=m_{\alpha\hspace{0,03cm}\cap\hspace{0,03cm}\alpha'}+m_{(\alpha\hspace{0,03cm}\cup\hspace{0,03cm}\alpha')-(\alpha\hspace{0,03cm}\cap\hspace{0,03cm}\alpha')}-\ell\delta_{\alpha,\alpha\hspace{0,03cm}\cap\hspace{0,03cm}\alpha'}.$$ 
In the second case \eqref{4r3} and \eqref{4r4} transform to

$$m_{\beta}=m_{(\alpha\hspace{0,03cm}\cap\hspace{0,03cm}\alpha')}+m_{ \beta-(\alpha\hspace{0,03cm}\cap\hspace{0,03cm}\alpha')}-\ell\delta_{\alpha,\alpha\hspace{0,03cm}\cap\hspace{0,03cm}\alpha'}$$
and
$$m_{\alpha\hspace{0,03cm}\cup\hspace{0,03cm}\alpha'}=m_{\beta}+m_{\alpha\hspace{0,03cm}\cup\hspace{0,03cm}\alpha'-\beta}-\ell\delta_{\alpha,\alpha\hspace{0,03cm}\cap\hspace{0,03cm}\alpha'}$$
Hence we have once more the desired property
$$m_{\alpha\hspace{0,03cm}\cup\hspace{0,03cm}\alpha'}=m_{(\alpha\hspace{0,03cm}\cap\hspace{0,03cm}\alpha')}+m_{ \beta-(\alpha\hspace{0,03cm}\cap\hspace{0,03cm}\alpha')}+m_{\alpha\hspace{0,03cm}\cup\hspace{0,03cm}\alpha'-\beta}-2\ell\delta_{\alpha,\alpha\hspace{0,03cm}\cap\hspace{0,03cm}\alpha'}=m_{\alpha\hspace{0,03cm}\cap\hspace{0,03cm}\alpha'}+m_{(\alpha\hspace{0,03cm}\cup\hspace{0,03cm}\alpha')-(\alpha\hspace{0,03cm}\cap\hspace{0,03cm}\alpha')}-\ell\delta_{\alpha,\alpha\hspace{0,03cm}\cap\hspace{0,03cm}\alpha'}$$
In the last case we get 
$$m_{ \beta}=m_{(\alpha\hspace{0,03cm}\cap\hspace{0,03cm}\alpha')}+m_{\beta-(\alpha\hspace{0,03cm}\cap\hspace{0,03cm}\alpha')}-\ell\delta_{\alpha,\alpha\hspace{0,03cm}\cap\hspace{0,03cm}\alpha'}$$ 
and
$$m_{\beta}=m_{(\alpha\hspace{0,03cm}\cup\hspace{0,03cm}\alpha')}+m_{\beta-(\alpha\hspace{0,03cm}\cup\hspace{0,03cm}\alpha')}-\ell\delta_{\alpha,\alpha\hspace{0,03cm}\cup\hspace{0,03cm}\alpha'}$$
This gives 
$$m_{(\alpha\hspace{0,03cm}\cup\hspace{0,03cm}\alpha')}=m_{(\alpha\hspace{0,03cm}\cap\hspace{0,03cm}\alpha')}+m_{\beta-(\alpha\hspace{0,03cm}\cap\hspace{0,03cm}\alpha')}-m_{\beta-(\alpha\hspace{0,03cm}\cup\hspace{0,03cm}\alpha')}-\ell\delta_{\alpha,\alpha\hspace{0,03cm}\cap\hspace{0,03cm}\alpha'}+\ell\delta_{\alpha,\alpha\hspace{0,03cm}\cup\hspace{0,03cm}\alpha'}
$$
$$=m_{(\alpha\hspace{0,03cm}\cap\hspace{0,03cm}\alpha')}+m_{(\alpha\hspace{0,03cm}\cup\hspace{0,03cm}\alpha')-(\alpha\hspace{0,03cm}\cap\hspace{0,03cm}\alpha')}-\ell\delta_{\alpha,\alpha\hspace{0,03cm}\cap\hspace{0,03cm}\alpha'}$$
and we are done. \par 
So we can suppose in the rest of this case that $(\alpha\hspace{0,03cm}\cup\hspace{0,03cm}\alpha')-\beta\notin R$. We consider two possibilities 
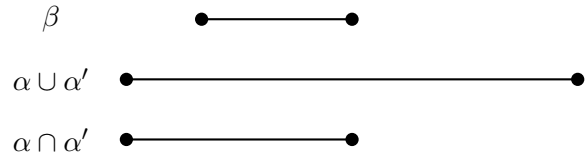
\begin{enumerate}
    \item $$\begin{tikzpicture}[thick]
\draw (0,0) -- (2, 0);
\node (A) at (-2,0) [] {$\beta$};
\node (R) at (0,0) [] {};
\node (P) at (2,0) [] {};
\node (E) at (-1,-1.6) [] {};
\node (G) at (2, -1.6) [] {};
\fill (R) circle (2.5pt);
\fill (P) circle (2.5pt);
\draw (-1,-0.8) -- (5, -0.8);
\node (C) at (-2,-0.8) [circle] {$\alpha\cup\alpha'$};
\node (D) at (-2,-1.6) [circle] {$\alpha\cap \alpha'$};
\node (RE) at (-1,-0.8) [] {};
\node (PE) at (5,-0.8) [] {};
\fill (RE) circle (2.5pt);
\fill (PE) circle (2.5pt);
\fill (E) circle (2.5pt);
\fill (G) circle (2.5pt);
\draw (-1,-1.6) -- (2, -1.6);
\end{tikzpicture}
$$
\item  
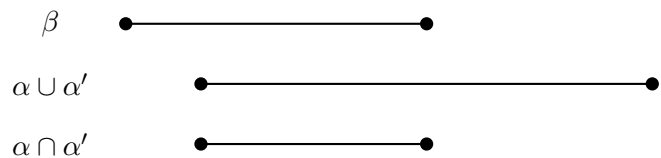
$$\begin{tikzpicture}[thick]
\draw (-2,0) -- (2, 0);
\node (A) at (-3,0) [] {$\beta$};
\node (R) at (-2,0) [] {};
\node (P) at (2,0) [] {};
\node (E) at (-1,-1.6) [] {};
\node (G) at (2, -1.6) [] {};
\fill (R) circle (2.5pt);
\fill (P) circle (2.5pt);
\draw (-1,-0.8) -- (5, -0.8);
\node (C) at (-3,-0.8) [circle] {$\alpha\cup\alpha'$};
\node (D) at (-3,-1.6) [circle] {$\alpha\cap \alpha'$};
\node (RE) at (-1,-0.8) [] {};
\node (PE) at (5,-0.8) [] {};
\fill (RE) circle (2.5pt);
\fill (PE) circle (2.5pt);
\fill (E) circle (2.5pt);
\fill (G) circle (2.5pt);
\draw (-1,-1.6) -- (2, -1.6);
\end{tikzpicture}
$$
\end{enumerate}
In the first case 
$$m_{\alpha\hspace{0,03cm}\cap\hspace{0,03cm}\alpha'}=m_{ \beta}+m_{(\alpha\hspace{0,03cm}\cap\hspace{0,03cm}\alpha')-\beta}-\ell\delta_{\alpha,\alpha\hspace{0,03cm}\cup\hspace{0,03cm}\alpha'}$$
and
$$m_{\alpha\hspace{0,03cm}\cup\hspace{0,03cm}\alpha'}=m_{\beta}+m_{\gamma_1}+m_{\gamma_2}-2\ell \delta_{\alpha,\alpha\hspace{0,03cm}\cap\hspace{0,03cm}\alpha'}.$$
This gives (note that $\gamma_1=(\alpha\hspace{0,03cm}\cap\hspace{0,03cm}\alpha')-\beta$,\ $\gamma_2=(\alpha\hspace{0,03cm}\cup\hspace{0,03cm}\alpha')-(\alpha\hspace{0,03cm}\cap\hspace{0,03cm}\alpha'$))
$$m_{\alpha\hspace{0,03cm}\cup\hspace{0,03cm}\alpha'}=m_{\alpha\hspace{0,03cm}\cap\hspace{0,03cm}\alpha'}+m_{\gamma_2}-2\ell \delta_{\alpha,\alpha\hspace{0,03cm}\cap\hspace{0,03cm}\alpha'}+\ell\delta_{\alpha,\alpha\hspace{0,03cm}\cup\hspace{0,03cm}\alpha'}$$
which is a contradiction. \par 
In the second case
$$m_{\beta}=m_{(\alpha\hspace{0,03cm}\cap\hspace{0,03cm}\alpha')}+m_{\beta-(\alpha\hspace{0,03cm}\cap\hspace{0,03cm}\alpha')}-\ell\delta_{\alpha,\alpha\hspace{0,03cm}\cap\hspace{0,03cm}\alpha'}$$
 and 
$$m_{(\alpha\hspace{0,03cm}\cup\hspace{0,03cm} \alpha')\hspace{0,03cm}\cup\hspace{0,03cm} \beta}=m_{\alpha\hspace{0,03cm}\cup \hspace{0,03cm}\alpha'}+m_{(\alpha\hspace{0,03cm}\cup\hspace{0,03cm} \alpha')\hspace{0,03cm}\cup \hspace{0,03cm}\beta-(\alpha\hspace{0,03cm}\cup\hspace{0,03cm} \alpha')}-\ell \delta_{\alpha,\alpha\hspace{0,03cm}\cup \hspace{0,03cm}\alpha'}=m_{\beta}+m_{(\alpha\hspace{0,03cm}\cup \hspace{0,03cm}\alpha')\hspace{0,03cm}\cup\hspace{0,03cm} \beta-\beta}-\ell \delta_{\alpha',\alpha\hspace{0,03cm}\cup \hspace{0,03cm}\alpha'}$$
This gives (note that $\beta-(\alpha\hspace{0,03cm}\cap\hspace{0,03cm}\alpha')=(\alpha\hspace{0,03cm}\cup\hspace{0,03cm} \alpha')\cup\hspace{0,03cm} \beta-(\alpha\hspace{0,03cm}\cup\hspace{0,03cm} \alpha')$)
$$m_{\alpha\hspace{0,03cm}\cup \hspace{0,03cm}\alpha'}=m_{(\alpha\hspace{0,03cm}\cap\hspace{0,03cm}\alpha')}+\ell \delta_{\alpha,\alpha\hspace{0,03cm}\cup \hspace{0,03cm}\alpha'}-2\ell\delta_{\alpha,\alpha\hspace{0,03cm}\cap\hspace{0,03cm}\alpha'}+m_{(\alpha\hspace{0,03cm}\cup \hspace{0,03cm}\alpha')-(\alpha\hspace{0,03cm}\cap \hspace{0,03cm}\alpha')}$$
which is again a contradiction. 

\textit{Case 2}: In this case we assume $\alpha-\alpha'\notin R$ and $\mathrm{supp}(\alpha')\subseteq \mathrm{supp}(\alpha)$ or $\mathrm{supp}(\alpha)\subseteq \mathrm{supp}(\alpha')$. We set in the rest of this case for simplicity $X=\alpha\cup\alpha'$ and $Y=\alpha\cap \alpha'$. We first consider the constellation of the roots as follows

$$
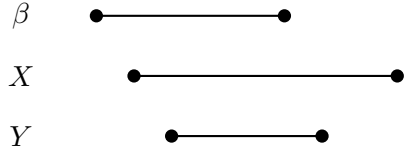
\begin{tikzpicture}[thick]
\draw (-1,0) -- (1.5, 0);
\node (A) at (-2,0) [] {$\beta$};
\node (R) at (-1,0) [] {};
\node (P) at (1.5,0) [] {};
\node (E) at (0,-1.6) [] {};
\node (G) at (2, -1.6) [] {};
\fill (R) circle (2.5pt);
\fill (P) circle (2.5pt);
\draw (-0.5,-0.8) -- (3, -0.8);
\node (C) at (-2,-0.8) [circle] {$X$};
\node (D) at (-2,-1.6) [circle] {$Y$};
\node (RE) at (-0.5,-0.8) [] {};
\node (PE) at (3,-0.8) [] {};
\fill (RE) circle (2.5pt);
\fill (PE) circle (2.5pt);
\fill (E) circle (2.5pt);
\fill (G) circle (2.5pt);
\draw (0,-1.6) -- (2, -1.6);
\end{tikzpicture}
$$
We have 
$$m_{\alpha\hspace{0,03cm}\cup\hspace{0,03cm} \beta}=m_{\alpha}+m_{\alpha\hspace{0,05cm} \cup\hspace{0,03cm} \beta-\alpha}-\ell=m_{\beta}+m_{\alpha\hspace{0,03cm}\cup \hspace{0,03cm} \beta-\beta}.$$
$$m_{\alpha'\hspace{0,03cm}\cup\hspace{0,03cm} \beta}=m_{\alpha'}+m_{\alpha'\hspace{0,05cm} \cup\hspace{0,03cm} \beta-\alpha'}=m_{\beta}+m_{\alpha'\hspace{0,03cm}\cup \hspace{0,03cm} \beta-\beta}-\ell.$$ 
Note that $\gamma_2:=X\hspace{0,03cm}\cup\hspace{0,03cm} \beta-Y\hspace{0,03cm}\cup\hspace{0,03cm} \beta$ is a positive root and we get from a straightforward calculation
$$m_{\gamma_2}+m_{Y\hspace{0,03cm}\cup\hspace{0,03cm} \beta-\beta}=m_{X\hspace{0,03cm}\cup\hspace{0,03cm} \beta-\beta}+\delta_{X,\alpha}\ell$$
Similarly $\gamma_1:=(Y\hspace{0,03cm}\cup\hspace{0,03cm} \beta-Y)-(X\hspace{0,03cm}\cup\hspace{0,03cm} \beta-X)$ is a positive root 
and $X-Y=\gamma_1+\gamma_2$. We get 
$$m_{\gamma_2}-2\delta_{X,\alpha'}(\ell+m_{\gamma_2})=m_{\alpha}-m_{\alpha'}-\ell+m_{\alpha\hspace{0,05cm} \cup\hspace{0,03cm} \beta-\alpha}-m_{\alpha'\hspace{0,05cm} \cup\hspace{0,03cm} \beta-\alpha'}=m_{\alpha}-m_{\alpha'}-m_{\gamma_1}$$
which forces $X=\alpha$ and thus
$$m_{\gamma_2}=m_{\alpha}-m_{\alpha'}-\ell+m_{\alpha\hspace{0,05cm} \cup\hspace{0,03cm} \beta-\alpha}-m_{\alpha'\hspace{0,05cm} \cup\hspace{0,03cm} \beta-\alpha'}=m_{\alpha}-m_{\alpha'}-m_{\gamma_1}\Rightarrow m_{\gamma_1}+m_{\gamma_2}+m_{\alpha'}=m_{\alpha}$$
So we are done in this case. \par

Now suppose

$$
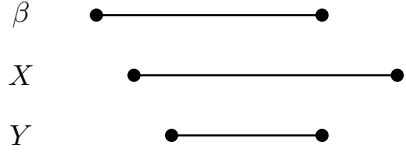
\begin{tikzpicture}[thick]
\draw (-1,0) -- (2, 0);
\node (A) at (-2,0) [] {$\beta$};
\node (R) at (-1,0) [] {};
\node (P) at (2,0) [] {};
\node (E) at (0,-1.6) [] {};
\node (G) at (2, -1.6) [] {};
\fill (R) circle (2.5pt);
\fill (P) circle (2.5pt);
\draw (-0.5,-0.8) -- (3, -0.8);
\node (C) at (-2,-0.8) [circle] {$X$};
\node (D) at (-2,-1.6) [circle] {$Y$};
\node (RE) at (-0.5,-0.8) [] {};
\node (PE) at (3,-0.8) [] {};
\fill (RE) circle (2.5pt);
\fill (PE) circle (2.5pt);
\fill (E) circle (2.5pt);
\fill (G) circle (2.5pt);
\draw (0,-1.6) -- (2, -1.6);
\end{tikzpicture}
$$
which gives
$$m_{X\hspace{0,03cm}\cup\hspace{0,03cm} \beta}=m_{X}+m_{X\hspace{0,05cm} \cup\hspace{0,03cm} \beta-X}-\delta_{X,\alpha}\ell=m_{\beta}+m_{X\hspace{0,03cm}\cup \hspace{0,03cm} \beta-\beta}-\delta_{X,\alpha'}\ell$$
$$ m_{\beta}=m_Y+m_{\beta-Y}-\delta_{X,\alpha'}\ell$$
Now substituting implies ($\gamma_2=(X\hspace{0,03cm}\cup \hspace{0,03cm} \beta)-\beta$,\ $\beta-Y=(X\hspace{0,03cm}\cup \hspace{0,03cm} \beta)-X+\gamma_1$)
$$m_X=m_{\beta}-\delta_{X,\alpha'}\ell-m_{X\hspace{0,03cm}\cup \hspace{0,03cm} \beta-\beta}+m_{\gamma_2}+\delta_{X,\alpha}\ell$$ 
and we get 
$$m_X=m_Y+m_{\beta-Y}-2\delta_{X,\alpha'}-m_{X\hspace{0,03cm}\cup \hspace{0,03cm} \beta-\beta}+m_{\gamma_2}+\delta_{X,\alpha}\ell$$
which is the desired property. \par 
Now suppose the following situation
$$
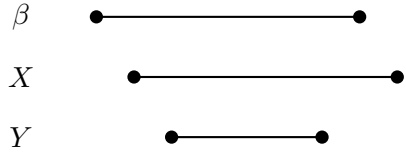
\begin{tikzpicture}[thick]
\draw (-1,0) -- (2.5, 0);
\node (A) at (-2,0) [] {$\beta$};
\node (R) at (-1,0) [] {};
\node (P) at (2.5,0) [] {};
\node (E) at (0,-1.6) [] {};
\node (G) at (2, -1.6) [] {};
\fill (R) circle (2.5pt);
\fill (P) circle (2.5pt);
\draw (-0.5,-0.8) -- (3, -0.8);
\node (C) at (-2,-0.8) [circle] {$X$};
\node (D) at (-2,-1.6) [circle] {$Y$};
\node (RE) at (-0.5,-0.8) [] {};
\node (PE) at (3,-0.8) [] {};
\fill (RE) circle (2.5pt);
\fill (PE) circle (2.5pt);
\fill (E) circle (2.5pt);
\fill (G) circle (2.5pt);
\draw (0,-1.6) -- (2, -1.6);
\end{tikzpicture}
$$
We have 
\begin{align}\label{1112}m_{X\hspace{0,03cm}\cup\hspace{0,03cm} \beta}&=m_{X}+m_{X\hspace{0,05cm} \cup\hspace{0,03cm} \beta-X}-\delta_{X,\alpha}\ell=m_{\beta}+m_{X\hspace{0,03cm}\cup \hspace{0,03cm} \beta-\beta}-\delta_{X,\alpha'}\ell&\\&\notag=m_{\gamma_1}+m_{Y}+m_{\gamma_2}-3\delta_{X,\alpha'}\ell+m_{X\hspace{0,03cm}\cup \hspace{0,03cm} \beta-\beta}.\end{align}
Since $\gamma_1=\tilde{\gamma}_1+(X\hspace{0,03cm}\cup \hspace{0,03cm} \beta-X)$ and $\tilde{\gamma}_2=\gamma_2+(X\hspace{0,03cm}\cup \hspace{0,03cm} \beta-\beta)$ we get (note that $m_{X\hspace{0,05cm} \cup\hspace{0,03cm} \beta-X}=m_{\gamma_1}-m_{\tilde{\gamma}_1}+\delta_{X,\alpha}\ell$; otherwise we would contradict \eqref{1112})
$$m_{X\hspace{0,03cm}\cup\hspace{0,03cm} \beta}=m_{X}+m_{\gamma_1}-m_{\tilde{\gamma}_1}=m_{\gamma_1}+m_{Y}+m_{\gamma_2}-3\delta_{X,\alpha'}\ell+m_{X\hspace{0,03cm}\cup \hspace{0,03cm} \beta-\beta}$$
$$=m_{\gamma_1}+m_{Y}+m_{\gamma_2}+(m_{\tilde{\gamma}_2}-m_{\gamma_2})-2\delta_{X,\alpha'}\ell$$
Hence 
$$m_{X}=m_{Y}+m_{\tilde{\gamma}_2}+m_{\tilde{\gamma}_1}-2\delta_{X,\alpha'}\ell$$ and we are done in this case also. \par 
Now suppose
$$
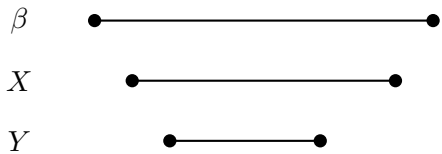
\begin{tikzpicture}[thick]
\draw (-1,0) -- (3.5, 0);
\node (A) at (-2,0) [] {$\beta$};
\node (R) at (-1,0) [] {};
\node (P) at (3.5,0) [] {};
\node (E) at (0,-1.6) [] {};
\node (G) at (2, -1.6) [] {};
\fill (R) circle (2.5pt);
\fill (P) circle (2.5pt);
\draw (-0.5,-0.8) -- (3, -0.8);
\node (C) at (-2,-0.8) [circle] {$X$};
\node (D) at (-2,-1.6) [circle] {$Y$};
\node (RE) at (-0.5,-0.8) [] {};
\node (PE) at (3,-0.8) [] {};
\fill (RE) circle (2.5pt);
\fill (PE) circle (2.5pt);
\fill (E) circle (2.5pt);
\fill (G) circle (2.5pt);
\draw (0,-1.6) -- (2, -1.6);
\end{tikzpicture}
$$

We write 
$$m_{\beta}=m_{\gamma_1}+m_{X}+m_{\gamma_2}-2\delta_{X,\alpha}\ell=m_{\tilde{\gamma}_1}+m_{Y}+m_{\tilde{\gamma}_2}-2\delta_{X,\alpha'}\ell$$
and obtain
$$m_{X}=(m_{\tilde{\gamma}_1}-m_{\gamma_1})+m_{Y}+(m_{\tilde{\gamma}_2}-m_{\gamma_2})+2\delta_{X,\alpha}\ell-2\delta_{X,\alpha'}\ell=m_{\tilde{\gamma}_1-\gamma_2}+m_{Y}+m_{\tilde{\gamma}_2-\gamma_2}-2\delta_{X,\alpha'}\ell$$
and this case is done. \par 
Now suppose 
$$
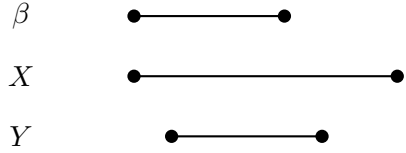
\begin{tikzpicture}[thick]
\draw (-0.5,0) -- (1.5, 0);
\node (A) at (-2,0) [] {$\beta$};
\node (R) at (-0.5,0) [] {};
\node (P) at (1.5,0) [] {};
\node (E) at (0,-1.6) [] {};
\node (G) at (2, -1.6) [] {};
\fill (R) circle (2.5pt);
\fill (P) circle (2.5pt);
\draw (-0.5,-0.8) -- (3, -0.8);
\node (C) at (-2,-0.8) [circle] {$X$};
\node (D) at (-2,-1.6) [circle] {$Y$};
\node (RE) at (-0.5,-0.8) [] {};
\node (PE) at (3,-0.8) [] {};
\fill (RE) circle (2.5pt);
\fill (PE) circle (2.5pt);
\fill (E) circle (2.5pt);
\fill (G) circle (2.5pt);
\draw (0,-1.6) -- (2, -1.6);
\end{tikzpicture}
$$
This gives 
$$m_{Y\hspace{0,03cm}\cup\hspace{0,03cm} \beta}=m_{Y}+m_{Y\hspace{0,05cm} \cup\hspace{0,03cm} \beta-Y}-\delta_{X,\alpha'}\ell=m_{\beta}+m_{Y\hspace{0,03cm}\cup \hspace{0,03cm} \beta-\beta}-\delta_{X,\alpha}\ell$$
and $m_{X}=m_{\beta}+m_{X-\beta}-\delta_{X,\alpha'}\ell$.
Hence 
$$m_{Y\hspace{0,03cm}\cup\hspace{0,03cm} \beta}+m_{\gamma_2}=m_{Y}+m_{Y\hspace{0,05cm} \cup\hspace{0,03cm} \beta-Y}+m_{\gamma_2}-\delta_{X,\alpha'}\ell=m_{\beta}+m_{Y\hspace{0,03cm}\cup \hspace{0,03cm} \beta-\beta}-\delta_{X,\alpha}\ell+m_{\gamma_2}$$
and we are done if $m_{X}+\delta_{X,\alpha'}\ell=m_{Y\hspace{0,03cm}\cup\hspace{0,03cm} \beta}+m_{\gamma_2}$.
Otherwise 
$$m_{X}+\delta_{X,\alpha}\ell=m_{Y}+m_{\gamma_1}+m_{\gamma_2}-\delta_{X,\alpha'}\ell=m_{\beta}+m_{Y\hspace{0,03cm}\cup \hspace{0,03cm} \beta-\beta}-\delta_{X,\alpha}\ell+m_{\gamma_2}$$
and we end in a contradiction
$$m_{X-\beta}=m_{Y\hspace{0,03cm}\cup \hspace{0,03cm} \beta-\beta}+m_{\gamma_2}-2\delta_{X,\alpha}\ell+\delta_{X,\alpha'}\ell.$$
Now suppose that 
$$
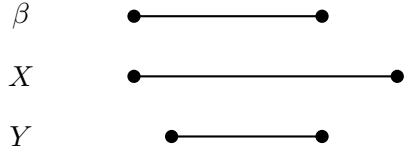
\begin{tikzpicture}[thick]
\draw (-0.5,0) -- (2, 0);
\node (A) at (-2,0) [] {$\beta$};
\node (R) at (-0.5,0) [] {};
\node (P) at (2,0) [] {};
\node (E) at (0,-1.6) [] {};
\node (G) at (2, -1.6) [] {};
\fill (R) circle (2.5pt);
\fill (P) circle (2.5pt);
\draw (-0.5,-0.8) -- (3, -0.8);
\node (C) at (-2,-0.8) [circle] {$X$};
\node (D) at (-2,-1.6) [circle] {$Y$};
\node (RE) at (-0.5,-0.8) [] {};
\node (PE) at (3,-0.8) [] {};
\fill (RE) circle (2.5pt);
\fill (PE) circle (2.5pt);
\fill (E) circle (2.5pt);
\fill (G) circle (2.5pt);
\draw (0,-1.6) -- (2, -1.6);
\end{tikzpicture}
$$
which gives
$$m_{X}=m_{\beta}+m_{X-\beta}-\ell \delta_{\alpha',X}, \ \ m_{\beta}=m_{Y}+m_{\beta-Y}-\ell \delta_{\alpha',X}$$
and thus
$$m_{X}=m_{Y}+m_{\beta-Y}+m_{X-\beta}-2\ell \delta_{\alpha',X}$$
which is the desired property since $\gamma_1=\beta-Y$ and $\gamma_2=X-\beta$. \par

Now suppose that 
$$
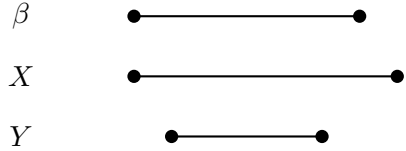
\begin{tikzpicture}[thick]
\draw (-0.5,0) -- (2.5, 0);
\node (A) at (-2,0) [] {$\beta$};
\node (R) at (-0.5,0) [] {};
\node (P) at (2.5,0) [] {};
\node (E) at (0,-1.6) [] {};
\node (G) at (2, -1.6) [] {};
\fill (R) circle (2.5pt);
\fill (P) circle (2.5pt);
\draw (-0.5,-0.8) -- (3, -0.8);
\node (C) at (-2,-0.8) [circle] {$X$};
\node (D) at (-2,-1.6) [circle] {$Y$};
\node (RE) at (-0.5,-0.8) [] {};
\node (PE) at (3,-0.8) [] {};
\fill (RE) circle (2.5pt);
\fill (PE) circle (2.5pt);
\fill (E) circle (2.5pt);
\fill (G) circle (2.5pt);
\draw (0,-1.6) -- (2, -1.6);
\end{tikzpicture}
$$
In this case we have $m_{\beta}=m_{\gamma_1}+m_{\gamma_2}+m_Y-2\ell \delta_{X,\alpha'}$ and $m_{X}=m_{\beta}+m_{X-\beta}-\delta_{X,\alpha'}\ell$. Now noting that $\tilde{\gamma}_2=\gamma_2+(X-\beta)$ we obtain by substitution $$m_{X}=m_{\beta}+m_{X-\beta}-\delta_{X,\alpha'}\ell=m_{\gamma_1}+m_{\gamma_2}+m_Y+m_{X-\beta}-3\ell \delta_{X,\alpha'}=m_{\gamma_1}+m_{\tilde{\gamma}_2}+m_Y-2\ell \delta_{X,\alpha'}$$
and we are done.\par 

Now suppose that 
$$
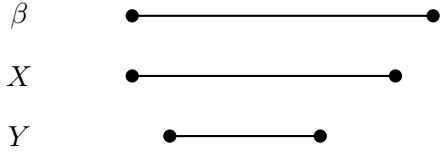
\begin{tikzpicture}[thick]
\draw (-0.5,0) -- (3.5, 0);
\node (A) at (-2,0) [] {$\beta$};
\node (R) at (-0.5,0) [] {};
\node (P) at (3.5,0) [] {};
\node (E) at (0,-1.6) [] {};
\node (G) at (2, -1.6) [] {};
\fill (R) circle (2.5pt);
\fill (P) circle (2.5pt);
\draw (-0.5,-0.8) -- (3, -0.8);
\node (C) at (-2,-0.8) [circle] {$X$};
\node (D) at (-2,-1.6) [circle] {$Y$};
\node (RE) at (-0.5,-0.8) [] {};
\node (PE) at (3,-0.8) [] {};
\fill (RE) circle (2.5pt);
\fill (PE) circle (2.5pt);
\fill (E) circle (2.5pt);
\fill (G) circle (2.5pt);
\draw (0,-1.6) -- (2, -1.6);
\end{tikzpicture}
$$
Then we have 
$$m_{\beta}=m_{\gamma_1}+m_Y+m_{\gamma_2}-2\ell\delta_{X,\alpha'},\ \ \ m_{\beta}=m_X+m_{\beta-X}-\delta_{X,\alpha}\ell.$$
Thus
$$m_X=m_{\gamma_1}+m_Y+m_{\gamma_2}-m_{\beta-X}-2\ell\delta_{X,\alpha'}+\delta_{X,\alpha}\ell$$
which is the desired property since $m_{\gamma_2}-m_{\beta-X}=m_{\tilde{\gamma}_2}-\delta_{X,\alpha}\ell$.\par

Now suppose that
$$
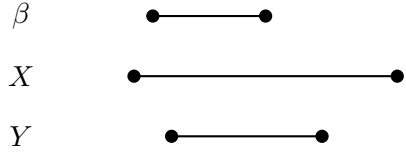
\begin{tikzpicture}[thick]
\draw (-0.25,0) -- (1.25, 0);
\node (A) at (-2,0) [] {$\beta$};
\node (R) at (-0.25,0) [] {};
\node (P) at (1.25,0) [] {};
\node (E) at (0,-1.6) [] {};
\node (G) at (2, -1.6) [] {};
\fill (R) circle (2.5pt);
\fill (P) circle (2.5pt);
\draw (-0.5,-0.8) -- (3, -0.8);
\node (C) at (-2,-0.8) [circle] {$X$};
\node (D) at (-2,-1.6) [circle] {$Y$};
\node (RE) at (-0.5,-0.8) [] {};
\node (PE) at (3,-0.8) [] {};
\fill (RE) circle (2.5pt);
\fill (PE) circle (2.5pt);
\fill (E) circle (2.5pt);
\fill (G) circle (2.5pt);
\draw (0,-1.6) -- (2, -1.6);
\end{tikzpicture}
$$
Then we have 
$$m_X=m_{\gamma_1}+m_{\beta}+m_{\gamma_2}-2\ell\delta_{X,\alpha'}$$
and $$m_{Y\hspace{0,03cm}\cup\hspace{0,03cm} \beta}=m_{Y}+m_{Y\hspace{0,05cm} \cup\hspace{0,03cm} \beta-Y}-\delta_{X,\alpha'}\ell=m_{\beta}+m_{Y\hspace{0,03cm}\cup \hspace{0,03cm} \beta-\beta}-\delta_{X,\alpha}\ell$$
Moreover, we have 
$$\tilde{\gamma}_1=\gamma_1+(Y\hspace{0,03cm}\cup\hspace{0,03cm} \beta-Y),\ \ \tilde{\gamma}_2=\gamma_2-(Y\hspace{0,03cm}\cup\hspace{0,03cm} \beta-\beta)$$
and together with
$$m_X=m_{\gamma_1}+m_{\gamma_2}-3\ell\delta_{X,\alpha'}+m_Y+m_{Y\hspace{0,03cm}\cup\hspace{0,03cm} \beta-Y}-m_{Y\hspace{0,03cm}\cup\hspace{0,03cm} \beta-\beta}+\delta_{X,\alpha}\ell$$
we get the desired property
$$m_X=m_{\tilde{\gamma}_1}+m_{Y}+m_{\tilde{\gamma}_2}-2\ell\delta_{X,\alpha'}.$$
Now we assume that
$$
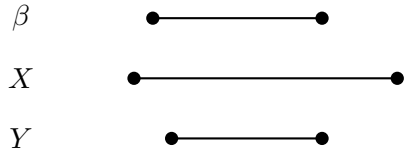
\begin{tikzpicture}[thick]
\draw (-0.25,0) -- (2, 0);
\node (A) at (-2,0) [] {$\beta$};
\node (R) at (-0.25,0) [] {};
\node (P) at (2,0) [] {};
\node (E) at (0,-1.6) [] {};
\node (G) at (2, -1.6) [] {};
\fill (R) circle (2.5pt);
\fill (P) circle (2.5pt);
\draw (-0.5,-0.8) -- (3, -0.8);
\node (C) at (-2,-0.8) [circle] {$X$};
\node (D) at (-2,-1.6) [circle] {$Y$};
\node (RE) at (-0.5,-0.8) [] {};
\node (PE) at (3,-0.8) [] {};
\fill (RE) circle (2.5pt);
\fill (PE) circle (2.5pt);
\fill (E) circle (2.5pt);
\fill (G) circle (2.5pt);
\draw (0,-1.6) -- (2, -1.6);
\end{tikzpicture}
$$
Then we have 
$$m_X=m_{\gamma_1}+m_{\beta}+m_{\gamma_2}-2\ell\delta_{X,\alpha'}$$
and $m_{\beta}=m_Y+m_{\beta-Y}-\delta_{X,\alpha'}\ell$. Note that $(\beta-Y)+\gamma_1=\tilde{\gamma}_1$ and $\tilde{\gamma}_2=\gamma_2$. We get once more the claimed equation
$$m_X=m_{\gamma_1}+m_{\beta}+m_{\gamma_2}-2\ell\delta_{X,\alpha'}=m_{\gamma_1}+m_Y+m_{\beta-Y}+m_{\tilde{\gamma}_2}-3\ell\delta_{X,\alpha'}=m_{\tilde{\gamma}_1}+m_Y+m_{\tilde{\gamma}_2}-2\ell\delta_{X,\alpha'}.$$
Now suppose that 
$$
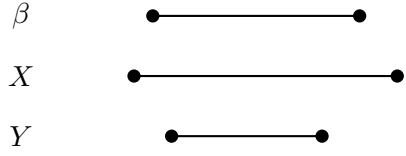
\begin{tikzpicture}[thick]
\draw (-0.25,0) -- (2.5, 0);
\node (A) at (-2,0) [] {$\beta$};
\node (R) at (-0.25,0) [] {};
\node (P) at (2.5,0) [] {};
\node (E) at (0,-1.6) [] {};
\node (G) at (2, -1.6) [] {};
\fill (R) circle (2.5pt);
\fill (P) circle (2.5pt);
\draw (-0.5,-0.8) -- (3, -0.8);
\node (C) at (-2,-0.8) [circle] {$X$};
\node (D) at (-2,-1.6) [circle] {$Y$};
\node (RE) at (-0.5,-0.8) [] {};
\node (PE) at (3,-0.8) [] {};
\fill (RE) circle (2.5pt);
\fill (PE) circle (2.5pt);
\fill (E) circle (2.5pt);
\fill (G) circle (2.5pt);
\draw (0,-1.6) -- (2, -1.6);
\end{tikzpicture}
$$
Then we have 
$$m_X=m_{\gamma_1}+m_{\beta}+m_{\gamma_2}-2\ell\delta_{X,\alpha'}$$
and 
$$m_{\beta}=m_Y+m_{\tilde{\gamma}_1}+m_{\tilde{\gamma}_2}-2\ell \delta_{X,\alpha'}$$  
and again we get the desired property by substituting the second equation into the first one. \par
Now suppose that
$$
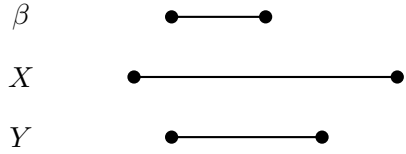
\begin{tikzpicture}[thick]
\draw (0,0) -- (1.25, 0);
\node (A) at (-2,0) [] {$\beta$};
\node (R) at (0,0) [] {};
\node (P) at (1.25,0) [] {};
\node (E) at (0,-1.6) [] {};
\node (G) at (2, -1.6) [] {};
\fill (R) circle (2.5pt);
\fill (P) circle (2.5pt);
\draw (-0.5,-0.8) -- (3, -0.8);
\node (C) at (-2,-0.8) [circle] {$X$};
\node (D) at (-2,-1.6) [circle] {$Y$};
\node (RE) at (-0.5,-0.8) [] {};
\node (PE) at (3,-0.8) [] {};
\fill (RE) circle (2.5pt);
\fill (PE) circle (2.5pt);
\fill (E) circle (2.5pt);
\fill (G) circle (2.5pt);
\draw (0,-1.6) -- (2, -1.6);
\end{tikzpicture}
$$
Then we have 
$$m_X=m_{\gamma_1}+m_{\beta}+m_{\gamma_2}-2\ell\delta_{X,\alpha'}$$
and $m_Y=m_{\beta}+m_{Y-\beta}-\delta_{X,\alpha}\ell$. We obtain
$$m_X=m_{\gamma_1}+m_{\beta}+m_{\gamma_2}-2\ell\delta_{X,\alpha'}=m_{\tilde{\gamma}_1}+m_{Y}-m_{Y-\beta}+\delta_{X,\alpha}\ell+m_{\gamma_2}-2\ell\delta_{X,\alpha'}$$
and the statement is clear with $m_{\tilde{\gamma}_2}=m_{\gamma_2}-m_{Y-\beta}-\delta_{X,\alpha}\ell$. 
Finally we suppose that 
$$
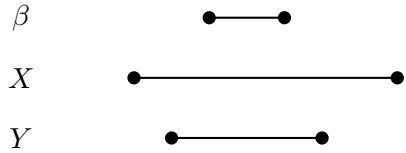
\begin{tikzpicture}[thick]
\draw (0.5,0) -- (1.5, 0);
\node (A) at (-2,0) [] {$\beta$};
\node (R) at (0.5,0) [] {};
\node (P) at (1.5,0) [] {};
\node (E) at (0,-1.6) [] {};
\node (G) at (2, -1.6) [] {};
\fill (R) circle (2.5pt);
\fill (P) circle (2.5pt);
\draw (-0.5,-0.8) -- (3, -0.8);
\node (C) at (-2,-0.8) [circle] {$X$};
\node (D) at (-2,-1.6) [circle] {$Y$};
\node (RE) at (-0.5,-0.8) [] {};
\node (PE) at (3,-0.8) [] {};
\fill (RE) circle (2.5pt);
\fill (PE) circle (2.5pt);
\fill (E) circle (2.5pt);
\fill (G) circle (2.5pt);
\draw (0,-1.6) -- (2, -1.6);
\end{tikzpicture}
$$
Then we have 
$$m_X=m_{\gamma_1}+m_{\beta}+m_{\gamma_2}-2\ell\delta_{X,\alpha'}$$
$$m_Y=m_{\tilde{\gamma}_1}+m_{\beta}+m_{\tilde{\gamma}_2}-2\ell\delta_{Y,\alpha'}$$
Hence 
$$m_X=m_{\gamma_1}+m_Y-m_{\tilde{\gamma}_1}-m_{\tilde{\gamma}_2}+m_{\gamma_2}-2\ell\delta_{X,\alpha'}+2\ell\delta_{Y,\alpha'}$$
$$=m_{\gamma_1-\tilde{\gamma}_1}+m_Y+m_{\gamma_2-\tilde{\gamma}_2}-2\ell \delta_{X,\alpha'}$$
which finishes the proof in this case.\par
Note that all other cases are just the reflected versions of the above cases. 

\textit{Case 3}: Here we suppose that $\alpha-\alpha'\notin R$
and
$\mathrm{supp}(\alpha)\nsubseteq \mathrm{supp}(\alpha')\nsubseteq \mathrm{supp}(\alpha)$. So we have to show in the rest of the calculations that 
$$m_{\alpha\hspace{0,03cm}\cup\hspace{0,03cm} \alpha'}=m_{\alpha}+m_{\alpha\hspace{0,05cm} \cup\hspace{0,03cm} \alpha'-\alpha}-\ell=m_{\alpha'}+m_{\alpha\hspace{0,03cm}\cup \hspace{0,03cm} \alpha'-\alpha'}$$
We assume by contradiction one of the following cases
\begin{enumerate}
    \item $m_{\alpha}+m_{\alpha\hspace{0,05cm} \cup\hspace{0,03cm} \alpha'-\alpha}=m_{\alpha'}+m_{\alpha\hspace{0,03cm}\cup \hspace{0,03cm} \alpha'-\alpha'}$
    
    \vspace{0.3cm}
    \item $m_{\alpha\hspace{0,03cm}\cup\hspace{0,03cm} \alpha'}=m_{\alpha}+m_{\alpha\hspace{0,05cm} \cup\hspace{0,03cm} \alpha'-\alpha}=m_{\alpha'}+m_{\alpha\hspace{0,03cm}\cup \hspace{0,03cm} \alpha'-\alpha'}-\ell$
\end{enumerate}
In the second case we have $\alpha'\succeq \alpha\succeq \beta$ and from the calculations above we can assume that $\beta\succeq \alpha'$ is induced from Definition~\ref{poset1}(iii); otherwise we know the transitivity already which would give $\alpha'\succeq \beta$. Hence 
$$m_{\alpha'}+m_{\alpha'\hspace{0,05cm} \cup\hspace{0,03cm} \beta-\alpha'}=m_{\beta}+m_{\alpha'\hspace{0,03cm}\cup \hspace{0,03cm} \beta-\beta}-\ell.$$
So we get 
$$m_{\alpha}+m_{\alpha\hspace{0,05cm} \cup\hspace{0,03cm} \alpha'-\alpha}=m_{\beta}+m_{\alpha'\hspace{0,03cm}\cup \hspace{0,03cm} \beta-\beta}-m_{\alpha'\hspace{0,05cm} \cup\hspace{0,03cm} \beta-\alpha'}+m_{\alpha\hspace{0,03cm}\cup \hspace{0,03cm} \alpha'-\alpha'}-2\ell$$
which contradicts $m_{\alpha}\geq m_{\beta}$ (see Remark~\ref{rem1}). Thus we can assume in the rest of this case that (1) holds:
\begin{equation}\label{ann12}m_{\alpha}+m_{\alpha\hspace{0,05cm} \cup\hspace{0,03cm} \alpha'-\alpha}=m_{\alpha'}+m_{\alpha\hspace{0,03cm}\cup \hspace{0,03cm} \alpha'-\alpha'}\end{equation}
\textit{Case 3.1}: In this case we assume that $\beta\succeq \alpha'$ is induced from Definition~\ref{poset1}(iii). Similarly as above we get 
\begin{equation}\label{raw1}m_{\alpha}+m_{\alpha\hspace{0,05cm} \cup\hspace{0,03cm} \alpha'-\alpha}=m_{\beta}+m_{\alpha'\hspace{0,03cm}\cup \hspace{0,03cm} \beta-\beta}-m_{\alpha'\hspace{0,05cm} \cup\hspace{0,03cm} \beta-\alpha'}+m_{\alpha\hspace{0,03cm}\cup \hspace{0,03cm} \alpha'-\alpha'}-\ell\end{equation}
\textit{Case 3.1.1}: If $\alpha\succeq \beta$ is also induced from Definition~\ref{poset1}(iii) we get further
\begin{equation}\label{raw}m_{\alpha\hspace{0,05cm} \cup\hspace{0,03cm} \beta-\beta}+m_{\alpha\hspace{0,05cm} \cup\hspace{0,03cm} \alpha'-\alpha}=m_{\alpha'\hspace{0,03cm}\cup \hspace{0,03cm} \beta-\beta}-m_{\alpha'\hspace{0,05cm} \cup\hspace{0,03cm} \beta-\alpha'}+m_{\alpha\hspace{0,03cm}\cup \hspace{0,03cm} \alpha'-\alpha'}+ m_{\alpha\hspace{0,05cm} \cup\hspace{0,03cm} \beta-\alpha}-2\ell\end{equation}
This will end in a contradiction by considering the possible constellations of roots. We consider the following two possible cases 
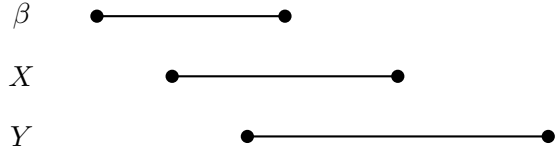
\begin{enumerate}
\item $$\begin{tikzpicture}[thick]
\draw (-1,0) -- (1.5, 0);
\node (A) at (-2,0) [] {$\beta$};
\node (R) at (-1,0) [] {};
\node (P) at (1.5,0) [] {};
\node (E) at (1,-1.6) [] {};
\node (G) at (5, -1.6) [] {};
\fill (R) circle (2.5pt);
\fill (P) circle (2.5pt);
\draw (0,-0.8) -- (3, -0.8);
\node (C) at (-2,-0.8) [circle] {$X$};
\node (D) at (-2,-1.6) [circle] {$Y$};
\node (RE) at (0,-0.8) [] {};
\node (PE) at (3,-0.8) [] {};
\fill (RE) circle (2.5pt);
\fill (PE) circle (2.5pt);
\fill (E) circle (2.5pt);
\fill (G) circle (2.5pt);
\draw (1,-1.6) -- (5, -1.6);
\end{tikzpicture}
$$
\item 
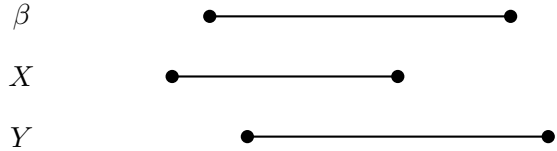
$$\begin{tikzpicture}[thick]
\draw (0.5,0) -- (4.5, 0);
\node (A) at (-2,0) [] {$\beta$};
\node (R) at (0.5,0) [] {};
\node (P) at (4.5,0) [] {};
\node (E) at (1,-1.6) [] {};
\node (G) at (5, -1.6) [] {};
\fill (R) circle (2.5pt);
\fill (P) circle (2.5pt);
\draw (0,-0.8) -- (3, -0.8);
\node (C) at (-2,-0.8) [circle] {$X$};
\node (D) at (-2,-1.6) [circle] {$Y$};
\node (RE) at (0,-0.8) [] {};
\node (PE) at (3,-0.8) [] {};
\fill (RE) circle (2.5pt);
\fill (PE) circle (2.5pt);
\fill (E) circle (2.5pt);
\fill (G) circle (2.5pt);
\draw (1,-1.6) -- (5, -1.6);
\end{tikzpicture}
$$
\end{enumerate}
where $\{X,Y\}=\{\alpha,\alpha'\}$. In the first case we see that $(X\cup Y-Y)+(\beta\cup X-X)=(\beta\cup Y-Y)$ and in the second case $(X\cup Y-Y)-(\beta\cup X-\beta)=(\beta\cup Y-Y)$ which we substitute into \eqref{raw}. This will lead to a contradiction which we demonstrate only in the first case when $X=\alpha$. So we have in this case
$$m_{\alpha\hspace{0,03cm}\cup \hspace{0,03cm} \alpha'-\alpha'}+ m_{\alpha\hspace{0,05cm} \cup\hspace{0,03cm} \beta-\alpha}=m_{\alpha'\hspace{0,05cm} \cup\hspace{0,03cm} \beta-\alpha'}$$
or 
$$m_{\alpha\hspace{0,03cm}\cup \hspace{0,03cm} \alpha'-\alpha'}+ m_{\alpha\hspace{0,05cm} \cup\hspace{0,03cm} \beta-\alpha}=m_{\alpha'\hspace{0,05cm} \cup\hspace{0,03cm} \beta-\alpha'}+\ell$$
but both equations obviously give a contradiction after substituting into \eqref{raw}.

\textit{Case 3.1.2}: If $\alpha\succeq \beta$ is induced from Definition~\ref{poset1}(i), then a straightforward inspection of the roots shows that we must have $$\alpha'\hspace{0,03cm}\cup \hspace{0,03cm} 
\beta-\alpha'=\alpha\hspace{0,03cm}\cup \hspace{0,03cm} \alpha'-\alpha$$ or 
$$\alpha'\hspace{0,03cm}\cup \hspace{0,03cm} 
\beta-\beta=\alpha\hspace{0,03cm}\cup \hspace{0,03cm} \alpha'-\alpha.$$
which contradicts \eqref{raw1}

\textit{Case 3.1.3}: If $\alpha\succeq \beta$ is induced from Definition~\ref{poset1}(ii), then we must have one of the following constellations:

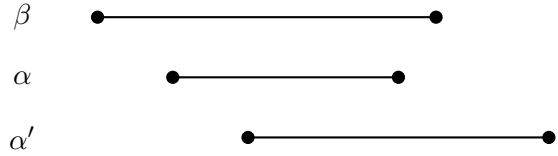
\begin{enumerate}
\item $$\begin{tikzpicture}[thick]
\draw (-1,0) -- (3.5, 0);
\node (A) at (-2,0) [] {$\beta$};
\node (R) at (-1,0) [] {};
\node (P) at (3.5,0) [] {};
\node (E) at (1,-1.6) [] {};
\node (G) at (5, -1.6) [] {};
\fill (R) circle (2.5pt);
\fill (P) circle (2.5pt);
\draw (0,-0.8) -- (3, -0.8);
\node (C) at (-2,-0.8) [circle] {$\alpha$};
\node (D) at (-2,-1.6) [circle] {$\alpha'$};
\node (RE) at (0,-0.8) [] {};
\node (PE) at (3,-0.8) [] {};
\fill (RE) circle (2.5pt);
\fill (PE) circle (2.5pt);
\fill (E) circle (2.5pt);
\fill (G) circle (2.5pt);
\draw (1,-1.6) -- (5, -1.6);
\end{tikzpicture}
$$
\item 
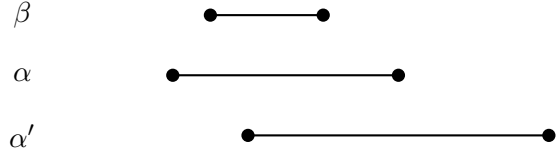
$$\begin{tikzpicture}[thick]
\draw (0.5,0) -- (2, 0);
\node (A) at (-2,0) [] {$\beta$};
\node (R) at (0.5,0) [] {};
\node (P) at (2,0) [] {};
\node (E) at (1,-1.6) [] {};
\node (G) at (5, -1.6) [] {};
\fill (R) circle (2.5pt);
\fill (P) circle (2.5pt);
\draw (0,-0.8) -- (3, -0.8);
\node (C) at (-2,-0.8) [circle] {$\alpha$};
\node (D) at (-2,-1.6) [circle] {$\alpha'$};
\node (RE) at (0,-0.8) [] {};
\node (PE) at (3,-0.8) [] {};
\fill (RE) circle (2.5pt);
\fill (PE) circle (2.5pt);
\fill (E) circle (2.5pt);
\fill (G) circle (2.5pt);
\draw (1,-1.6) -- (5, -1.6);
\end{tikzpicture}
$$
\end{enumerate}
Now \eqref{raw1} turns into

$$m_{\alpha\hspace{0,05cm} \cup\hspace{0,03cm} \alpha'-\alpha}=m_{\gamma_1}+m_{\gamma_2}+m_{\alpha'\hspace{0,03cm}\cup \hspace{0,03cm} \beta-\beta}-m_{\alpha'\hspace{0,05cm} \cup\hspace{0,03cm} \beta-\alpha'}+m_{\alpha\hspace{0,03cm}\cup \hspace{0,03cm} \alpha'-\alpha'}-3\ell$$
in the first case and into
$$m_{\gamma_1}+m_{\gamma_2}+m_{\alpha\hspace{0,05cm} \cup\hspace{0,03cm} \alpha'-\alpha}=m_{\alpha'\hspace{0,03cm}\cup \hspace{0,03cm} \beta-\beta}-m_{\alpha'\hspace{0,05cm} \cup\hspace{0,03cm} \beta-\alpha'}+m_{\alpha\hspace{0,03cm}\cup \hspace{0,03cm} \alpha'-\alpha'}-\ell$$
in the second case.
In the first case we have $$\alpha'\hspace{0,03cm}\cup \hspace{0,03cm} \beta-\alpha'=\gamma_1+(\alpha\hspace{0,03cm}\cup \hspace{0,03cm} \alpha'-\alpha')$$
leading to a contradiction and in the second case we have 
$$(\alpha\hspace{0,03cm}\cup \hspace{0,03cm} \alpha'-\alpha')=\gamma_1+(\alpha'\hspace{0,03cm}\cup \hspace{0,03cm} \beta-\alpha')$$
leading once more to a contradiction.

\textit{Case 3.2}: In this case we assume that $\beta\succeq \alpha'$ is induced from Definition~\ref{poset1}(ii). We can have the following six cases
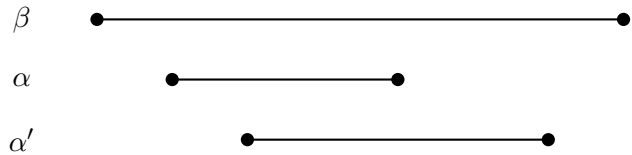
\begin{enumerate}
\item $$\begin{tikzpicture}[thick]
\draw (-1,0) -- (6, 0);
\node (A) at (-2,0) [] {$\beta$};
\node (R) at (-1,0) [] {};
\node (P) at (6,0) [] {};
\node (E) at (1,-1.6) [] {};
\node (G) at (5, -1.6) [] {};
\fill (R) circle (2.5pt);
\fill (P) circle (2.5pt);
\draw (0,-0.8) -- (3, -0.8);
\node (C) at (-2,-0.8) [circle] {$\alpha$};
\node (D) at (-2,-1.6) [circle] {$\alpha'$};
\node (RE) at (0,-0.8) [] {};
\node (PE) at (3,-0.8) [] {};
\fill (RE) circle (2.5pt);
\fill (PE) circle (2.5pt);
\fill (E) circle (2.5pt);
\fill (G) circle (2.5pt);
\draw (1,-1.6) -- (5, -1.6);
\end{tikzpicture}
$$
\item 
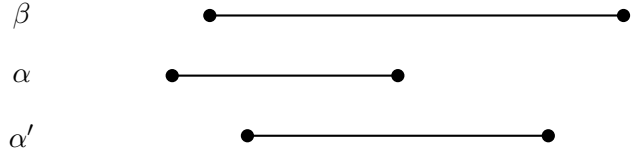
$$\begin{tikzpicture}[thick]
\draw (0.5,0) -- (6, 0);
\node (A) at (-2,0) [] {$\beta$};
\node (R) at (0.5,0) [] {};
\node (P) at (6,0) [] {};
\node (E) at (1,-1.6) [] {};
\node (G) at (5, -1.6) [] {};
\fill (R) circle (2.5pt);
\fill (P) circle (2.5pt);
\draw (0,-0.8) -- (3, -0.8);
\node (C) at (-2,-0.8) [circle] {$\alpha$};
\node (D) at (-2,-1.6) [circle] {$\alpha'$};
\node (RE) at (0,-0.8) [] {};
\node (PE) at (3,-0.8) [] {};
\fill (RE) circle (2.5pt);
\fill (PE) circle (2.5pt);
\fill (E) circle (2.5pt);
\fill (G) circle (2.5pt);
\draw (1,-1.6) -- (5, -1.6);
\end{tikzpicture}
$$
\item 
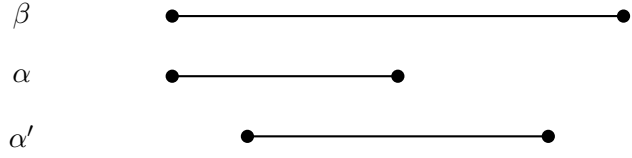
$$\begin{tikzpicture}[thick]
\draw (0,0) -- (6, 0);
\node (A) at (-2,0) [] {$\beta$};
\node (R) at (0,0) [] {};
\node (P) at (6,0) [] {};
\node (E) at (1,-1.6) [] {};
\node (G) at (5, -1.6) [] {};
\fill (R) circle (2.5pt);
\fill (P) circle (2.5pt);
\draw (0,-0.8) -- (3, -0.8);
\node (C) at (-2,-0.8) [circle] {$\alpha$};
\node (D) at (-2,-1.6) [circle] {$\alpha'$};
\node (RE) at (0,-0.8) [] {};
\node (PE) at (3,-0.8) [] {};
\fill (RE) circle (2.5pt);
\fill (PE) circle (2.5pt);
\fill (E) circle (2.5pt);
\fill (G) circle (2.5pt);
\draw (1,-1.6) -- (5, -1.6);
\end{tikzpicture}
$$
\item 
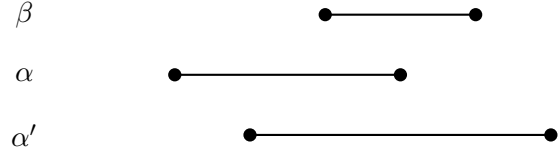
$$\begin{tikzpicture}[thick]
\draw (2,0) -- (4, 0);
\node (A) at (-2,0) [] {$\beta$};
\node (R) at (2,0) [] {};
\node (P) at (4,0) [] {};
\node (E) at (1,-1.6) [] {};
\node (G) at (5, -1.6) [] {};
\fill (R) circle (2.5pt);
\fill (P) circle (2.5pt);
\draw (0,-0.8) -- (3, -0.8);
\node (C) at (-2,-0.8) [circle] {$\alpha$};
\node (D) at (-2,-1.6) [circle] {$\alpha'$};
\node (RE) at (0,-0.8) [] {};
\node (PE) at (3,-0.8) [] {};
\fill (RE) circle (2.5pt);
\fill (PE) circle (2.5pt);
\fill (E) circle (2.5pt);
\fill (G) circle (2.5pt);
\draw (1,-1.6) -- (5, -1.6);
\end{tikzpicture}
$$
\item 
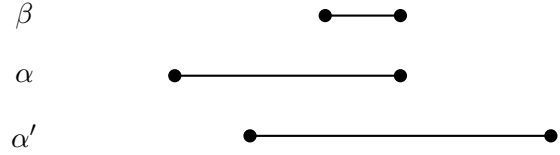
$$\begin{tikzpicture}[thick]
\draw (2,0) -- (3, 0);
\node (A) at (-2,0) [] {$\beta$};
\node (R) at (3,0) [] {};
\node (P) at (2,0) [] {};
\node (E) at (1,-1.6) [] {};
\node (G) at (5, -1.6) [] {};
\fill (R) circle (2.5pt);
\fill (P) circle (2.5pt);
\draw (0,-0.8) -- (3, -0.8);
\node (C) at (-2,-0.8) [circle] {$\alpha$};
\node (D) at (-2,-1.6) [circle] {$\alpha'$};
\node (RE) at (0,-0.8) [] {};
\node (PE) at (3,-0.8) [] {};
\fill (RE) circle (2.5pt);
\fill (PE) circle (2.5pt);
\fill (E) circle (2.5pt);
\fill (G) circle (2.5pt);
\draw (1,-1.6) -- (5, -1.6);
\end{tikzpicture}
$$
\item 
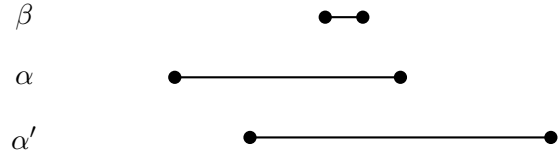
$$\begin{tikzpicture}[thick]
\draw (2,0) -- (2.5, 0);
\node (A) at (-2,0) [] {$\beta$};
\node (R) at (2.5,0) [] {};
\node (P) at (2,0) [] {};
\node (E) at (1,-1.6) [] {};
\node (G) at (5, -1.6) [] {};
\fill (R) circle (2.5pt);
\fill (P) circle (2.5pt);
\draw (0,-0.8) -- (3, -0.8);
\node (C) at (-2,-0.8) [circle] {$\alpha$};
\node (D) at (-2,-1.6) [circle] {$\alpha'$};
\node (RE) at (0,-0.8) [] {};
\node (PE) at (3,-0.8) [] {};
\fill (RE) circle (2.5pt);
\fill (PE) circle (2.5pt);
\fill (E) circle (2.5pt);
\fill (G) circle (2.5pt);
\draw (1,-1.6) -- (5, -1.6);
\end{tikzpicture}
$$
\end{enumerate}
Recall that we have equation \eqref{ann12}. In the first case 
$$m_{\beta}=m_{\alpha'}+m_{\gamma_1}+m_{\gamma_2}=m_{\alpha}+m_{\tilde{\gamma}_1}+m_{\tilde{\gamma}_2}-2\ell.$$ 
Substitution into \eqref{ann12} gives 
$$m_{\gamma_1}+m_{\gamma_2}-m_{\tilde{\gamma}_2}+2\ell+m_{\alpha\hspace{0,05cm} \cup\hspace{0,03cm} \alpha'-\alpha}=m_{\alpha\hspace{0,03cm}\cup \hspace{0,03cm} \alpha'-\alpha'}+m_{\tilde{\gamma}_1}$$
which is a contradiction since $\gamma_1=\tilde{\gamma}_1+(\alpha\hspace{0,05cm} \cup\hspace{0,03cm} \alpha'-\alpha')$.\par
In the second case
$$m_{\beta}=m_{\alpha'}+m_{\gamma_1}+m_{\gamma_2}=m_{\alpha}+m_{\alpha\hspace{0,05cm} \cup\hspace{0,03cm} \beta-\alpha}-m_{\alpha\hspace{0,05cm} \cup\hspace{0,03cm} \beta-\beta}-\ell.$$ 
Substitution into \eqref{ann12} gives 
$$m_{\gamma_1}+m_{\gamma_2}+m_{\alpha\hspace{0,05cm} \cup\hspace{0,03cm} \beta-\beta}+\ell+m_{\alpha\hspace{0,05cm} \cup\hspace{0,03cm} \alpha'-\alpha}=m_{\alpha\hspace{0,03cm}\cup \hspace{0,03cm} \alpha'-\alpha'}+m_{\alpha\hspace{0,05cm} \cup\hspace{0,03cm} \beta-\alpha}$$
which is a contradiction since $\gamma_1=(\alpha\hspace{0,05cm} \cup\hspace{0,03cm} \alpha'-\alpha')+(\alpha\hspace{0,05cm} \cup\hspace{0,03cm} \beta-\alpha)$.\par
In the third case we have 
$$m_{\beta}=m_{\alpha'}+m_{\gamma_1}+m_{\gamma_2}=m_{\alpha}+m_{\beta-\alpha}-\ell.$$
Substitution into \eqref{ann12} gives 
$$m_{\gamma_1}+m_{\gamma_2}+\ell+m_{\alpha\hspace{0,05cm} \cup\hspace{0,03cm} \alpha'-\alpha}=m_{\alpha\hspace{0,03cm}\cup \hspace{0,03cm} \alpha'-\alpha'}+m_{\beta-\alpha}$$
which is a contradiction since $\gamma_1=(\alpha\hspace{0,03cm}\cup \hspace{0,03cm} \alpha'-\alpha')$.\par In the fourth case we have 
$$m_{\beta}=m_{\alpha'}-m_{\gamma_1}-m_{\gamma_2}+2\ell=m_{\alpha}+m_{\alpha\hspace{0,05cm} \cup\hspace{0,03cm} \beta-\alpha}-m_{\alpha\hspace{0,05cm} \cup\hspace{0,03cm} \beta-\beta}-\ell.$$
Substitution into \eqref{ann12} gives
$$-m_{\gamma_2}+3\ell-m_{\alpha\hspace{0,05cm} \cup\hspace{0,03cm} \beta-\alpha}+m_{\alpha\hspace{0,05cm} \cup\hspace{0,03cm} \beta-\beta}
+m_{\alpha\hspace{0,05cm} \cup\hspace{0,03cm} \alpha'-\alpha}=m_{\alpha\hspace{0,03cm}\cup \hspace{0,03cm} \alpha'-\alpha'}+m_{\gamma_1}$$
which is a contradiction since $\gamma_1+(\alpha\hspace{0,03cm}\cup \hspace{0,03cm} \alpha'-\alpha')=\alpha\hspace{0,05cm} \cup\hspace{0,03cm} \beta-\beta$.\par 
In the fifth case we have 
$$m_{\beta}=m_{\alpha'}-m_{\gamma_1}-m_{\gamma_2}+2\ell=m_{\alpha}-m_{\alpha-\beta}.$$
Substitution into \eqref{ann12} gives
$$-m_{\gamma_2}+2\ell+m_{\alpha-\beta}+m_{\alpha\hspace{0,05cm} \cup\hspace{0,03cm} \alpha'-\alpha}=m_{\alpha\hspace{0,03cm}\cup \hspace{0,03cm} \alpha'-\alpha'}+m_{\gamma_1}$$
which is a contradiction since $\gamma_1+(\alpha\hspace{0,03cm}\cup \hspace{0,03cm} \alpha'-\alpha')=\alpha-\beta$.\par 
In the sixth case we have 
$$m_{\beta}=m_{\alpha'}-m_{\gamma_1}-m_{\gamma_2}+2\ell=m_{\alpha}-m_{\tilde{\gamma}_1}-m_{\tilde{\gamma_2}}.$$
Substitution into \eqref{ann12} gives 
$$-m_{\gamma_2}+2\ell+m_{\tilde{\gamma}_1}+m_{\tilde{\gamma_2}}
+m_{\alpha\hspace{0,05cm} \cup\hspace{0,03cm} \alpha'-\alpha}=m_{\alpha\hspace{0,03cm}\cup \hspace{0,03cm} \alpha'-\alpha'}+m_{\gamma_1}$$
which is a contradiction since $\tilde{\gamma_1}=\gamma_1+(\alpha\hspace{0,03cm}\cup \hspace{0,03cm} \alpha'-\alpha')=\alpha-\beta$.

\textit{Case 3.3}: In this case we assume that $\beta\succeq \alpha'$ is induced from Definition~\ref{poset1}(i). Remember once more that we have equation \eqref{ann12}.\par
If $\alpha\succeq \beta$ is induced from Definition~\ref{poset1}(iii), then a simple case consideration shows that either of the following equations hold
\begin{itemize}
    \item $\alpha\cup \alpha'-\alpha=\beta\cup \alpha-\alpha$
    \item $\alpha\cup \alpha'-\alpha'=\beta\cup \alpha-\beta$
\end{itemize}
We get
$$m_{\beta}+m_{\alpha\hspace{0,05cm} \cup\hspace{0,03cm} \beta-\beta}+\ell+m_{\alpha\hspace{0,05cm} \cup\hspace{0,03cm} \alpha'-\alpha}=m_{\alpha'}+m_{\alpha\hspace{0,03cm}\cup \hspace{0,03cm} \alpha'-\alpha'}+m_{\alpha\hspace{0,05cm} \cup\hspace{0,03cm} \beta-\alpha}$$
Now using the observation above we obtain an equation of the form 
$$m_{\beta}+m_{A}+\ell=m_{\alpha'}+m_{B}$$
which is a contradiction if we substitute further $m_{\alpha'}=m_{\beta}+m_{\alpha'-\beta}-\ell$ and $m_{\alpha'}=m_{\beta}-m_{\beta-\alpha'}$ respectively depending  whether $\beta-\alpha'$ is positive or negative.\par
If $\alpha\succeq \beta$ is induced from Definition~\ref{poset1}(i), then a simple case consideration shows that
$\alpha-\beta\in R^+$ if and only if $\alpha'-\beta\in R^+$ and moreover
$\alpha\hspace{0,05cm} \cup\hspace{0,03cm} \alpha'-\alpha'=\alpha-\beta$ if $\alpha-\beta\in R^+$ and $\alpha\hspace{0,05cm} \cup\hspace{0,03cm} \alpha'-\alpha'=\beta-\alpha'$ otherwise. 
Hence \eqref{ann12}
leads to a contradiction. \par So we can suppose that $\alpha\succeq \beta$ is induced from Definition~\ref{poset1}(ii) and we have either of the following situations
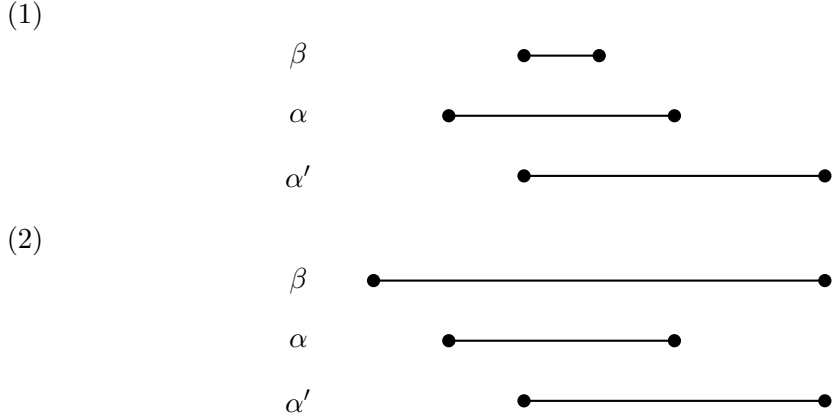
\begin{enumerate}
\item $$\begin{tikzpicture}[thick]
\draw (1,0) -- (2, 0);
\node (A) at (-2,0) [] {$\beta$};
\node (R) at (2,0) [] {};
\node (P) at (1,0) [] {};
\node (E) at (1,-1.6) [] {};
\node (G) at (5, -1.6) [] {};
\fill (R) circle (2.5pt);
\fill (P) circle (2.5pt);
\draw (0,-0.8) -- (3, -0.8);
\node (C) at (-2,-0.8) [circle] {$\alpha$};
\node (D) at (-2,-1.6) [circle] {$\alpha'$};
\node (RE) at (0,-0.8) [] {};
\node (PE) at (3,-0.8) [] {};
\fill (RE) circle (2.5pt);
\fill (PE) circle (2.5pt);
\fill (E) circle (2.5pt);
\fill (G) circle (2.5pt);
\draw (1,-1.6) -- (5, -1.6);
\end{tikzpicture}
$$
\item $$\begin{tikzpicture}[thick]
\draw (-1,0) -- (5, 0);
\node (A) at (-2,0) [] {$\beta$};
\node (R) at (-1,0) [] {};
\node (P) at (5,0) [] {};
\node (E) at (1,-1.6) [] {};
\node (G) at (5, -1.6) [] {};
\fill (R) circle (2.5pt);
\fill (P) circle (2.5pt);
\draw (0,-0.8) -- (3, -0.8);
\node (C) at (-2,-0.8) [circle] {$\alpha$};
\node (D) at (-2,-1.6) [circle] {$\alpha'$};
\node (RE) at (0,-0.8) [] {};
\node (PE) at (3,-0.8) [] {};
\fill (RE) circle (2.5pt);
\fill (PE) circle (2.5pt);
\fill (E) circle (2.5pt);
\fill (G) circle (2.5pt);
\draw (1,-1.6) -- (5, -1.6);
\end{tikzpicture}
$$
\end{enumerate}
In the first case we get
$$m_{\alpha}=m_{\beta}+m_{\gamma_1}+m_{\gamma_2},\ \ m_{\alpha'}=m_{\beta}+m_{\alpha'-\beta}-\ell.$$
Substituting into \eqref{ann12} gives
$$m_{\gamma_1}+m_{\gamma_2}+m_{\alpha\hspace{0,05cm} \cup\hspace{0,03cm} \alpha'-\alpha}=m_{\alpha'-\beta}-\ell+m_{\alpha\hspace{0,03cm}\cup \hspace{0,03cm} \alpha'-\alpha'}$$
which is a contradiction since $\gamma_1=\alpha\hspace{0,03cm}\cup \hspace{0,03cm} \alpha'-\alpha'$. In the second case 
$$m_{\beta}=m_{\alpha}+m_{\gamma_1}+m_{\gamma_2}-2\ell=m_{\alpha'}+m_{\beta-\alpha'}.$$
Substituting into  \eqref{ann12} gives a contradiction
$$2\ell+m_{\beta-\alpha'}+m_{\alpha\hspace{0,05cm} \cup\hspace{0,03cm} \alpha'-\alpha}=m_{\alpha\hspace{0,03cm}\cup \hspace{0,03cm} \alpha'-\alpha'}+m_{\gamma_1}+m_{\gamma_2}$$
since $\gamma_2=\alpha\hspace{0,03cm}\cup \hspace{0,03cm} \alpha'-\alpha$.

\bibliographystyle{plain}
\bibliography{bibfile}
\end{document}